\documentclass[english,11pt]{article}
\usepackage[german,french,english]{babel}
\usepackage[cp850]{inputenc}
\usepackage{latexsym,graphicx, fancybox}
\usepackage{graphicx}
\usepackage{booktabs}

\usepackage{amssymb, amsmath, amsfonts}
\usepackage{amsmath, amssymb}
\usepackage{color}
\usepackage{tikz}
\usepackage{latexsym}

\usepackage[colorlinks=true, allcolors=blue]{hyperref}    
\usepackage{tocloft}        

\usepackage{mathrsfs}

\usepackage{float}   
\usepackage{graphicx}  
\usepackage{hyperref}  
\usepackage{latexsym}  
\usepackage{xcolor}   
\usepackage{lipsum}   
\usepackage{subcaption}

\usepackage{amssymb}  
\usepackage{amsmath}  
\usepackage{amsbsy}
\usepackage{amsthm}
\usepackage{bbm}    
\usepackage{eucal}   
\usepackage{mathrsfs}  

\usepackage{multirow}  
\usepackage{ulem}    

\usepackage[numbers,sort&compress]{natbib} 

\usepackage{etoolbox}
\apptocmd{\thebibliography}{\setlength{\itemsep}{0.01pt}}{}{}

\usepackage{microtype}
	\usepackage{amsthm}
	
	\newcommand{\dd}{\,\mathrm{d}}
	\newcommand{\R}{\mathbb{R}}
	\newcommand{\N}{\mathbb{N}}
	\newcommand{\E}{\mathbb{E}}
	\renewcommand{\P}{\mathbb{P}}

	\usepackage{authblk}
	
	\newtheorem{Theorem}{Theorem}[section]
	\newtheorem{Definition}[Theorem]{Definition}
	\newtheorem{Proposition}[Theorem]{Proposition}
	\newtheorem{Assumption}[Theorem]{Assumption}
	
	\newtheorem{Lemma}[Theorem]{Lemma}
	\newtheorem{Corollary}[Theorem]{Corollary}
	
	\newtheorem{Example}[Theorem]{Example}
	
	\newtheorem{assumption}{Assumption}[section]

	\usepackage{thm-restate}
	\usepackage{float}
	\usepackage{hyperref}
	
\begin{document}

\selectlanguage{english}

\title{\bf 
	On Inhomogeneous Affine Volterra Processes: Stationarity and Applications to the Volterra Heston Model\footnote{Emmanuel Gnabeyeu and Mathieu Rosenbaum gratefully acknowledge support from the ILB Chair \textit{Artificial Intelligence and Quantitative Methods for Finance} at University Paris Dauphine-PSL.} 
}

\author{
	Emmanuel Gnabeyeu\\
	LPSM, Sorbonne Universit\'e and Universit\'e Paris Cit\'e, Paris, France.\\
	E-mail: {\tt emmanuel.gnabeyeu\_mbiada@sorbonne-universite.fr}
	\and
	Gilles Pag\`es\\
	LPSM,  Sorbonne Universit\'e, case 158, 4, pl. Jussieu, F-75252 Paris Cedex 5, France. \\
	E-mail: {\tt gilles.pages@sorbonne-universite.fr}
	\and
	Mathieu Rosenbaum\\
	CEREMADE, \textit{Universit\'e Paris Dauphine-PSL}, Paris, France. \\
	E-mail: {\tt mathieu.rosenbaum@dauphine.psl.eu}
} 
\date{November 1, 2025}

\maketitle
\vspace{-.90cm} 
\renewcommand{\abstractname}{Abstract}
\begin{abstract}
	In this paper, we 
	study the finite-time behavior of solutions to scaled inhomogeneous affine Stochastic Volterra equations (SVIEs) through the lens of
	a weaker notion of stationarity, referred to as a \textit{fake stationary regime}, in which 
	all marginal
	distributions share the same expectation and variance.
	As a first application, we introduce the \textit{fake stationary Volterra Heston model} and derive a closed-form expression for its characteristic function.
	Having established this finite-time proxy for stationarity, we then investigate the asymptotic (long-time) behavior to assess whether genuine stationary regimes emerge in the limit.
	Using an extension of the exponential-affine transformation formula for these processes, we establish the existence in the long run of limiting distributions, which, unlike in the classical affine diffusion setting, may depend on the initial random state; we further show that this dependence vanishes when the constant mean imposed initially satisfies a condition involving the resolvent of the kernel and the mean-reversion parameters.
	We then proceed to the construction of stationary processes associated with these limiting distributions. However, the dynamics in this long-term regime are analytically intractable, and 
	the original process need not be stationary in the classical sense on finite horizons. 
	This highlights the relevance of finite-time analysis through the lens of the aforementioned 
	 \textit{fake stationarity}, which offers a tractable approximation to stationary behavior in genuinely non-stationary Volterra systems.	  
\end{abstract}



\vspace{-.3cm}
\section{Introduction}
\vspace{-.2cm}
\subsection{Literature review}
We consider the class of convolution kernels i.e. kernels $K: \{(s,t)\!\in \mathbb{R}_+^2: 0\le s < t \} \to \mathbb{R}_+$ satisfying  \(\forall\, s,\, t\ge 0, \; s<t, \; K(s,t)= K(0,t-s)\) and we are chiefly interested in \textit{inhomogeneous affine Volterra stochastic integral equations} which naturally arise in the modeling of systems with memory, (including) in mathematical finance, physics, and biology:
\vspace{-.2cm}
\begin{equation}\label{eq:intro}
	\begin{cases}
		X_t = X_0 \varphi(t) + \int_0^t K(s,t) b(s, X_s) \, ds + \int_0^t K(s,t) \sigma(s, X_s) \, dW_s, \\
		X_0 : (\Omega, \mathcal{F}, \mathbb{P}) \to (\mathbb{R}, \mathcal{B}(\mathbb{R})), \quad X_0\perp\!\!\!\perp W.
	\end{cases}
	\vspace{-.1cm}
\end{equation}
Here, $\varphi$ is a deterministic continuous function, $K(s,t)$ a deterministic kernel modeling the memory or hereditary structure of the system and the process $(W_t)_{t \geq 0}$, an $\mathbb{R}$-valued Brownian motion independent of $X_0$, both defined on a probability space $(\Omega, \mathcal{A}, \P)$ with $\mathcal{F}_t \supset \mathcal{F}_{t}^{X_0, W}$ a filtration satisfying the usual conditions.
 The drift \(b\) and the diffusion coefficient \(\sigma\)
are time-inhomogeneous and affine in the state variable in the sense that the functions $b$ and $a := \sigma \sigma^\top$ belong to $ \mathrm{Pol}_1(\mathbb{R})$, where $\mathrm{Pol}_n(\mathbb{R})$ denote the subspace of the ring $\mathrm{Pol}(\mathbb{R})$ of all polynomials on $\mathbb{R}$, consisting of polynomials of degree at most $n\in\N$.

 Special cases of affine processes have long been instrumental in various domains of stochastic analysis, particularly in mathematical finance and biological modeling. Seminal contributions include the foundational works of~\cite{Feller1951}, ~\cite{KawazuWatanabe1971}, and ~\cite{CoxIngersollRoss1985}. The theoretical framework of affine processes was further advanced through significant milestones achieved by ~\cite{duffie2003affine}, as well as ~\cite{filipovic2009affine}.
More recently, the study of affine processes within the context of Volterra equations has attracted considerable attention~\cite{JaissonR2016,ElEuchR2018,abi2019affine,ackermann2022}. 
In this setting, the dynamics depend not only on the current time $s$ but also on the terminal time $t$ via the kernel $K$ appearing in equation~\eqref{eq:intro}. This dependence introduces memory effects, thereby violating the Markov property. The aim is to mimic the behavior of affine stochastic differential equations (SDEs) driven by fractional Brownian motion (\(W^H_\cdot\)), while avoiding the technical challenges associated with their analysis, particularly those arising from \textit{``high-order" rough path theory} or \textit{regularity structures}, which become especially intricate as the Hurst parameter \( H \)  playing a crucial role in the path's regularity, approaches zero. This type of dynamic has attracted significant attention in the financial community over the past decade, particularly because it enables consistent modeling of financial markets across multiple time scales from order book dynamics to the pricing of derivative products despite its non-Markovian nature.
More recently, the roughness phenomenon in the behavior of high-frequency volatility data has come under the spotlight. Empirical studies suggest that the log-volatility is well modeled by a fractional Brownian motion with a Hurst parameter less than \( 0.5\); see, for instance, 
\cite{GatheralJR2018} where the fractional versions of affine models
are able to reproduce the slope of short-term skew observed in option markets.

These empirical findings have motivated the development of rough extensions of the celebrated Heston stochastic volatility model~\cite{Heston1993}, whose widespread practical success has made it a natural starting point for such extensions. However, in the classical Heston~\cite{Heston1993} model equipped with its Cox-Ingersoll-Ross (CIR) driven stochastic
volatility process, it has long been observed by practitioners that
when the variance process is initialized by a deterministic value at the origin,
usually its long-run mean, the model often leads to poor calibration since the volatility dynamic goes through a transient regime for short maturities, then it deploys until attaining its stationary regime (a $\gamma-$distribution ).
A natural remedy to get rid of this transient effect is to consider directly a Heston model in which the volatility is taken under its stationary
regime (see e.g.~\cite{PagesPanloup2009}), which not only improves calibration across maturities (see e.g.~\cite{LemMonPag2022}), but also reduces by one the number of free parameters. A similar limitation is expected to arise in rough volatility models driven by rough Volterra versions of the CIR dynamics.
Consequently, the search for some stationarity properties of such rough CIR processes~\eqref{eq:intro} can
be of interest to make possible a calibration consistent for both short and long maturities, and even constitute fundamental problems in stochastic analysis. This question has recently attracted considerable attention for Volterra CIR-type processes: Existing ergodicity results rely mainly on Markovian lifting methods (see, e.g.,~\cite{Jacquieretal2022})
 or on the analysis of Volterra-type Riccati equations (see, e.g.,~\cite{friesen2022volterra}).
 These works establish the existence of stationary solutions, without providing an explicit characterization of their dynamics. More recently, \cite{Pages2024,EGnabeyeu2025} introduced the notion of finite-time \textit{fake
 stationary regime} for SVIEs with general diffusion coefficients and proved that, under this regime, the limiting process is weakly $L^2$-stationary.
 Our aim is to investigate this property in the specific setting of~\eqref{eq:intro} over a finite horizon and in the long run. This regime is  characterized by a functional equation satisfied by the coefficients of the SVIEs~\eqref{eq:intro}, which adjusts them to impose the desired mean and variance structure.
This framework aims to overcome the transient effects inherent in classical affine processes, thereby providing a unified and consistent description of their dynamics.

\vspace{-.3cm}
\subsection{Our contribution}
We investigate the finite-horizon \textit{fake
	stationary regime} of affine SVIEs~\eqref{eq:intro}, introduce a (fake) stationary version of the Volterra Heston model of~\cite{abi2019affine}, and prove that as $u\to+\infty$, the time-shifted process \((X_{t+u})_{t \geq 0}\)
functionally weakly converges to a limiting continuous process, which unlike in~\cite{Pages2024,EGnabeyeu2025, GnabeyeuPages2026b} is strongly stationary (i.e. in the classical sense, where its finite-dimensional distribution is invariant under time shifts). We first establish a \textit{measure-extended} affine transformation formula for~\eqref{eq:intro} associated to a  measure-forced inhomogeneous Riccati--Volterra equation and allowing for explicit control over the finite-dimensional distributions of the process.
We then prove the existence of limiting distributions and construct the associated stationary processes, explicitly deriving the characteristic functions of their finite-dimensional distributions; all of which may retain a dependence on the initial random state, and we identify a condition involving the resolvent of the kernel and the mean-reversion parameters under which this dependence vanishes. However, the pathwise dynamics
of the limiting stationary processes $X^\infty$ (an SVIE for $X^\infty$) remain open. 
This motivates the study of \textit{fake stationary regimes}, which provide a tractable framework for investigating both finite- and long-time behaviors in settings where classical stationarity is either unavailable~\cite{Pages2024} or analytically intractable~\cite{friesen2022volterra,Jacquieretal2022}.
As a main application, we introduce the \textit{fake stationary Volterra Heston model} and derive a semi-closed-form representation of its characteristic function in terms of a fractional Riccati equation, making the model  highly
tractable and 
amenable to efficient option pricing and volatility  fitting via \textit{Fourier inversion methods}. 
\vspace{-.3cm}
\subsection{Structure of the paper and Notations} 
\noindent {\sc \textbf{Organization of the Work.}}  
This work is organized as follows.
In Section~\ref{sect-tools}, we present the main assumptions and collect some preliminaries needed for the study of affine SVIEs~\eqref{eq:intro}, including the so-called resolvent of the convolution kernel and the forward process.
In Section~\ref{sect-CharacteristicFunction} we establish that the conditional measure-extended Laplace functional of the 
SVIEs~\eqref{eq:intro} can
be represented via an exponential-affine form by a solution of an inhomogeneous measure-extended Riccati-Volterra integral equation. We further analyze this Riccati-Volterra equation and establish existence, local $L^p$-integrability and regularity of its solution.
Section~\ref{sect-FakeStationarity} is devoted to establishing \textit{``fake stationarity regime"} for affine SVIEs~\eqref{eq:intro}
 in finite time; we show that, under appropriate functional equations satisfied by its coefficients, all marginal distributions share the same mean and variance. We then introduce the \textit{fake stationary Volterra Heston model}, a natural extension of the rough Heston model introduced by~\cite{ElEuchR2018, el2019characteristic} and its Volterra version in~\cite{abi2019affine}, and widely used in mathematical finance. We also provide its characteristic function, namely its Fourier-Laplace exponential-affine representation.
Building on the Fourier-Laplace transform representation developed 
in Section~\ref{sect-CharacteristicFunction}, Section~\ref{sect-longRun} establishes the weak convergence of the conditional law of $X_t$ given $X_0$ towards a limiting random distribution $\pi_{X_0}$ as $t \to \infty$. We also construct
the stationary process associated with the limiting conditional law $\pi_{X_0}$
 and characterize when the limit $\pi_{X_0}$ actually depends on the initial random state or distribution $X_0$. We then analyse the functional long-run asymptotics in the fake stationarity regime.
Finally, in Section~\ref{sect-Application}, we provide a broad class of applications, including a numerical display of the finite-time \textit{fake stationary fractional CIR} process and its long-run confluence, followed by the SPX volatility surface calibration in the \textit{fake stationary rough Heston model}. 


\medskip
\noindent {\sc \textbf{Notations.}} 

\smallskip
\noindent $\bullet$ Denote for any $T>0$, $\mathbb{T} = [0, T] \subset \mathbb{R}_+$. 
$\mathrm{Leb}_{\mathbb{R}}$ (resp.\ $\mathrm{Leb}_{\mathbb{R}_+}$)
denotes the Lebesgue measure on $(\mathbb{R}, {\cal B}or(\R))$
(resp.\ $(\mathbb{R}_+=[0,\infty),{\cal B}or(\R_+))$).
 Let $\mathbb{H} \in \{\mathbb{R},\mathbb{R}_+,\R^d, \mathbb{C}\}$ etc. 

\noindent $\bullet$ $\mathbb{X} := {\cal C}([0,T], \mathbb H)  (\text{resp.} \; {\mathcal C_0}([0,T], \mathbb H))$ denotes the set of continuous functions(resp. null at 0)  from $[0,T]$ to $\mathbb H $ and ${\cal B}or({\cal C}_d)$ denotes the  Borel $\sigma$-field of ${\cal C}_d$ induces by the $\sup$-norm topology.

\smallskip 
\noindent $\bullet$ For $p\in(0,+\infty)$, $L_{\mathbb H}^p(\P)$ or simply $L^p(\P)$ denote the set of  $\mathbb H$-valued random vectors $X$  defined on a probability space $(\Omega, {\cal A}, \P)$ such that $\|X\|_p:=(\E[\|X\|_{\mathbb H}^p])^{1/p}<+\infty$. 

\smallskip 
\noindent $\bullet$ Let $\mathcal{P}(\mathbb H)$ denote the space of Borel
probability measures on $\mathbb H$.

\smallskip 
\noindent $\bullet$
We let $\mathbb M_{\mathrm{lf}}(\R_+, \mathbb{H})$ denote the space of $\mathbb{H}$-valued locally finite measures, so that $|m|([0,T])<\infty$
for every $T<\infty$ i.e. the space of all $\mathbb{H}$-valued $(\R_+, {\cal B}or(\R))$-measurable functions \(m\) on \(\mathbb{R}_+\) such that the restriction \(m|_{[0, T]}\), for any \(T > 0\), is a \(\mathbb{H}\)-valued finite measure (i.e. the restriction $m|_{[0,T]}$ with $T > 0$ is well-defined). For \(m \in \mathbb M_{\mathrm{lf}}(\R_+, \mathbb{H})\) and a compact set \(E \subset \mathbb{R}_+\), we define the total variation of \(m\) on \(E\) by:
\centerline{$|m|(E) := \sup \left\{ \sum_{j=1}^N |m(E_j)| : \{E_j\}_{j=1}^N \text{ is a finite measurable partition of } E \right\}.$}

\noindent Now we define \(\mathbb{H}_{-} := \{ u \in \mathbb{H} : \operatorname{Re}(u) \leq 0 \}\) and we denote by \( M_{\mathrm{lf}}^{-}(\R_+, \mathbb{H}) \subset M_{\mathrm{lf}}(\R_+, \mathbb{H}) \) the subset of \(\mathbb{H}\)-valued locally finite measures \(m \in M_{\mathrm{lf}}(\R_+, \mathbb{H})\) such that \(\operatorname{Re}(m) \leq 0\).
 
\smallskip 
\noindent $\bullet$ For $f, g \!\in {\cal L}_{\R_+,loc}^1 (\R_+, {\rm Leb}_1)$, we define their convolution by $(f*g)(t) = (f\stackrel{dt}{*}g)_t= \int_0^tf(t-s)g(s)ds$, $t\ge 0$.  We will frequently use Young's convolution inequality
which states that for any \(T>0\), \( f \in L^p([0,T]) \), \( g \in L^q([0,T]) \), and \( 1 \leq p, r, q \leq \infty \) such that \(\frac{1}{p} + \frac{1}{q} = 1 + \frac{1}{r} \), the convolution \( f * g \) belongs to \( L^r([0,T]) \), and Young's inequality writes: \(\quad \| f * g \|_{L^r([0,T])} \leq \| f \|_{L^p([0,T])} \cdot \| g \|_{L^q([0,T])}.\)

\smallskip 
\noindent $\bullet$ Convolution between a function and a measure. Let \(f : (0, T] \to \mathbb{R}\) be a function and \(m \in \mathbb M_{\mathrm{lf}}(\R_+, \mathbb{H})\). Their convolution (whenever the integral is well-defined) is defined by
\vspace{-.2cm}
\begin{equation}\label{eq:convolmeasure}
	(f * m)(t)= \int_{[0,t)} f(t - s) \, dm(s) = \int_{[0,t)} f(t - s) \, m(ds) = (f\stackrel{m}{*}\mathbf{1})_t, \quad t \in (0, T].
	\vspace{-.3cm}
\end{equation}
\vspace{-.2cm}
\smallskip 
\noindent $\bullet$ For a random variable/vector/process $X$, we denote by $\mathcal L(X)$ or $[X]$ its law or distribution. 

\smallskip 
\noindent $\bullet$ $X\perp \! \! \!\perp Y$  stands for independence of random variables, vectors or processes $X$ and $Y$.  

\noindent $\bullet$ For a measurable function \( \varphi: \mathbb{R}^+ \to \mathbb{R} \), we denote: 

\smallskip
\centerline{$
	\forall p \geq 1, \quad \| \varphi \|^p_{L^p([0,T])} := \int_0^{T} |\varphi(u)|^p \, du,  \; \|\varphi\|_{\infty,T}:= \sup_{u\in [0,T]}|\varphi(u)|\; \text{and} \; \displaystyle \|\varphi\|_{\infty}=\|\varphi\|_{\sup} := \sup_{u\in \mathbb{R}^+}|\varphi(u)|.
	$}

\smallskip 
\noindent $\bullet$ $\Gamma(a) = \int_0^{+\infty} u^{a-1} e^{-u} \, du, \quad a > 0, \; 
\text{and} \; 
B(a, b) = \int_0^1 u^{a-1} (1 - u)^{b-1} \, du, \quad a, b > 0.$
We recall the classical identities:
$\Gamma(a + 1) = a \, \Gamma(a) 
\quad \text{and} \quad 
B(a, b) = \frac{\Gamma(a)\Gamma(b)}{\Gamma(a + b)}$
and we set \(\R_+=[0,+\infty)\), \(\R_-=(-\infty,0]\).

The results of this paper are developed in a one-dimensional setting to keep the presentation fairly short without cumbersome notation. 

\vspace{-.2cm}
\section{A primer on Volterra processes with convolution kernels and tools.}\label{sect-tools}
\vspace{-.2cm}
From now on, to alleviate notations, we denote $K(t) := K(0, t)$ so that $K(s, t) = K(t - s)$ and  we restrict our attention to
  the special case of a \textit{scaled} inhomogeneous stochastic Volterra equation 
  associated to the convolution kernel $K:\R_+\to \R_+$ satisfying~\eqref{eq:contKtilde} and~\eqref{eq:Kcont} of assumption~\ref{assump:kernelVolterra} further on:
  \vspace{-.2cm}
\begin{equation}\label{eq:Volterrameanrevert}
	X_t = X_0\varphi(t) +\int_0^t K(t-s)(\mu(s)-\lambda X_s)ds + \int_0^t K(t-s)\varsigma(s)\sigma(X_{ s})dW_s, \quad X_0\perp\!\!\!\perp W.
\end{equation}
where $\lambda>0$, $\varphi : \mathbb{T} \to \R_+ $, $\mu :\mathbb{T}\to \R$ are Borel functions, \(\varsigma\) a deterministic Borel function and $\sigma: \mathbb{R} \to  \R$ is \(\gamma-\) H\"older continuous in $x$ (\(\gamma \in [\tfrac{1}{2}, 1]\)), locally uniformly in $t\!\in \mathbb{T}$.
We always work under the assumption below, which applies to the inhomogeneous Volterra equation~\eqref{eq:Volterrameanrevert}.
\begin{assumption}[On Volterra Equations with convolution kernels]\label{assump:kernelVolterra}
	\noindent
	\begin{enumerate}
		\item[(i)] Assume that the kernel $K$ 
		satisfy for every \(T>0\) :
		\begin{itemize}
			\item The integrability assumption: The following is satisfied for some $\widehat\theta\in (0,1]$.
			\begin{equation}\label{eq:contKtilde}
				(\widehat {\cal K}^{cont}_{\widehat \theta})\;\;\exists\,\widehat\kappa< +\infty,\;\forall\bar{\delta}\!\in (0,T],\; \widehat \eta(\delta) :=  \sup_{t\in [0,T]} \left[\int_{(t-\bar{\delta})^+}^t \hskip-0,25cm K\big(t-u\big)^2 du\right]^{\frac12}\le \widehat \kappa \,\bar{\delta}^{\,\widehat \theta}.
			\end{equation}
			\item  The continuity assumption: \(({\cal K}^{cont}_{\theta}) \;\;  
			\exists\, \kappa< +\infty,\; \exists \; \theta\in (0,1] \; \text{such that}\; \forall \,\bar{\delta}{\in (0, T)}\)
			\begin{equation}\label{eq:Kcont}
				({\cal K}^{cont}_{T,\theta}) \; \forall \,\bar{\delta}{\,\in (0, T)},\; \eta(\bar{\delta}):= \sup_{t\in [0,T]} \left[\int_0^t |K(\big(s+\delta)\wedge T\big)-K(s)|^2ds \right]^{\frac 12} \le  \kappa\,\bar{\delta}^{\theta}.
			\end{equation}
		\end{itemize}  
		\item[(ii)] Assume that the state dependent diffusion coefficient $\mathbb{R}\ni x\mapsto \sigma(x)$ is of linear growth and H\"older continuous of order \(\gamma\) for some \(\gamma \in [\tfrac{1}{2}, 1]\).
		Hence, there are constants \(C^\prime_{\sigma}, C_{\sigma} > 0\) such that 
		\begin{equation}\label{eq:LLcoef}  |\sigma(x)| \leq C^\prime_{\sigma}(1 + |x|),
		\quad \text{and }\quad 
		|\sigma(x) - \sigma(y)| \leq C_{\sigma}|x - y|^{\gamma}
		\; \text{hold for all } x, y \in \mathbb{R}.
		\end{equation}
		
		\item[(iii)] Finally, assume that
		the input curve $t \to x_0(t) :=X_0 \varphi(t)$ is absolutely continuous, $(\mathcal F_t)$-adapted and for some 
		\( p \in (0, +\infty) \) such that \( X_0 \in L^p(\mathbb{P}) \), there exists $\delta > 0$ such that for any $T > 0$,  
		\begin{equation}\label{eq:Input} 
		\mathbb{E} \,\!\Big(\sup_{t \in [0,T]} |x_0(t)|^p\Big) < +\infty,\quad 
		\mathbb{E}\!\big[\,|x_0(t') - x_0(t)|^p\,\big] 
		\le C_{T,p} \Big( 1 + \mathbb{E}\,\big[\sup_{t \in [0,T]} |x_0(t)|^p\big] \Big) |t' - t|^{\delta p}.
	   \end{equation}
	\end{enumerate}
\end{assumption}
\noindent Notice that the drift term \(b(t,x)=\mu(t)-\lambda x\) in~\eqref{eq:Volterrameanrevert} is clearly globally Lipschitz continuous in $x\in \R $, uniformly in time
 and has 
 linear growth  in the state variable $x$, though not necessarily uniformly in time. Since \(\varsigma\) is locally bounded, the same $\gamma$-H\"older continuity and linear-growth property holds for the full diffusion coefficient \((t,x)\mapsto\sigma(t,x)=\varsigma(t)\sigma(x)\) locally uniformly
in time.

\smallskip
\noindent Under Assumption~\ref{assump:kernelVolterra}, if \(\mu\) is bounded, standard existence results for stochastic Volterra equations yield a continuous solution \( (X_t)_{t\geq0} \) to~\eqref{eq:Volterrameanrevert}; when \(\gamma=1\), pathwise uniqueness yields a unique strong solution  adapted to \( (\mathcal{F}^{X_0, W}_t)_{t\geq0} \) (see, e.g., \cite{ZhangXi2010, GnabeyeuPages2026b} ), whereas for \(\gamma\in[1/2,1)\) (see, e.g., \cite{PromelScheffels2023,PromelScheffels2023b} ), only weak existence is required in what follows. 
These results are applied successively on finite time intervals \([0,T]\), \(T\in\mathbb N\), and the resulting solutions are glued together to obtain a solution on \(\mathbb R_+\).
In particular, for every \(T>0\) and every admissible \(p>0\),  a continuous solution on $\R_+$ to  Equation~\eqref{eq:Volterrameanrevert}, starting from   $X_0$ satisfies (among other properties),
{\small   
	\begin{equation}\label{eq:L^p-supBound}
		\forall\, T>0, \; \exists \,C_{_{T,p} }>0,\quad \mathbb E\!\Big[\sup_{t\in[0,T]}|X_t|^p\Big] \le C_{_{T,p} } \Big( 1 + \sup_{t \in [0,T]} |\varphi(t)|^p \mathbb{E}\big[|X_0|^p\big] \Big)<\infty.
	\end{equation}
}
Consequently, the linear growth of the state-dependent diffusion coefficient in Assumption~\ref{assump:kernelVolterra}(ii) also gives $\mathbb{E}[ \sup_{t \in [0,T]}  |\sigma(t, X_{ t})|^p]=\mathbb E\,\!\Big[\sup_{t\in[0,T]}
|\varsigma(t)\sigma(X_t)|^p\Big] < C_{T,\sigma} (1 +  \sup_{t \in [0,T]} |\varphi(t)|^p \mathbb{E}[|X_0|^p]) <\infty$  for some constant $C_{T,\sigma}>0$,
which provides the integrability required below for applications of both deterministic (regular) and stochastic Fubini theorems.

\medskip\noindent{\bf Remark:}
	We will mainly focus on the affine case in sections~\ref{sect-CharacteristicFunction},~\ref{sect-FakeStationarity},~\ref{subsect:LimitD}-\ref{subsect:longRun} and~\ref{sect-Application}, for which \(\gamma=\frac12\), 
	and assume that the solution remains nonnegative on \(\R_+\). This is in particular the case, when \(\sigma(x):=\nu\sqrt{x}\) (see sections~\ref{subsect:FS-RH} and~\ref{sect-Application}) as shown in~\cite{EGnabeyeuR2025}  for the inhomogeneous \(\alpha-\) fractional square root process (corresponding to \(K=K_{1,\alpha,0}\) in Example~\ref{Ex:SolventGammaKernel} further on ) building upon scaling limit theorems.
	There, positivity 
	is established by identifying the process as the joint time-space scaling limit of a sequence of intensity processes associated with nearly unstable Hawkes processes. Importantly, the mean-reversion function $\mu$
	may be unbounded near zero
	(see,
	Example~\ref{Ex:FakeStatAutonome_Lm}, Part~2(a) and~\cite[Example 2.7-$(2)$]{EGnabeyeuR2025}).

\smallskip\noindent
We next introduce key tools including (functional) resolvent of a Kernel and Wiener--Hopf equation that will be important to our analysis.

\subsection{Resolvent of convolution Kernels}
Throughout this work, we will adopt
 the functional resolvent introduced in \cite{gripenberg1990} and also discussed or presented in works such as \cite{abi2019affine}.
Let $K: (0,\infty) \to \R_+$ be a convolution kernel satisfying ~\eqref{eq:contKtilde},~\eqref{eq:Kcont}.     
For every $\lambda>0$,  the \textit{ resolvent or Solvent core} $f_{\lambda}$ associated to $\lambda K$, known as the \textit{resolvent of the second kind of $\lambda K$} is defined as the unique solution to the deterministic Volterra equation
\vspace{-.1cm}
\begin{equation}\label{eq:flambda-eq}
	\forall\,  t\ge 0,\quad f_{\lambda}(t) + \lambda \int_0^t K(t-s)f_{\lambda}(s)ds = \lambda K(t).
\end{equation}
or, equivalently, written in terms of convolution, 
$f_{\lambda} +\lambda K *f_{\lambda}=\lambda   K.$ Since $K$ is locally integrable, the resolvent $f_{\lambda}$ always exists.
This equation is also known as a \textit{renewal equation}. Its solution  always satisfies $f_{\lambda}(0^+)=\lambda K (0^+)$ and admits the formal \textit{Neumann series expansion} 
\begin{equation}\label{eq:Resolvent3}
	f_{\lambda} = \sum_{n\ge 1} (-1)^{n-1} \lambda^n K^{*n}.
\end{equation}
where \(K^{*n}\) denotes the \(n\)-th convolution of \(K\) or the $n$-fold $*$
product of \(K\) defined by $ K^{*n}(t) = \int_0^t K(t - s) \cdot K^{*(n-1)}(s) \, ds$, with the convention, $K^{*0}= \delta_0$ (Dirac mass at $0$) so that $ K^{*1} = K $.

\smallskip \noindent{\bf Remarks} 1. If the (non-negative) kernel $K$  satisfies 
\begin{equation}\label{eq:Kcontrol}
	0\le K(t)\le Ce^{bt }t^{a-1}\mbox{  for some }\; a>0,\, \; b,\; C\geq0\!\in \R_+. 
\end{equation} then, by induction  \(K^{*n}(t) \le C^n e^{bt}\frac{\Gamma(a)^n}{\Gamma(na)} t^{na-1},\) so that for such kernels, using Stirling's formula, one readily
checks that the series in~\eqref{eq:Resolvent3} converges normally (i.e. locally uniformly and absolutely) on
$(0,\infty)$. In particular,
the resolvent $f_\lambda$ is (well-defined and) continuous on $(0,\infty)$.

2. If the kernel $K$ has a finite Laplace transform i.e.  \(\forall\, t>0,\; 0<L_K(t)<+\infty\), then,
taking the Laplace transform of both sides of the above equality~\eqref{eq:flambda-eq}, we have that :
$L_{f_{\lambda}}(t)(1+\lambda L_K(t))=\lambda L_K(t) $, $t>0$.
This yields:
\vspace{-.1cm}
\begin{equation}\label{eq:Lapl_Resolvent}
	L_{f_{\lambda}}(t) = \frac{\lambda L_K(t)}{1+\lambda L_{K}(t)}.
\end{equation}
which is $(0, 1)$-valued.
Consequently, for $\lambda \geq 0, $ $L_{f_{\lambda}}(t) \equiv 0$ if and only if 
$L_K(t) \equiv 0$ i.e. if and only if $K=0$ by the injectivity of Laplace transform.

\begin{Example}[Resolvent associated to the $\rho-$exponential $\alpha-$fractional Kernel]\label{Ex:SolventGammaKernel}
	Let the Kernel $K$ be the $\rho$-exponential $\alpha$-fractional integration kernel or Gamma kernel given by $K(t) = K_{b, \alpha , \rho }(t):= b e^{-\rho t}\frac{ t^{\alpha - 1}}{\Gamma(\alpha)}  \cdot \mathbf{1}_{(0,\infty)}(t)
	$, for \(b > 0, \alpha \in(\frac12,1) \) and \(\rho \geq 0\). One checks that these kernels also satisfy~\eqref{eq:contKtilde} and~\eqref{eq:Kcont} 
	(with $\widehat{\theta}=\theta= (\alpha-\frac 12)\wedge 1$) and trivially~\eqref{eq:Kcontrol}.
	It follows from the easy identity \(K_{b,\alpha, \rho} * K_{b,\alpha^\prime, \rho} = b K_{b,\alpha + \alpha^\prime, \rho}\) that $K_{b,\alpha,\rho}^{*k}(t)
	=
	b^k e^{-\rho t}
	\frac{t^{k\alpha-1}}{\Gamma(k\alpha)}$ for every $t>0$ and the Neumann series expansion provided in equation~\eqref{eq:Resolvent3} implies that the resolvent
	reads:
	{\small
		\begin{equation}\label{eq:DerivSolventGammaKernel}
			f_{b, \alpha, \rho, \lambda}(t) = \sum_{k\ge 1} (-1)^{k-1} \lambda^k K_{b, \alpha, \rho}^{(*k)} =  -\sum_{k\ge1}
			(-1)^k\lambda^k
			b^k e^{-\rho t}
			\frac{t^{k\alpha-1}}{\Gamma(k\alpha)}= b\lambda e^{-\rho t}
			t^{\alpha-1}
			E_{\alpha,\alpha}(-b\lambda t^\alpha).
		\end{equation}
	}
	where 
	$E_{\alpha,\beta}(z)
	=
	\sum_{k=0}^\infty
	\frac{z^k}{\Gamma(\alpha k+\beta)},$ the Mittag--Leffler function. Moreover, $f_{b, \alpha, \rho, \lambda}$ is completely monotone.
%
\end{Example}
This example is particularly relevant for applications, especially in Finance, 
to the class of rough stochastic volatility models (see~\cite{el2019characteristic,ElEuchR2018,abi2019affine}), since in that case \(\alpha=H+\frac{1}{2}\in\big(\frac{1}{2},1\big),\)
where \(H\in(0,\frac12)\) denotes the Hurst coefficient representative of the 
roughness regime.
\subsection{Wiener-Hopf equations and forward Process}\label{subsect:wiener_hopf}	
We will always work under the following assumption in the rest of the paper.
\begin{assumption}[Kernel $K$ and resolvent $f_{\lambda}$ of $\lambda K$]\label{ass:resolvent} Throughout the paper,  we assume that the following assumption on the kernel $K$ and the (subsequently continuous) resolvent $f_\lambda$ of $\lambda K$ is satisfied:
	 \begin{equation}\label{eq:hypoRlambda}
	 	(\mathcal{K F}_\lambda)\quad  \left\{
	 	\begin{aligned}
	 		\text{(i)}\quad & K \text{ satisfies~\eqref{eq:contKtilde},~\eqref{eq:Kcont},~\eqref{eq:Kcontrol} and}\; \forall\, t>0,\; 0<L_K(t)<+\infty,\\
	 		\text{(ii)}\quad & f_\lambda \in L^1(\mathbb{R}_+) \text{ and } a:= 1-\int_{0}^{+\infty} f_\lambda(t)\,dt \in [0,1),
	 	\end{aligned}
	 	\right.
	 \end{equation}
	for some $\lambda>0$ under consideration, and assume that $\ell_\infty:=\lim_{t\to +\infty}  (\varphi-f_\lambda * \varphi)(t)\in\mathbb{R}$. 
\end{assumption}
\noindent{\bf Remark:} This assumption is satisfied by the $\rho-$exponential $\alpha-$fractional integration Kernel of Example~\ref{Ex:SolventGammaKernel} so that the function $f_\lambda$
has a finite Laplace transform given by~\eqref{eq:Lapl_Resolvent} and defined (at least) on $\mathbb{R}_+$ as an
$\overline{\mathbb{R}}_+$-valued function.
Indeed, the Fourier-Laplace transform associated to $K_{b, \alpha, \rho}$ always exists: Setting $v=(t+\rho)u$, we obtain
{\small
	\[L_{K_{b,\alpha,\rho}}(t)=\int_0^\infty be^{-(t+\rho)u}\frac{u^{\alpha-1}}{\Gamma(\alpha)}du=\frac{b(t+\rho)^{-\alpha}}{\Gamma(\alpha)}\int_0^\infty e^{-v}v^{\alpha-1}dv=b(t+\rho)^{-\alpha}.\]
}
so that using~\eqref{eq:Lapl_Resolvent}, we obtain
$L_{f_{b,\alpha,\rho,\lambda}}(t) = \frac{\lambda b(t+\rho)^{-\alpha}}{1 + \lambda b(t+\rho)^{-\alpha}}$. Consequently,
\begin{equation}\label{eq:limitResolGamma}
	a:= 1-\int_{0}^{+\infty}f_{b,\alpha,\rho,\lambda}(s)ds = 1- L_{f_{b,\alpha,\rho,\lambda}}(0)= \frac{\rho^{\alpha}}{\rho^{\alpha} + \lambda b } \in [0, 1).
	\vspace{-.3cm}
\end{equation}

\bigskip\noindent It is worth noting that, under Assumption~\ref{assump:kernelVolterra}~(i), if the kernel $K$ is completely monotone and not identically zero ($K\not\equiv 0$), then, for every $\lambda>0$, the resolvent $f_\lambda$ is completely monotone on $(0,\infty)$ by~\cite[Ch.~5, Theorem~3.1]{gripenberg1990}. In particular $f_\lambda \geq 0$ and $\displaystyle \int_0^\infty f_\lambda(t)\,dt \leq 1$ so that Assumption~\ref{ass:resolvent}~$(ii)$ holds true.

We now come to the main result of these preliminaries.
\begin{Proposition}[Wiener-Hopf transform and Forward Process]\label{lm:fwdprocess} For all \( s, t \in [0, T] \), we call \(\xi_t(s) := \mathbb{E}[X_s \mid \mathcal{F}_t]\) the \textit{ Forward process } of \( X \).
	Assume that assumptions ~\ref{assump:kernelVolterra} and ~\ref{ass:resolvent} are satisfied and that \( X = (X_t)_{t \in [0,T]} \) is a continuous weak solution of~\eqref{eq:intro}. 
	\((X_t)_{t\geq0}\) is solution of equation~\eqref{eq:intro} if and only if it is the solution of
	\vspace{-.3cm}
	\begin{equation}\label{eq:Volterrameanrevert2}
		X_t= X_0\big(\varphi(t)- \int_0^t f_{\lambda}(t-s)\varphi(s)ds\big) +\frac{1}{\lambda}\int_0^t f_{\lambda}(t-s)\mu(s)ds + \frac{1}{\lambda}\int_0^t f_{\lambda}(t-s)\sigma(s, X_s) \,dW_s.
	\end{equation}
	Then \(\xi_t(s)\) is an \(\mathcal{F}_t\) -martingale, and for all \( s, t \in [0, T] \) such that \(t\leq s\), we have
	\vspace{-.4cm}
	\begin{equation}\label{eq:Forward1}
		\mathbb{E}[X_s \mid \mathcal{F}_t] = X_0\varphi(s) +
		\int_0^s K(s-r) \left(\mu(r)-\lambda\mathbb{E}[X_r \mid \mathcal{F}_t]\right) \, \mathrm{d}r +
		\int_0^{t} K(s-r) \sigma(r, X_r) \, \mathrm{d}W_r.
	\end{equation}
	\vspace{-.4cm}
	Equivalently,
	\begin{equation}\label{eq:Forward2}
		\mathbb{E}[X_s \mid \mathcal{F}_t] =  X_0(\varphi(s)- \int_0^s f_{\lambda}(s-r) \varphi(r) \, dr ) + \frac{1}{\lambda} \int_0^s f_{\lambda}(s-r) \mu(r) \, \mathrm{d}r + \frac{1}{\lambda} \int_0^t f_{\lambda}(s-r) \sigma(r, X_{r}) \, \mathrm{d}W_r.
	\end{equation}
	\vspace{-.3cm}
	Moreover, the forward process \(\xi_t(s)\)
	satisfies the stochastic differential equation:
	\begin{equation}\label{eq:DiffForward} 
		\mathrm{d}\xi_t(s) = \frac1\lambda f_\lambda(s - t)\, \sigma(t,X_t)\, \mathrm{d}W_t,\; \xi_0(s) =  X_0(\varphi(s)- \int_0^s f_{\lambda}(s-r) \varphi(r) \, dr ) + \frac{1}{\lambda} \int_0^s f_{\lambda}(s-r) \mu(r) \, \mathrm{d}r.
	\end{equation}
	where \(\xi_0(s) :=\mathbb{E}[X_s \mid \mathcal{F}_0 ]\) is the initial condition (the expected process at future time $s$).
\end{Proposition}
\noindent {\bf Remark:}
It is easy to see that the initial forward process curve \( t \mapsto \xi_0(t) = \mathbb{E}[X_t\mid \mathcal{F}_0] \) is the solution to the Volterra equation

\centerline{$
	\xi_0(t) = X_0 \varphi(t) +
	\int_0^t K(t - s) (\mu(s) - \lambda \xi_0(s))\, ds,
	$}
\noindent or equivalently, if the kernel \(K\) is the $\alpha$-fractional kernel $ K(t) =K_{\alpha}(t) = \frac{t^{\alpha - 1}}{\Gamma(\alpha)} \mathbf{1}_{\mathbb{R}_+}(t)$, $ \alpha \in ( \frac{1}{2}, 1)$, the initial forward process is the solution to the fractional equation
\vspace{-.3cm}
\[
\xi_0(t) - X_0\varphi(t) + \lambda\, I^{\alpha}(\xi_0)(t) = I^{\alpha}(\mu)(t),
\]
where \( I^r \) denotes the fractional integral of order \( r \in (0,1] \).\footnote{Recall that the fractional integral of order \( r \in (0, 1] \) of a function \( f \) is \(I^r f(t) = \frac{1}{\Gamma(r)} \int_0^t (t - s)^{r - 1} f(s) \, ds,\)
	while the fractional derivative of order \( r \in [0, 1) \) is defined as \(	D^r f(t) = \frac{1}{\Gamma(1 - r)} \frac{d}{dt} \int_0^t (t - s)^{-r} f(s) \, ds,\)
	whenever the integrals exist.} As a result, setting \(\mu(t) = D^{\alpha}(\xi_0(t) - X_0\varphi(t)) + \lambda \xi_0(t)\)
ensures that the model whose variance process is defined by \((X_t)_{t\geq0}\) is consistent with any given initial forward process curve \( \xi_0(u) = \mathbb{E}[X_u\mid \mathcal{F}_0]\). Here, \( D^r \) denotes the fractional derivative of order \( r \in [0,1) \). (see also Proposition 3.1 and Remark 3.2 in \cite{ElEuchR2018}).

\section{Measure-extended Fourier-Laplace Functional for Inhomogeneous Affine Volterra Processes}
\label{sect-CharacteristicFunction}
\vspace{-.4cm}

In this section, we establish the representation result for the conditional Fourier-Laplace functional of the time-inhomogeneous affine Volterra equation and prove that it is exponential-affine in the past path. 
More generally, we consider the time-inhomogeneous affine Volterra equation where the diffusion coefficient is given by
\vspace{-.2cm}
\begin{equation}
	\sigma(t,x) = \varsigma(t)\,\sigma(x)
	\quad \text{with} \quad
	\sigma(x) :=\sqrt{\kappa_0+\kappa_1\,x}\,\quad\text{and} \quad  \kappa_1\geq0\,, \kappa_0\geq0.
	\label{eq:sigma-definition}
\end{equation}
and we assume throughout that the resulting Volterra
SDEs~\eqref{eq:Volterrameanrevert}
has (at least) one continuous \textit{non-negative} weak solution \(X = (X_t)_{t \geq 0}\) defined on some filtered stochastic basis \((\Omega, \mathcal{F}, (\mathcal{F}_t)_{t \geq 0}, \mathbb{P}, (W_t)_{t \geq 0})\) and satisfying~\eqref{eq:L^p-supBound}, e.g. as the $C$-weak limit of Hawkes processes as illustrated in 
\cite{EGnabeyeuR2025}. 



From now on, we set \(x_0(\cdot)\equiv X_0 \varphi (\cdot)\) and aim at analysing the measure-extended conditional Fourier-Laplace functional for any measure $m \in  \mathbb{M}_{\mathrm{lf}}^{-}(\R_+, \mathbb{C})\subset  \mathbb{M}_{\mathrm{lf}}(\R_+, \mathbb{C})$,  the subset of $\mathbb{C}-$-valued set functions $m \in  \mathbb{M}_{\mathrm{lf}}(\R_+, \mathbb{C})$ such that \(\operatorname{Re}(m)\) is negative on \(\R_+\):
\vspace{-.2cm}
\begin{equation}
	\E_{} \left[
	\exp \big( \int_0^T X_{T-s} \, m(\mathrm{d}s)
	\big)
	\,\Big|\, \mathcal{F}_t
	\right],
	\quad t \in [0, T].
	\label{eq:fourier-laplace}
\end{equation}
Remarkably, the functional in equation~\eqref{eq:fourier-laplace} admits a representation that can be expressed in terms of the solution of an associated time-inhomogeneous measure-extended Riccati-type Volterra equation.
This type of representation has been investigated in the convolution setting by \cite*[Theorem 4.3]{abi2019affine} and, in the non-convolution case with more general Volterra time-inhomogeneous diffusions, by \cite*[Theorem 2.1]{ackermann2022}. These results extend the classical exponential-affine transform for standard affine diffusions (see, e.g.,~\cite{duffie2003affine}) to the Volterra framework.  
In our context, however, the structure differs significantly: the initial condition is not deterministic but is instead given by a random function \( x_0 \) which evolves deterministically for \( t > 0 \). Concretely, for \( t > 0 \), the value \( x_0(t) \) is \( \mathcal{F}_0 \)-measurable, where \( \mathcal{F}_0 \) encodes the information up to time \( t = 0 \). Moreover, we consider a measure-extended version i.e. we extend the affine transform formula in the spirit of \cite{duffie2003affine,abi2019affine,ackermann2022} from $f \in L^1_{\mathrm{loc}}(\mathbb{R}_+; \mathbb{C})$ to $ m \in \mathbb{M}_{\mathrm{lf}}(\R_+, \mathbb{C})$ and provide a direct proof of both the affine transform representation and the existence of solutions to the associated measure-extended Riccati-Volterra equation via the Banach fixed-point theorem, bypassing the approximation and Sobolev-Slobodeckij compactness arguments employed in~\cite{friesen2022volterra}.
This provides the key tool to characterize the finite-dimensional distributions of the stationary process via a  measure-extended Riccati--Volterra equation.
To state the main formula in a synthetic form, let us define and then consider for a measure $m \in \mathbb{M}_{\mathrm{lf}}(\R_+, \mathbb{C})$, the following measure-extended Riccati--Volterra equation:
\begin{equation}\label{eq:measureFLplce}
	\begin{aligned}
		\forall \, m \in \mathbb{M}_{\mathrm{lf}}(\R_+, \mathbb{C}),\quad \psi(t) &= \int_{[0,t)} K(t-s)\,m(\dd s) + \int_0^t K(t-s)\,F(T-s,\psi(s))\,\dd s, \quad 0\leq t \leq T\\
		F(t, \psi) &= -\lambda \psi + \frac{\kappa_1}{2} \varsigma^2(t)\psi^2 \quad (t,\psi)\in \R_+\times \mathbb{C}.
	\end{aligned}
\end{equation}
where \( \lambda \in \mathbb{R}_+ \), and \( \varsigma : \mathbb{R}_+ \to \R \) is a given continuous function.

\medskip
\noindent {\bf Remark:}
Equation~\eqref{eq:measureFLplce} is written in a forward form. An equivalent expression in backward form is:
\vspace{-.2cm}
\begin{equation}\label{eq:measureFLplce_}
	\psi(T - t) = \int_t^T K(s-t) \, dm(T - s) + \int_t^T F(s, \psi(T - s)) K(s-t) \, \dd s.
\end{equation}
This 
formulation~\eqref{eq:measureFLplce_} is essential in problems where the system's behavior is determined by a known final state, allowing for the determination of the system's evolution by integrating backwards in time.
\begin{Lemma}\label{lem:RiccatiWithF} The inhomogeneous measure-extended Riccati--Volterra equation~\eqref{eq:measureFLplce} is equivalent to
	\begin{equation} \label{RicVolSqrt2}
		\psi(t) = \frac1\lambda\int_{[0,t)} f_{\lambda}(t-s) \, \,m(\mathrm{d}s) + \frac{1}{\lambda}\frac{\kappa_1}{2}\int_0^t \varsigma^2(T-s) \psi^2(s) f_{\lambda}(t-s) \, \mathrm{d}s, \quad t\in [0,T].
	\end{equation}
	where \(f_{\lambda} \) is solution to the equation \( 
	f_{\lambda} +\lambda K *f_{\lambda}=\lambda   K.\)
\end{Lemma} 
\noindent {\bf Proof:} Equation~\eqref{eq:measureFLplce} can be interpreted as a Wiener-Hopf equation with $x(t) = \psi(t)$ and
\[g(t) = \int_{[0,t)} K(t-s)\,m(\dd s) + \frac{\kappa_1}{2} \int_0^t K(t-s)\,\varsigma^2(T-s)\,\psi^2(s)\,\dd s = \Big(\big(m + \frac{\kappa_1}{2}\varsigma^2(T-\cdot)\,\psi^2 \big)*K\Big)(t).\]
From \cite[Proposition 2.4 (a)]{EGnabeyeu2025}, it follows that the expression for \( \psi \) is given for all \(t\geq0\) by:
\vspace{-0.3cm}
\begin{align*}
	&\, \psi(t) = g(t) - \int_0^t f_\lambda(t-s) g(s) \, \dd s =   \Big(\big(m + \frac{\kappa_1}{2}\varsigma^2(T-\cdot)\,\psi^2 \big)*K\Big)(t) \\
	&- \int_0^t f_\lambda(t-s)\Big[    \Big(\big(m + \frac{\kappa_1}{2}\varsigma^2(T-\cdot)\,\psi^2 \big)*K\Big)(s)\Big] \dd s = \big(m*K\big)(t) + \frac{\kappa_1}{2}\Big(\big(\varsigma^2(T-\cdot)\,\psi^2\big)*K \Big)(t) \\
	&\quad - \int_0^t f_{\lambda}(t-s) (m*K)(s) \, \dd s - \frac{\kappa_1}{2}\int_0^t f_{\lambda}(t-s) \left(K * \big(\varsigma^2(T-\cdot)\,\psi^2\big)\right)(s) \, \dd s\\
	&= \big(m*( K -f_\lambda*K )\big)(t) + \frac{\kappa_1}{2} \big((\varsigma^2(T-\cdot)\,\psi^2)*( K -f_\lambda*K )\big)(t)=\frac1\lambda\Big(m*f_\lambda + \frac{\kappa_1}{2} \big(\varsigma^2(T-\cdot)\,\psi^2\big)*f_\lambda\Big)(t)
\end{align*}
where we used commutativity and associativity (via regular Fubini's theorem) of convolution, and the 
last equality coming from the definition~\eqref{eq:flambda-eq} of \(f_{\lambda}\) namely \( 
f_{\lambda} +\lambda K *f_{\lambda}=\lambda   K.\)
The converse is true by convoluting equation~\eqref{RicVolSqrt2} by \( K \) and arguing similarly or by using the equivalence in the Wiener-Hopf equation since the solution is uniquely defined on \( \mathbb{R}_+ \) up to \( \text{d}t \)-almost sure equality.

In this section, we will work under the following assumptions:
	\begin{Assumption}[Integrability and Uniform H\"older Continuity]\label{ass:int_holregul}
	Assume that the kernel \( K \) is completely monotone, not identically zero and satisfies equations~\eqref{eq:contKtilde} and~\eqref{eq:Kcont} of Assumption~\ref{assump:kernelVolterra}. Moreover, assume that the resolvent \( f_{\lambda} \) of \( \lambda K \) satisfies:
	\begin{enumerate}
		\item[(i)] \textit{Integrability:} The condition below is satisfied for some $\beta>1$
		\begin{equation}\label{eq:contf_lambda}
			\int_0^{+\infty}
			\left(
			f_\lambda(t)^2+f_\lambda(t)^{2\beta}
			\right)\,dt<+\infty.
		\end{equation}
		
		\item[(ii)] \textit{H\"older Continuity:} There exists \( \vartheta \in (0, 1] \), \( C < +\infty \) such that
		\begin{equation}\label{eq:contHolder}
			\max_{i=1,2} \left[ \int_0^{+\infty} |f_{\lambda}(u + \bar{\delta}) - f_{\lambda}(u)|^i \, du \right]^{\frac{1}{i}} \leq C \bar{\delta}^{\vartheta}.
		\end{equation}
		
	\end{enumerate}
	\noindent Finally assume that  the initial random state  $X_0\in L^p(\mathbb P)$ for some $p>2$
	such that $p\big(\vartheta\wedge\frac{\beta-1}{2\beta}\big)>1,$ and that
	the function \( \varsigma : \mathbb{R}_+ \to \mathbb{R} \) is continuous, bounded on \( [0,\infty) \) and $ \lim_{t \to +\infty} \varsigma^2(t)= \varsigma^2_\infty<\infty$.
\end{Assumption}
\noindent {\bf Remark on the regularity:}  By Assumption~\ref{ass:int_holregul}, the resolvent \( f_{\lambda} \) is completely monotone (See also \cite[Ch. 5]{gripenberg1990}).
Let \( \alpha \in \left( \frac12, 1\right) \) and \(\rho\geq0\).
According to \cite[Proposition 5.1]{EGnabeyeu2025} (see also \cite[Proposition 5.1]{Pages2024} for the \(\alpha\)-fractional integration kernel), the resolvent of the  \(\rho\)-exponential \(\alpha\)-fractional integration kernel $K_{\alpha,\rho}$ of Example~\ref{Ex:SolventGammaKernel} satisfies Assumption~\ref{ass:int_holregul} for every \(\beta\in (1,\frac{1}{2(1-\alpha)}) \) with
\( \vartheta \in (0, \alpha - \frac{1}{2}) \) and $\lim_{
 	\vartheta\to\alpha-\frac12,
 	\beta\to\frac{1}{2(1-\alpha)}
 }
 \big(\frac{\beta-1}{2\beta}\wedge\vartheta\big)
 =
 \frac p2(2\alpha-1)>1$ whenever $p>\frac{2}{2\alpha-1}>2$.
 
\medskip\noindent
The results of this section will play a central role in
section~\ref{sect-longRun}, where we characterize the 
distributional properties of the limiting and stationary process 
associated with inhomogeneous affine Volterra equations.
\vspace{-.6cm}
\subsection{Analysis of the measure-extended generalized Riccati--Volterra Equation}
\noindent We derive the existence of a solution to the Riccati--Volterra equation \eqref{eq:measureFLplce}.
First note that, for any \(T>0\), from the definition of the convolution of a measure \(m \in \mathbb{M}_{\mathrm{lf}}(\R_+, \mathbb{C})\) and a function \(f : (0, T] \to \mathbb{R}\) in equation ~\eqref{eq:convolmeasure}, it is straightforward to check that for each \(p \in [1, \infty]\),
\(\|f * m\|_{L^p([0, T])} \leq \|f\|_{L^p([0, T])} \, |m|([0, T]).\)
Furthermore, if \(f\) is continuous on \([0, T]\) and either \(f(0)=0\) or \(m\)
is atomless, then the convolution \(f * m\) is also continuous on \([0, T]\).
\vspace{-.2cm}
\begin{Theorem}[Existence of solution to the Riccati--Volterra Equation]\label{theorem:riccati-extension}
	Under Assumption~\ref{ass:int_holregul}, 

	\medskip
	\noindent {(a)} For every $m \in \mathbb{M}^-_{\mathrm{lf}}(\R_+, \mathbb{C})$, 
	the measure-extended Riccati--Volterra integral equation~\eqref{eq:measureFLplce}-\eqref{RicVolSqrt2} admits a unique global solution $\psi = 	\psi(\cdot,m) \in L^2([0,T];\mathbb{C}) $, satisfying $\operatorname{Re}(\psi(t,m)) \leq 0$ for all \( t \in [0,T] \).

	\medskip
	\noindent {(b)} ($L^p$-bounds and regularity):  For each $p\in[1,2\beta]$ and for every $T > 0$, the unique solution $\psi = \psi(\cdot,m)$ 
	of~\eqref{eq:measureFLplce}-\eqref{RicVolSqrt2} satisfies the $L^p([0,T])$ a priori estimate
	\begin{equation}\label{eq:LpBounds}
		\left\|\psi(\cdot,m)\right\|_{L^p([0,T])}
		\leq \frac2\lambda \|f_\lambda\|_{L^p([0,T])} |m| ([0,T]) + \frac{\kappa_1}{2\lambda^3}
		\|\varsigma^2\|_{\infty,T} \|f_\lambda\|_{L^p([0,T])}  \|f_\lambda\|^2_{L^2([0,T])}|m| ([0,T])^2. 
	\end{equation}
	\noindent {(c)} (Continuity): If \(\beta\geq \frac32\) and \(K\in L^3_{\mathrm{loc}}\), the function \( \psi(\cdot, m) \) solution to~\eqref{eq:measureFLplce}-\eqref{RicVolSqrt2} is continuous at each \( t_0 \geq 0 \) for which the convolution \( (K * m)(\cdot):=  \int_{[0,\cdot)} K(\cdot-s)\,m(\dd s)\) is continuous at \( t_0 \).	
\end{Theorem}
\medskip
For clarity and conciseness, the proof of the above theorem is deferred to Appendix \ref{app:lemmata}, where the main technical results are presented. 
It is worth noting that
Theorem \ref{theorem:riccati-extension}~$(b)$ applied for \( T = \infty \) still provides the desired integrability as sketched below.
\begin{Corollary}\label{corol:riccati-extension}
	Let us consider a measure \( m \in \mathbb{M}^{-}_{\mathrm{lf}}(\R_+, \mathbb{C}) \) with \( |m|(\mathbb{R}_+) < \infty \). Then, under the same assumptions as in Theorem~\ref{theorem:riccati-extension}, there exists a unique global solution
	 \( \psi_\infty:= \psi_\infty(\cdot, m)\in L^1(\mathbb{R}_+; \mathbb{C}_{-}) \cap L^2(\mathbb{R}_+; \mathbb{C}_{-}) \), 
	 to the following  measure-extended Riccati--Volterra equation: 
	\vspace{-.2cm}
	\begin{align}
		\qquad\; \psi_\infty(t,m) &= \int_{[0,t)} K(t-s)\,m(\dd s) + \int_0^t K(t-s)\,F_\infty(\psi_\infty(s,m))\,\dd s, \quad t \geq 0
		\label{eq:measureFLplce2} \\
		F_\infty(\psi) &:= -\lambda \psi + \frac{\kappa_1}{2} \varsigma^2_\infty\psi^2 \quad\text{and}\quad \varsigma^2_\infty:= \lim_{t \to +\infty} \varsigma^2(t).
		\label{eq:FVarsigma_infty}
	\end{align} 
     Moreover, for every $T>0$, let $\psi_T$ denote the
     solution of the finite-horizon
     Riccati--Volterra equation~\eqref{eq:measureFLplce}
     on $[0,T]$, and let $\psi_\infty$ denote
     the unique solution of the limiting autonomous equation
     \eqref{eq:measureFLplce2}.
     Then
     \begin{equation}\label{eq:Stabil}
     \mathbf 1_{[0,T]}\psi_T
     \longrightarrow\psi_\infty
     \quad\text{in }
     L^1(\R_+;\mathbb C)\cap L^2(\R_+;\mathbb C)
     \qquad\text{as }T\to\infty.
     \end{equation}
\end{Corollary}

\subsection{Conditional Fourier-Laplace functional of inhomogeneous affine Volterra processes}

Similarly to the classical square-root process, there is a semi-explicit form for the
Fourier-Laplace transform of the inhomogeneous affine Volterra process, i.e., it is an affine process
on \( \mathbb{R}_+ \). 
The following theorem establishes the weak existence and uniqueness of \( \mathbb{R}_+ -\) solutions to~\eqref{eq:measureFLplce}, together with an expression for their Fourier-Laplace transform ~\eqref{eq:fourier-laplace}
in terms of the Riccati--Volterra equation~\eqref{eq:measureFLplce}.

\begin{Theorem}\label{T:VolSqrt}
	Let \(X = (X_t)_{t \geq 0}\) be a continuous $\R_+$-valued weak solution to the Volterra
	SDEs~\eqref{eq:Volterrameanrevert} with diffusion coefficient as in~\eqref{eq:sigma-definition}. Fix \(T>0\) and 
	suppose that Assumption~\ref{ass:int_holregul} holds.
	Consider a measure $m \in \mathbb{M}^{-}_{\mathrm{lf}}(\R_+, \mathbb{C})$
	such that
	 the measure-extended Riccati--Volterra equation ~\eqref{eq:measureFLplce} admits a unique global solution $\psi = \psi(\cdot,m) \in L^2([0,T];\mathbb{C}_-)$, 
	then: 
	
	\medskip
	\noindent \textit{1.} \label{T:VolSqrt:2}  The following exponential-affine transform formula holds for the measure-extended Fourier-Laplace transform of \( (X_t)_{t\in[0,T]} \):
	\vspace{-.2cm}
	{\small
		\begin{align}\label{eq:laplace}
			\mathbb{E}_{}\left[ \exp\left( \int_0^{T-t} X_{T-s} \, \,m( \mathrm{d}s) \right) \Big| \mathcal{F}_t \right]
			= \exp\left( \int_0^{T-t} \xi_t(T-s) \, m( \mathrm{d}s) + \frac{1}{2}
			\int_t^T  \varsigma^2(s)\sigma^2(\xi_t(s))\, \psi^2(T - s)\, \mathrm{d}s \right).
		\end{align}
	} 
	where the process \( \xi_t(s) \) is given by: \(\xi_t(s) = \xi_0(s) + \frac1\lambda \int_0^t  f_\lambda(s - r)\, \varsigma(r) \sigma(X_r)\, \mathrm{d}W_r, \quad \text{for } t \in [0,s],\) and
	\vspace{-.2cm} 
	\[\xi_0(s) = x_0(s)- \int_0^s f_{\lambda}(s-r) x_0(r)\, dr  + \frac{1}{\lambda} \int_0^s f_{\lambda}(s-r) \mu(r) \, \mathrm{d}r.\]
	\noindent \textit{2.}  
	The Fourier-Laplace functional of any weak solution $X$ to the inhomogeneous affine Volterra
	equation \eqref{eq:Volterrameanrevert}
	 admits the exponential-affine representation
	 {\small	
	 	\begin{align}
	 		&\ \mathbb{E}_{}\Big[ \exp\big( \int_0^T X_{T-s} \,m(\dd s) \big) \Big| X_0\Big]
	 		= \exp\Big( \int_0^T \mathbb{E}_{}[X_{T-s}\big| X_0] \,m(\dd s)
	 		+ \frac{1}{2} \int_0^T \varsigma^2(s)\, \sigma^2(\mathbb{E}_{}[X_s\big| X_0])\, \psi^2(T - s,m) \, \dd s \Big) \label{eq: MeasureFLplaceZero}
	 		\\ 
	 		&= \exp\Big( X_0 (\varphi* m)(T) + (\mu* \psi(\cdot,m))(T)
	 		+ X_0 \int_0^T F(s, \psi(T-s,m))\, \varphi(s)\, \dd s
	 		+ \frac{\kappa_0}{2} \int_0^T \varsigma^2(s) \psi^2(T-s,m) \, \dd s
	 		\Big). \label{eq: MeasureFLplaceZero II}
	 	\end{align}
	 }
	 	In particular, the stochastic Volterra equation \eqref{eq:Volterrameanrevert}-\eqref{eq:Volterrameanrevert2} admits a unique (in law) weak solution. 	
\end{Theorem}
\noindent{\bf Remark:} Note that the law of $X$ is determined by the Fourier-Laplace
transforms \(\E\big[\exp\left(\sum_{i=1}^n u_i\; X_{t_i}\right)\big]\) where $0 \leq t_1 < \dots < t_n$, $u_1,\dots, u_n \in \mathbb{C}_-$, and \(n\in \N\). Uniqueness in law in Theorem~\ref{T:VolSqrt}~$(2)$ thus follows by applying~\eqref{eq: MeasureFLplaceZero}-\eqref{eq: MeasureFLplaceZero II}  for the particular choice $m_{t_{1},\dots,t_{n}}(ds)=\sum_{i=1}^{n}u_{i}\delta_{t_{n}-t_{i}}(ds)$.

\smallskip \noindent
The main steps of the proof of Theorem~\ref{T:VolSqrt} are sketched in Appendix~\ref{app:lemmata}. 
Building on the preceding results, we will establish in section~\ref{sect-longRun} that the time-shifted process \((X_{t+u})_{t \geq 0}\) converges, as \(u \to \infty\), in finite-dimensional distributions to a stationary 
process. 
However, this limiting behavior does not imply that the original process \(X\) is itself stationary, nor does it yield information about the dynamics of the limiting and stationary process, an aspect that remains an open and challenging problem. To address this gap, we instead rely on a weaker yet analytically tractable framework: \textit{fake stationarity} in the sense of \cite{Pages2024, EGnabeyeu2025}. This notion allows us to capture key features of both the short- and long-term behavior of \(X\), despite the absence of full finite-time stationarity (see, e.g,~\cite{Pages2024}) or explicit dynamical characterization of the limit (see, e.g,~\cite{friesen2022volterra}).

\section{Finite-time fake stationary regimes of 
	affine Volterra processes.}\label{sect-FakeStationarity}
	We consider the stochastic affine Volterra integral equation~\eqref{eq:Volterrameanrevert} with the diffusion coefficient given in equation~\eqref{eq:sigma-definition}. We are interested in the so-called fake stationary regime, in the sense of~\cite{Pages2024}, which have previously been studied in the setting of quadratic diffusion coefficients for singular (rough in~\cite{EGnabeyeu2025}) and regular (long-memory in~\cite{GnabeyeuPages2026b}) kernels. In contrast to the quadratic settings, where the functional weak limiting distribution is shown to be at least $L^2$-stationary, the affine case is particularly appealing since, as we will demonstrate in Section~\ref{subsect:longRun}, this limiting distribution is in fact strongly stationary (i.e. in the classical sense, where the finite-dimensional distribution of the process is invariant under time shifts (see, e.g.,~\cite{friesen2022volterra, Jacquieretal2022} in the context of autonomous affine drift and diffusion coefficients)).
  Let $m_0\in\R^*_+$ and $v_0\in \R_+$.
	\begin{Definition}[Fake stationarity regimes]
		Let \( (X_t)_{t \geq 0} \) be a solution to the scaled Volterra equation in its form (\ref{eq:Volterrameanrevert}) starting from any \( X_0 \in L^2 (\P) \). Let \( \sigma(t, x) = \varsigma(t) \sigma(x) \) given in equation~\eqref{eq:sigma-definition}, where \(\varsigma\) is a (locally) bounded Borel function.
		\begin{enumerate}
			\item The process $(X_t)_{t\ge 0}$ exhibits a fake stationary regime of type I if it has constant mean, variance, and expected diffusion coefficient over time i.e.:
			{\small
				\begin{equation}\label{eq:fs1}
					\forall\, t\ge 0, \quad \mathbb{E}[X_t] =m_0, \quad \text{Var}(X_t) = v_0 \ge 0 \quad \mbox{and} \quad \bar{\sigma}^2(t) := \mathbb{E}[\sigma^2(X_{t})] := \bar{\sigma}^2_0 \ge 0.
				\end{equation}
			}
			\item The process $(X_t)_{t\ge 0}$ exhibits a fake stationary regime of type II if \( (X_t)_{t \geq 0} \) has the same marginal distribution, i.e., \( X_t \overset{d}{=}X_0 \) for every \( t \geq 0 \). 
		\end{enumerate}
	\end{Definition}
	\noindent These fake stationary regimes are characterized by an explicit closed-form expression for the deterministic forcing term $\varphi$, chosen such that $\forall\, t> 0, \; \mathbb{E}[X_t] =m_0$ as soon as $\mathbb{E}[X_0] =m_0$ and given by
	\begin{equation}
			\forall\, t\ge 0, \quad \varphi(t)   =1 + \lambda \int_0^t K(t-s) \Big( 1- \frac{\mu(s)}{\lambda m_0}  \Big) \, \dd s. \label{eq:CondMean}
	\end{equation}
	 together with a functional equation for the mean-reversion function $\mu$ and the multiplicative stabilizing term $\varsigma :=\varsigma_{\lambda,c}$, $\lambda,c>0$, which adjusts the volatility structure accordingly.
\subsection{Fake stationary regimes of affine Volterra processes} 
	 We henceforth assume throughout the paper that the function $\varphi$ called the forcing term, is specified by~\eqref{eq:CondMean}. Under these assumptions, for any initial condition $X_0 \in L^2(\mathbb{P})$, the Volterra SDE~\eqref{eq:Volterrameanrevert}--\eqref{eq:Volterrameanrevert2} admits the equivalent representation
	\begin{equation}
		X_t
		=
		X_0
		-
		\frac{X_0-m_0}{\lambda m_0}
		\int_0^t f_{\lambda}(t-s)\mu(s)\,ds
		+
		\frac{1}{\lambda}
		\int_0^t f_{\lambda}(t-s)\varsigma(s)\sigma(X_s)\,dW_s,
		\qquad X_0\perp\!\!\!\perp W,
		\label{eq:ConstMean}
	\end{equation}
	\begin{assumption}[Master Equation]\label{ass:MasterEq}  Let \(m_0\neq0\), $\lambda, c>0$. 
		We assume that the mean-reversion function \(\mu\) together with the multiplicative factor \(\varsigma\) in the diffusion term satisfy the following master equation on \((0,+\infty)\):
		\vspace{-.3cm}
		\begin{equation}\label{eq:MasterEq}
			\textit{($E_{\lambda, c}$)}: \quad\forall\, t\ge 0, \quad \Big(
			1-\frac{(f_\lambda * \mu)(t)}{\lambda m_0}
			\Big)^2
			=
			1-\frac{1}{c\lambda^2}
			(f_{\lambda}^2 * \varsigma^2)(t).
			\vspace{-.2cm}
		\end{equation}
	\end{assumption}
	\noindent From now on we define $h_{\lambda,c}(t)
	:=
	\sqrt{
		1-\frac{1}{c\lambda^2}
		(f_{\lambda}^2 * \varsigma^2)(t)
	}=1-\frac{(f_\lambda * \mu)(t)}{\lambda m_0}$. 
	Under the master identity~\eqref{eq:MasterEq} in assumption~\ref{ass:MasterEq}, the stochastic Volterra equation~\eqref{eq:ConstMean} starting from any initial condition
	\(X_0\in L^2(\mathbb P)\) admits the equivalent representation
	\vspace{-.3cm}
	\begin{equation}\label{eq:ConstMeanVar}
		X_t = X_0 -\big(X_0 - m_0\big) 	\big(
		1-h_{\lambda,c}(t)\big) +  \frac{1}{\lambda}\int_0^t f_{\lambda}(t-s)\varsigma(s)\sigma(X_{s})dW_s.
		\vspace{-.2cm}
	\end{equation}
	In this case, the corresponding time-inhomogeneous affine Volterra equation~\eqref{eq:ConstMeanVar} is referred to as a \textit{stabilized affine Volterra equation} or 
	as a \textit{fake stationary affine Volterra equation} if, in addition, $X_0$ is chosen such that $\mathbb{E}[X_0]=m_0$.
	
	\begin{Example}[Exponential and Fractional kernels]\label{Ex:FakeStatAutonome_Lm}\textcolor{white}{.}\\
		\noindent \textit{ 1. Exponential Kernel.} Let $\rho>0$ and consider the exponential kernel $K_{\rho}(t) = e^{-\rho t}\mbox{\bf 1}_{\R_+}(t)$, with
		$f_{\rho, \lambda}$ defined by Equation~\eqref{eq:DerivSolventGammaKernel} of Example~\ref{Ex:SolventGammaKernel}.
		\begin{itemize}
			\item[$(a)$] If \(\varsigma\equiv1\), under the condition $c\geq \frac{1}{2(\lambda + \rho)}$, Equation~\eqref{eq:MasterEq} is satisfied with
			$\mu(t)=m_0\Big[
			(\lambda+\rho)\bigl(
			1- \big(1-
			\frac{1-e^{-2(\lambda+\rho)t}}{2c(\lambda+\rho)}\big)^{\frac12}
			\bigr)
			+
			\frac{e^{-2(\lambda+\rho)t}}
			{2c} \bigl(1-\frac{1-e^{-2(\lambda+\rho)t}}{2c(\lambda+\rho)}\bigr)^{-\frac12}
			\Big]$.	
			It is straightforward that \(\mu\) is continuous on \(\mathbb{R}_+\). 
			One checks that \(\mu(0)=\frac{m_0}{2c}\) and that  $\mu_\infty:=\lim_{t\to+\infty}\mu(t)
			=
			m_0 (\lambda+\rho)
			\big(
			1-(1-
			\frac{1}{2c(\lambda+\rho)})^{\frac12}
			\big).$
			\item[$(b)$] If \(\mu\equiv\lambda m_0\), applying Leibniz's rule to the master equation yields \(\varsigma_{\rho,\lambda,c}^2(t)
			=
			\frac{1}{\lambda^2}
			\Bigl(
			g_{\lambda,c}'(t)
			+
			2(\lambda+\rho)g_{\lambda,c}(t)
			\Bigr),\)
			where \(g_{\lambda,c}(t)
			:=
			c\lambda^2
			\big[
			1-
			\big(
			1-
			\frac{\lambda}{\lambda+\rho}
			\bigl(1-e^{-(\lambda+\rho)t}\bigr)
			\big)^2
			\big].\)
			A straightforward computation then shows that \(\varsigma_{\rho,\lambda,c}^2(t)
			=
			\frac{2c\lambda}{\lambda+\rho}
			\big(
			\lambda+2\rho-\rho e^{-(\lambda+\rho)t}
			\big).\)
			Moreover, it is readily verified that
			\(\varsigma_{\rho,\lambda,c}^2\) is continuous, strictly positive and non-decreasing on \(\mathbb{R}_+\). Furthermore,  \(	\varsigma_{\rho,\lambda,c}^2(0)=2c\lambda,\)
			and \(\varsigma_{\rho,\lambda,c,\infty}^2
			:=
			\lim_{t\to\infty}
			\varsigma_{\rho,\lambda,c}^2(t)
			=
			\frac{2c\lambda(\lambda+2\rho)}{\lambda+\rho}.\)
		\end{itemize}
		\medskip
		\noindent \textit{ 2. Fractional integration Kernel.}
		Let $\alpha \in (\frac12, 1)$ and consider the Kernel $K_{\alpha}(t) = \frac{u^{\alpha-1}}{\Gamma(\alpha)} \mbox{\bf 1}_{\R_+}(t)$, with
		$f_{\alpha, \lambda}$ defined by Equation~\eqref{eq:DerivSolventGammaKernel} of Example~\ref{Ex:SolventGammaKernel}. 
		
		\begin{itemize}
			\item[$(a)$] If \(\varsigma\equiv1\), Equation~\eqref{eq:MasterEq} is satisfied provided $c\geq \frac{1}{\lambda^2} \int_0^\infty f_{\alpha,\lambda}^2(s)\,ds$, with $\mu(t)
			=
			g_{\alpha,\lambda}(t)
			+ \frac{1}{\lambda} D^{\alpha} g_{\alpha,\lambda}(t)$ where $g_{\alpha,\lambda}(t)
			=
			\lambda m_0
			\big(1-(1-\frac{1}{c\,\lambda^2}
			\int_0^t f_{\alpha,\lambda}^2(s)\,ds)^{\frac12}
			\big)$
			and \( D^{\alpha} := \frac{d}{dt} I^{1 - \alpha} \) where \( D^\alpha \) and \( I^{1 - \alpha} \) denote, respectively, the Riemann--Liouville fractional derivative of order \( \alpha \), and the Riemann--Liouville fractional integral of order \( 1 - \alpha \) (see \cite[Chapter~2]{samko}).
			One checks that  \(\mu\) is continuous on \((0,+\infty)\) and that $\lim_{t\to+\infty}\mu(t)
			=
			\lambda m_0
			\big(
			1-(1-\frac{1}{c\,\lambda^2}
			\int_0^\infty f_{\alpha,\lambda}^2(s)\,ds)^{\frac12}
			\big)$. 
			\item[$(b)$]  If \(\mu \equiv \lambda m_0\),
			the function \( \varsigma^2 = \varsigma^2_{\alpha,\lambda,c} \) exists as a non-negative, non-decreasing concave function, on \( (0, +\infty) \), such that \(\varsigma^2_{\alpha,\lambda,c}(0)=0\) and $\lim_{t \to +\infty} \varsigma^2_{\alpha,\lambda,c}(t) = \frac{c\lambda^2}{\|f_{\alpha,\lambda}\|^2_{L^2(\text{Leb}_1)}}$ (see e.g.~\cite[Proposition 5.3]{EGnabeyeu2025}). Moreover, it follows from \cite[Proposition 6.2]{Pages2024} that, by setting \(a_k = \frac{1}{\Gamma(\alpha k + 1)},
			b_k = \frac{1}{\Gamma(\alpha(k + 1))}, \; k \geq 0\), the function $\varsigma^2_{\alpha,\lambda, c}$ admits the following power series expansion:
			\(\varsigma^2_{\alpha,\lambda, c}(t) = c \lambda^{2-\frac1\alpha}\varsigma_\alpha^2(\lambda^{\frac1\alpha} t)\) where \(\varsigma_{\alpha}^2(t):= 2\,t^{1-\alpha}\sum_{k\ge 0} (-1)^k c_k t^{\alpha k}\) and  the coefficients $(c_k)_{k\geq0}$ are defined by the following recurrence formula:\[
			c_0=\frac{\Gamma(\alpha)^2}{\Gamma(2\alpha-1)\Gamma(2-\alpha)} \quad \textit{ and for every} \quad k\ge 1,
			\]
			{\small 
				\begin{equation}\label{eq:ck2}
					c_k = \frac{\Gamma(\alpha)^2\Gamma(\alpha(k+1))}{\Gamma(2\alpha-1)\Gamma(\alpha k+2-\alpha)}\left[ (a*b)_k- \alpha(k+1)\sum_{\ell=1}^kB\big(\alpha(\ell+2)-1,\alpha(k-\ell-1)+2\big) (b^{*2})_{\ell}  c_{k-\ell}  \right].
				\end{equation}
			}
			where for two sequences of real numbers \( (u_k)_{k \geq 0} \) and \( (v_k)_{k \geq 0} \), the Cauchy product is defined as \( (u * v)_k = \sum_{\ell = 0}^k u_\ell v_{k - \ell} \) and \( B(a, b) = \int_0^1 u^{a-1}(1 -
			u)^{b-1} du\) denoting the beta function. 
		\end{itemize}
	\end{Example}
\noindent {\bf Remark:} 
Note that the fractional power series \( \sum_{k \geq 0} c_k t^{\alpha k} \) defined by the coefficients~\eqref{eq:ck2}, is obtained by exploiting regular variation and applying Tauberian theorems to the Laplace transform of the functional equation~\eqref{eq:MasterEq}. Its convergence radius \( \rho_\alpha := \left( \liminf_k \left( |c_k|^{1/k} \right) \right)^{-1/\alpha} \), is infinite. 
As a consequence, the expansion which define \(\zeta_{\alpha,\lambda,c}^2\) converges for all \( t \in \mathbb{R}^+ \), and in fact, for all \( t \in \mathbb{R} \). 

\medskip\noindent
We now come to the main result of this section, which  shows what are the consequences of the three constraints in Equation~\eqref{eq:fs1} in the context of the Volterra SDE~\eqref{eq:Volterrameanrevert}--\eqref{eq:Volterrameanrevert2}.
\begin{Proposition}[(Fake stationary regimes (types I and II) and first asymptotics) ]\label{prop:mainStab}
	Let $\lambda, c > 0$, $m_0 \in \mathbb{R}^*$ and
	assume $(\mathcal{K F}_\lambda)$ is in force.
	Let \(X = (X_t)_{t \geq 0}\) be the one-dimensional solution of the (stabilized) Volterra Equation~\eqref{eq:ConstMeanVar} (Equation~\eqref{eq:Volterrameanrevert} with conditions~\eqref{eq:CondMean} and~\eqref{eq:MasterEq} 
	in force) with constant $ \lambda, c, m_0$, $\sigma$ given by~\eqref{eq:sigma-definition} and  \( X_0 \in L^2(\mathbb{P}) \) with
	 \(\mathbb{E}[X_0] = m_0\) and \(
	\mathrm{Var}(X_0) = v_0=c\sigma^2(m_0)\). 
	\begin{enumerate}
		\item  If $\sigma(m_0) =0$ (i.e. the diffusion coefficient  $ \sigma$ is degenerated), then the solution $X_t = m_0$ $\mathbb{P}$-$a.s.$  represents a trivial fake stationary (degenerate) regime (of type I).
		\item  If $\kappa_1 = 0$ (Volterra Ornstein-Uhlenbeck process with constant $\sigma(x) = \sqrt{\kappa_0}$), then the solution $(X_t)_{t \geq 0}$ has a constant mean $m_0$ and variance $v_0$. Consequently:
		\begin{itemize}
			\item[$(a)$] The process exhibits a fake stationary regime of type I i.e.
			
			\centerline{$
				\forall t \geq 0, \quad \mathbb{E}[X_t] = m_0, \quad \text{Var}(X_t) = v_0 =c\kappa_0.
				$}
			\item[$(b)$] Furthermore,  if $X_0 \sim \nu^* :=  \mathcal{N}\left( m_0, v_0\right)$, this represents a fake stationary regime of type II, since in this case,  $X_t \sim X_0$ for all $t \geq 0$. ($(X_t)_{t \geq 0}$ is a Gaussian process with a fake stationary regime of type II.). $\nu^*$ is the 1-marginal distribution.
		\end{itemize}
		\item If $\kappa_1 \neq 0$ and $\sigma(m_0) \neq 0$, the solution \( (X_t)_{t \geq 0} \) to the Volterra equation~\eqref{eq:ConstMeanVar} has a fake stationary regime of type I, in the sense that\[
		\forall\, t\ge 0, \quad 
		\mathbb{E}[X_t] = m_0, \; \text{Var}(X_t) = v_0= c\sigma^2(m_0),\; \text{and}\;\; \mathbb{E}[\sigma^2(X_t)] = \bar{\sigma}_0^2 =\sigma^2(m_0).\]
	\end{enumerate}
\end{Proposition}
\noindent {\bf Remark:} (\textit{On Inhomogeneous Volterra square-root process.})
In particular, if $\kappa_0=0$ , then \(\sigma(x) = \nu \sqrt{x}\), where $\nu =  \sqrt{\kappa_1} $.
The resulting Volterra equation has a fake stationary regime of type I, in the sense that

\centerline{$
	\forall\, t\ge 0, \quad 
	\mathbb{E}[X_t] = m_0, \; \text{Var}(X_t) = v_0= c\nu^2m_0,\; \text{and}\;\; \mathbb{E}[\sigma^2(X_t)] = \bar{\sigma}_0^2 = \nu^2m_0.
	$}

\medskip
\noindent {\bf Proof of Proposition \ref{prop:mainStab}.}
Let us consider a weak solution on the whole non-negative 
real line of the Stochastic Volterra Equation~\eqref{eq:ConstMeanVar} with volatility term $\sigma(t,x)= \varsigma_{\lambda,c}(t)\sigma(x)$ given in~\eqref{eq:sigma-definition} and
starting from any $X_0 \in L^2(\mathbb{P})$ such that $\mathbb{E}[X_0] = m_0$ and $Var[X_0] = v_0$.
We check that, $\mathbb{E}[X_t]=\mathbb{E}[X_0] = m_0$ for every $t > 0$. As for the variance, with the constant mean
and taking advantage of
the identity~\eqref{eq:MasterEq}, one checks that 
the variance satisfies the identity
\begin{equation}\label{eq:conditionVar}
	\operatorname{Var}(X_t) - c\sigma^2(m_0)
	=  h_{\lambda,c}^2(t) \bigl(v_0 - c\sigma^2(m_0) \bigr),
	\qquad t \geq 0.
\end{equation}
Since $h_{\lambda,c}(t) < 1$ for at least some $t > 0$ large enough, we obtain that $\operatorname{Var}(X_t)$ is  constant whenever $v_0 = c\sigma^2(m_0)$.
Conversely, one checks
that this constant value $v_0$ for the variance solves the above equation~\eqref{eq:conditionVar}. 
This proves claims $2.(a)$ and $3.$
For the first claim, if $\sigma(m_0)=0$, we get $v_0 = 0$ and consequently, $\operatorname{Var}(X_t) = 0$ for every $t \geq 0$. 
It follows that $X_t = m_0$, $\mathbb{P}$-a.s. 
For the second claim $2.(b)$, with $\nu^*$ as the starting value, $(X_t)_{t \geq 0}$ is a Gaussian process with constant mean and variance, which characterize its one-dimensional marginal distributions.
 \qed
\vspace{-.3cm}
\subsection{The fake stationary Volterra Heston model and its characteristic functions.}\label{subsect:FS-RH}
We now examine an affine Volterra process with state space $\mathbb{R} \times \mathbb{R}_+$, which can be interpreted as an extension of the Volterra Heston \cite{abi2019affine} and the classical Heston \cite{Heston1993} stochastic volatility model widely used in financial mathematics. It represents a special case of the more general inhomogeneous Volterra-Heston model introduced in \cite{ackermann2022}, where the diffusion coefficient is time-dependent and separable in the state variable and time, except that conditions~\eqref{eq:CondMean} and~\eqref{eq:MasterEq} are in force. We refer to this as the \textit{fake stationary Volterra Heston model}. That new model, also emerges for example
as the continuous-time limit of a suitably rescaled non-Markovian linear Hawkes model for high frequency assets prices, as proposed in
~\cite{EGnabeyeuR2025}. In this setting, we define the process $X = (\log S, V)$, where $S$ denotes the asset price and $V$ its variance process, governed by
\vspace{-.2cm}
\begin{equation}\label{E:HestonlogS_V}
	\begin{aligned}
		&\qquad\qquad\qquad \frac{dS_t}{S_t}
		= \sqrt{V_t}\,(\sqrt{1-\rho^2}\,dW^{(1)}_t+\rho\,dW^{(2)}_t), 
		\qquad S_0\in(0,\infty), \\[0.1em]
		& V_t = V_0\varphi(t) +  \int_0^t K(t-s) \left( (\mu(s) - \lambda V_s)ds 
		+  \nu \varsigma(s)\sqrt{V_s}\,dW^{(2)}_s\right), V_0\perp\!\!\!\perp (W_1, W_2), \qquad \varsigma=\varsigma_{\lambda,c}.
	\end{aligned}
\end{equation}
with $\varphi$ given by~\eqref{eq:CondMean} and the deterministic function $\mu$ together with $\varsigma$ satisfies the master identity~\eqref{eq:MasterEq} in assumption~\ref{ass:MasterEq} for $\lambda, c>0$,
 $(W_1, W_2)$ is a two-dimensional standard Brownian motion with independent components, $\rho \in [-1, 1]$, and  
  $\nu \in \R_+$ 
   More precisely, the process \(V\) in~\eqref{E:HestonlogS_V} can take the form~\eqref{eq:ConstMean}-\eqref{eq:ConstMeanVar}.
Note that the joint dynamics in equation~\eqref{E:HestonlogS_V} for the asset's price process
$S$ (e.g., the SPX) and its spot variance are expressed in this form (risk-neutral prices) because 
volatility models are typically formulated in terms of the so-called forward prices, that is, 
the process $F_t = e^{(r-q)t} S_t$ with $r$ and $q$ the interest and dividend rates.
Once $V$ is determined, the asset price process $S$ follows accordingly. Moreover  by applying It\^o's formula, one can verify that for every $t \in [0,T]$, the log-price process satisfies

\begin{equation}\label{eq:log_S_t}
	\log(S_t)=\log(S_0)+\int_0^t\sqrt{V_s}\left(\sqrt{1-\rho^2}dW^{(1)}_s+\rho dW^{(2)}_s\right)-\int_0^t\frac{V_s}{2}\dd s.
\end{equation}
Hence, $X = (\log S, V)$ constitutes an affine Volterra process which  evolves according to the system

\begin{equation}\label{eq:Systemlog_S_V}
	\begin{split}
		\begin{pmatrix}
			\log(S_t) \\ V_t
		\end{pmatrix}
		&= \begin{pmatrix}
			\log(S_0) \\ V_0 \varphi(t)
		\end{pmatrix}
		+ \int_0^t
		\begin{pmatrix}
			1 & 0 \\ 0 & K(t - u)
		\end{pmatrix}
		\left[
		\begin{pmatrix}
			0 \\ \mu(u)
		\end{pmatrix}
		+
		\begin{pmatrix}
			0 \\ 0
		\end{pmatrix}
		\log(S_u)
		+
		\begin{pmatrix}
			- \frac{1}{2} \\ -\lambda
		\end{pmatrix}
		V_u
		\right]
		du \\
		&\quad + \int_0^t
		\begin{pmatrix}
			1 & 0 \\ 0 & K(t - u)
		\end{pmatrix}
		\begin{pmatrix}
			\sqrt{1 - \rho^2} & \rho \\
			0 & \nu \varsigma(u)
		\end{pmatrix}
		\sqrt{V_u} \, dW_u, \quad t \in [0, T],\,\quad \varsigma = \varsigma_{\lambda,c}.
	\end{split}
\end{equation}
Fix \(T>0\), \(u=(u_1,u_2)\in(\mathbb C^2)^*\), and
\(f=(f_1,f_2)\in L^1_{\mathrm{loc}}(\mathbb R_+;(\mathbb C^2)^*)\).  For the fake stationary Volterra-Heston model, the Riccati-Volterra equation \eqref{eq:measureFLplce} with \(m \in \mathbb{M}_{\mathrm{lf}}(\R_+, (\mathbb C^2)^*)\) given by \(m(\dd s) = u\,\delta_0(\dd s) +f(s)\,\mathrm{Leb}_{\mathbb{R}_+}(\mathrm{d}s)\),
 in dimension 2 (see also \cite[Equation 12]{ackermann2022}) takes the form:
\begin{equation}
	\begin{split}\label{eq:riccati_heston_gen}
		\psi_1(t) & = u_1 + \int_0^t f_1(s) \dd s,\\
		\psi_2(t) & = u_2 K(t) + \int_0^t K(t-s)\;\Big( f_2(s) +  \frac{1}{2} (\psi^2_1(s)-\psi_1(s)) \\
		& \quad + \left( \rho \nu \varsigma(T-s) \psi_1(s) - \lambda \right)  \psi_2(s) + \frac{\nu^2}{2} \varsigma^2(T-s) \psi_2^2(s) \Big) \dd s,  \quad t \in [0,T],\,\quad \varsigma = \varsigma_{\lambda,c}.
	\end{split}
\end{equation} 

\begin{Proposition} \label{prop:existence_VolterraHestonS_V} Let \(T>0\).
	Suppose that Assumptions~\ref{assump:kernelVolterra} and \ref{ass:int_holregul} are in force.
	Consider a fake stationary Volterra Equation with constant $\lambda, c > 0$, $m_0 \in \mathbb{R}^*$ and conditions~\eqref{eq:CondMean} and~\eqref{eq:MasterEq}
	in force.	
	\begin{enumerate}
		\item\label{T:VolHeston:1} The stochastic Volterra system ~\eqref{E:HestonlogS_V} admits a \(\R\times\R_+-\)valued unique in law 
		continuous weak solution $(\log S, V)$ with values in $\mathbb{R} \times \mathbb{R}_+$, for any initial state $(\log S_0, V_0) \in \mathbb{R} \times \mathcal{P}(\mathbb{R}_+)$. Moreover, if $\mu\in L^2_{\rm loc}([0,\infty))$, the sample paths of $V$ are \( \left(\theta \wedge \widehat \theta - \frac{1}{p} - \eta \right) \)-H\"older pathwise continuous (modulo  \( P \)-indistinguability) for sufficiently small \( \eta > 0 \).
		
		\item\label{T:VolHeston:2} Let $u \in (\mathbb{C}^2)^*$ and $f \in L^1([0,T], (\mathbb{C}^2)^*)$,
		such that \(\text{${\rm Re\,} \psi_1 \in[0,1]$, ${\rm Re\,} u_2 \le0$, and ${\rm Re\,} f_2\le0$.}\)
		where $\psi_1$ solves the first relation in \eqref{eq:riccati_heston_gen}. Then the second equation of the Riccati--Volterra equation ~\eqref{eq:riccati_heston_gen} admits a unique solution $\psi_2 \in L^2([0,T], \mathbb{C})$ on \([0,T]\)
		with ${\rm Re\,}\psi_2 \le 0$. Furthermore, the exponential-affine representation~\eqref{eq: MeasureFLplaceZero II} is valid for $(\log S, V)$ for any initial state $(\log S_0, V_0) \in \mathbb{R} \times \mathcal{P}(\mathbb{R}_+)$ with \(\mathbb{M}_{\mathrm{lf}}(\R_+, (\mathbb C^2)^*) \ni m(\dd s) := u\,\delta_0(\dd s) +f(s)\,\mathrm{Leb}_{\mathbb{R}_+}(\mathrm{d}s)\).

		\item\label{T:VolHeston:3} The process $S$ solution of the first equation in~\eqref{E:HestonlogS_V} is a true martingale and can be written:
		\begin{equation}\label{eq:formula_S_inhom_Heston}
			S_t=S_0\exp\left(-\int_0^t \frac{V_s}{2}\dd s + \int_0^t \sqrt{V_s}\left(\sqrt{1-\rho^2}dW^{(1)}_s +\rho dW^{(2)}_s\right)\right), \quad t\in [0,T].
		\end{equation}
	\end{enumerate}
\end{Proposition}
The preceding results allow us to derive the Fourier--Laplace transform at time zero for the fake stationary Volterra-Heston model with fractional integration kernel. 
Such kernels are particularly relevant in capturing either the short or the long-term memory of volatility phenomena, as seen for instance in the rough Heston model ( \cite{el2019characteristic,ElEuchR2018}). The next result extends those of \cite{el2019characteristic,ElEuchR2018} and \cite[Example 7.2]{abi2019affine} by incorporating time-dependent drift $\mu(\cdot)$ and diffusion coefficient $\varsigma(\cdot)\sigma(\cdot)$ to give the full Fourier-Laplace functional with integrated log-price and variance.
Consider $(\log S_0, V_0) \in \mathbb{R} \times \mathcal{P}(\mathbb{R}_+)$, any initial state of the Volterra system ~\eqref{E:HestonlogS_V}.
\begin{Corollary}[Fake stationary Rough Heston model]\label{corol:StabilizedRHeston} Let \(T>0\).
	Let $\alpha\in (\frac{1}{2},1)$, and consider the \(\alpha-\) fractional integration kernel $K_{\alpha}(t)=\frac{t^{\alpha-1}}{\Gamma(\alpha)}$ for $t \in (0,T]$.  
	Suppose $u \in (\mathbb C^2)^*$ and $f \in L^1\left( [0,T], (\mathbb C^2)^* \right)$ satisfy the conditions $\Re \psi_1 \in [0,1]$, $\Re u_2 \leq 0$, and $\Re f_2 \leq 0$, where $\psi_1 = u_1 + \int_0^{\cdot} f_1(s)\,\dd s$. 	
	Then there exists a unique function $ \psi_2 \in L^2([0,T],\mathbb C)$ solving the fractional Riccati equation
	\begin{equation}\label{eq:riccati_heston_final}
		\begin{aligned}
			&\	(D^{\alpha} \psi_2)(t) = f_2(t) + \frac{1}{2}\left( u_1^2 - u_1 + (2u_1-1) \int_0^{t} f_1(s)\,\dd s + \left( \int_0^{t} f_1(s)\,\dd s \right)^2 \right)  + \frac{\nu^2}{2} \varsigma^2(T - t)\, \psi_2^2(t) \\
			&\qquad + \left( \rho \nu \varsigma(T - t) \left( u_1 + \int_0^{t} f_1(s)\,\dd s \right) - \lambda \right) \psi_2(t) , \quad t \in [0, T], \; \varsigma = \varsigma_{\lambda, c}, \;
			(I^{1 - \alpha} \psi_2)(0) = u_2.
		\end{aligned}
	\end{equation}
	leading to the full Fourier--Laplace representation for the integrated log-price and variance:
	{\small 
		\begin{equation}\label{eq:time_zero_frac_heston_all}
			\begin{split}
				& \E_{(\log(S_0),V_0)}\left[ \exp\left( u_1 \log(S_T) + u_2 V_T + \int_0^T f_1(T-u) \log(S_u) \,du   + \int_0^T f_2(T-u) V_u \,du \right) \right] \\
				& = 
				\exp\left( \phi(T) + \left(u_1 + \int_0^T f_1(s)\,\dd s \right)\log(S_0) + \left((I^{1-\alpha} \psi_2)(T) -\lambda \int_0^T \Big( \frac{\mu(u)}{\lambda m_0} - 1 \Big) \psi_2(T-u) \,du \right) V_0 \right).
			\end{split}
		\end{equation}
	}
	\noindent with \(\phi\) defined as $,  \forall \, t \in [0,T],\;\phi(t)= \int_0^t \mu(s)\psi_2(t-s)\,\dd s$, $\phi(0)=0$ and \( D^{\alpha} = \frac{d}{dt} I^{1 - \alpha} \) where \( D^\alpha \) and \( I^{1 - \alpha} \) denote, respectively, the Riemann--Liouville fractional derivative of order \( \alpha \), and the Riemann--Liouville fractional integral of order \( 1 - \alpha \) (see \citep[Chapter~2]{samko}). 
	In the particular case where \( \mu(t) = \mu_0 = \lambda m_0\;\forall t \geq 0 \) (so that \( \varphi \equiv 1 \) in\eqref{eq:CondMean}), 
	 the above Fourier--Laplace representation simplifies to
	\begin{equation}\label{eq:time_zero_frac_heston}
		\begin{split}
			& \E_{(\log(S_0),V_0)}\left[ \exp\left( u_1 \log(S_T) + u_2 V_T + \int_0^T f_1(T-u) \log(S_u) \,du   + \int_0^T f_2(T-u) V_u \,du \right) \right] \\
			& = 
			\exp\left( \phi(T) + \left(u_1 + \int_0^T f_1(s)\,\dd s \right)\log(S_0) +
			(I^{1-\alpha} \psi_2)(T) V_0 \right).
		\end{split}
	\end{equation}
\end{Corollary}

\noindent{\textbf{Practitioner corner:} 
	It is worth noting that the characteristic function of the \textit{fake stationary rough Heston model} introduced above can be computed by solving integral ordinary differential (ODE), specifically the so-called \textit{Riccati equations}. 
	See, for instance, the pioneering works~\cite{DiethelmFordFreed2002,DiethelmFordFreed2004}, which develops a predictor-corrector method for the numerical approximation of fractional ODE, known as the generalized Adams--Bashforth--Moulton method, or \textit{fractional Adams method}. 	 
	Plug the numerical solution of the riccati equation, into~\eqref{eq:time_zero_frac_heston}, to obtain the characteristic function. Then, classical methods can be applied to compute call (resp. put) option prices
; see \cite{CarrMadan1999, Lewis2001,Itkin2005} and the survey \cite{Schmelzle2010} (see, e.g. Figure~\ref{fig:Impl_Vol_gen}, which display the implied volatility surface generated by the model using an inverse Fourier transform).

	\vspace{-.4cm}
	\section{ 
		Long run confluence, limiting distributions and asymptotics}\label{sect-longRun}
		\noindent Throughout the rest of the paper we assume that the following assumption holds for the mean-reversion function.
		\begin{Assumption}[On the mean-reversion function]\label{ass:Mean_revertF}
		 We assume that the mean-reversion
		function \(\mu\in L^2_{\mathrm{loc}}(\mathbb R_+)\) is continuous on $(0,\infty)$ and converges at infinity to $\mu_\infty\in \R_+$. We also assume \(m_0\in \R^*_+\).
	\end{Assumption}
	\noindent One checks under Assumptions~\ref{ass:resolvent}~$(ii)$-\ref{ass:Mean_revertF}-\ref{eq:MasterEq} and condition~\eqref{eq:CondMean} that
	$\ell_\infty:=\lim_{t\to +\infty}  (\varphi-f_\lambda * \varphi)(t) =\lim_{t\to +\infty} 
	1-
	\frac{(f_\lambda * \mu)(t)}{\lambda m_0}= 1 - \frac{1-a}{\lambda m_0} \mu_\infty \in[0,1]$, so that \(\ell_\infty=0\) provided that \(m_0=(1-a)\frac{\mu_\infty}{\lambda}\).
	\subsection{ $L^p$-Moment control and $L^2$-Contraction.}\label{subsect-confluence}
	 For the sake of generality, we establish the results of this subsection within the broader framework of stochastic Volterra equations~\eqref{eq:Volterrameanrevert}
	 with affine drift and a space-dependent diffusion coefficient \(\sigma\), assumed to be locally \(\gamma\)-H\"older continuous in the spatial variable for some \(\gamma\in[\frac12,1]\). More precisely, we assume that there exists a constant
	 \([\sigma]_{\text{H\"ol}}>0\) such that
	 \vspace{-.1cm}
	 \begin{equation}\label{eq:subHolder} 
	 	|\sigma(x)-\sigma(y)|
	 	\le [\sigma]_{\text{H\"ol}} |x-y|^{\gamma},
	 	\vspace{-.1cm}
	 \end{equation}
	 for all \(t\in[0,T]\) and \(x,y\in\mathbb{R}\).
	 	The existence of strong solutions and pathwise uniqueness for such Stochastic Volterra equations with locally H\"older continuous diffusion coefficients are established in~\cite{PromelScheffels2023b}.
	 Our setting, described by \eqref{eq:Volterrameanrevert}-\eqref{eq:sigma-definition}, corresponds to the critical case \(\gamma=\frac12\). 
	
	\begin{Proposition}[$L^p$-Moment control and $\gamma$-H\"older $L^2$-Contraction]\label{prop:Momentctrl_contraction_Hol} 
		
		Consider the time-inhomogeneous affine Volterra equation~\eqref{eq:Volterrameanrevert}--\eqref{eq:ConstMeanVar} with locally \(\gamma\)-H\"older continuous state-dependent diffusion coefficients \(\sigma\).
		Assume assumptions~\ref{assump:kernelVolterra},~\ref{ass:Mean_revertF} and~\ref{ass:MasterEq} hold for some fixed $\lambda,c >0$ so that Equation $\textit{($E_{\lambda, c}$)}$ is in force. 		
		Then, one has:

	\smallskip
	     \noindent(a)({\em $L^p$-moments}.) Let $p\!\in [2,+\infty)$ and $(X_t)_{t\ge 0}$ the solution to the Volterra equation starting  from any random variable $X_0\!\in L^p(\P)$. 
	     Define \( \rho_{p,r} := 4p\,c\, \gamma\big(\frac{[\sigma]^2_{\text{H\"ol}}}{r}\big)^\frac{1}{\gamma} \) for sufficiently large constant \(r>0\) such that \(\rho_{p,r} \in(0,1)\). For $\gamma=1$, take $r=1$ and assume
	     that $\rho_{p,1}=4pc[\sigma]_{\mathrm{H\ddot{o}l}}^2<1$.
	      Then 	
			{\small
			\begin{align*}&\sup_{t\ge 0}\Big\|  |X_t-m_0|\Big\|_p\le \inf_{\epsilon\in (0, \frac{1}{\rho_{p,r}}-1)} \Big[ (1+\frac1\epsilon)\Big\|  |X_0-m_0|\Big\|_p  \\
				&\vee \frac{\Big(4p\,c\,[\sigma]^2_{\text{H\"ol}}(1+\epsilon)|\sigma(m_0)|^2+|\sigma(m_0)|\big(4p\,c\,[\sigma]^2_{\text{H\"ol}}(1+\epsilon)|\sigma(m_0)|^2+(1-\gamma)r^{\frac{1}{1-\gamma}}(1-4p\,c\,[\sigma]^2_{\text{H\"ol}}(1+\epsilon))\big)^{\frac12}\Big)}{\gamma\big(\frac{[\sigma]^2_{\text{H\"ol}}}{r}\big)^\frac{1}{\gamma}(1-4p\,c\,[\sigma]^2_{\text{H\"ol}}(1+\epsilon))\big((1-\gamma)r^{\frac{1}{1-\gamma}}+4p\,c\,[\sigma]^2_{\text{H\"ol}}(1+\epsilon)|\sigma(m_0)|^2\big)^{\frac12}} \Big] < +\infty.
			\end{align*}
			}
			
		\smallskip
		\noindent(b) ($\gamma$-H\"older $L^2$-Contraction) Let $X_0,X'_0\!\in L^2({\mathbb P})$ and consider the solutions to the SVIE, denoted $(X_t)_{t\ge 0}$ and $(X'_t)_{t\ge 0}$ starting from $X_0$ and $X'_0$ respectively and constructed with the same BM $(W_t)_{t\ge 0}$. 
		For  $c\!\in \big(0,\frac{1}{[\sigma]^2_{\text{H\"ol}}}\big)$, set $ \rho:=c[\sigma]^2_{\text{H\"ol}} $.  There exists a continuous non-negative function \( \varphi_{\infty}: \mathbb{R}^+ \to [0, 1] \), only depending on \( \lambda, c, K \) 
			such that 
			{\small
				$\lim_{t \to +\infty} \varphi_\infty(t) = \left\{
				\begin{array}{ll}
					& \overline{\ell}_\varphi^\star := \frac{\ell_\infty^2}{1-\rho(1-\ell_\infty^2)} \text{ if } \gamma = 1 \\
					&\underline{\ell}_{\varphi}^\star
					=\frac14
					\big(
					\rho(1-\ell_\infty^2)+
					(\rho^2(1-\ell_\infty^2)^2+4\ell_\infty^2)^{\frac12}
					\big)^2 \text{ if } \gamma = \frac12\\
					& \ell_\varphi^\star \in ]\overline{\ell}_\varphi^\star, \underline{\ell}_{\varphi}^\star[ \text{ if } \gamma \in ]\frac12,1[ 
				\end{array}, \right .$
			}
			\( \varphi_\infty(0) = 1 \) and
			\vspace{-.3cm}
			\begin{equation*}
			\forall t \geq 0, \; \E\, \big(| X_t - X^\prime_t|^2\big) \leq \varphi_\infty(t) \Psi \big(\E\, \big(| X_0 - X^\prime_0|^{2}\big)\big) , \; \text{where} \, \Psi : x \to \left\{
			\begin{array}{ll}
				& x \text{ if } \gamma = 1 \\
				& 1 \vee x \text{ if } \gamma \in [\frac12,1[ 
			\end{array} .\right .
			\end{equation*}
			In particular, if \(\ell_\infty=0\) (that is \(m_0:=(1-a)\frac{\mu_\infty}{\lambda}\)), then $\lim_{t \to +\infty} \varphi_\infty(t) = \left\{
			\begin{array}{ll}
				&  \rho^{\frac{1}{1-\gamma}} \text{ if } \gamma \in [\frac12,1[ \\
				& 0 \text{ if } \gamma = 1
			\end{array} .\right .$
	\end{Proposition}
	\noindent As a sanity check, we observe that the case \(\gamma=1\) recovers the result of \cite[Proposition 4.3]{EGnabeyeu2025}. In this particular Lipschitz case, the $L^2$-confluence, can be extended in $L^p$ using a banach fixed point argument (see, e.g.~\cite[Proposition 4.4]{EGnabeyeu2025}). The following Corollary is a straightforward consequence of Proposition~\ref{prop:Momentctrl_contraction_Hol}. 
	\begin{Corollary}[First asymptotics]\label{corol:Momentctrl_contraction_Hol} 
		Under the same assumptions as Proposition~\ref{prop:Momentctrl_contraction_Hol}~$(b)$, We have the following claims:
		\begin{itemize}
			\item[(a)] ($\mathcal{W}_2$-Contraction): With the $2$-Wasserstein distance between \(X\) and \(X^\prime\), for every integer $N \geq 1$ and every \(0 \leq t_1 < t_2 < \cdots < t_N < +\infty,\)
			we have
			\begin{equation}
				\forall t \geq 0, \; W_2\!\left(
				[X_{t+t_1},\ldots,X_{t+t_N}],
				[X'_{t+t_1},\ldots,X'_{t+t_N}]
				\right)
				\leq
				\sqrt{N}\,\sqrt{\widetilde{\varphi}_\infty(t)}\,\Psi \Big( W_2([X_0'], [X_0]) \Big)
			\end{equation}
			where \(\widetilde{\varphi}_\infty(t)
			=
			\sup_{u\geq t}\varphi_\infty(u),\)
			and $\varphi_\infty$ given by Proposition~\ref{prop:Momentctrl_contraction_Hol}~$(b)$.
			\item[(b)] (Lipschitz $L^2$-Confluence): In the setting \(\gamma = 1\), in particular, if \(\ell_\infty=0\), we have:
			\begin{itemize}
				\item  if \( X \) has a fake stationary regime of type I,  \( \mathbb{E} X_t' \to m_0 \), \( \text{Var}(X_t') \to v_0 \) as \(t \to +\infty\).
				
				\noindent And more generally finite-dimensional \( W_2 \)-convergence. 
				\item  In case if $X$ has a fake stationary regime of type II, its marginal distribution is unique.
			\end{itemize}
		\end{itemize}
	\end{Corollary}
	\noindent {\bf Remark on Lipschitz $L^2$-Confluence:}
		For the diffusion coefficient \( \sigma^2 \) given in~\eqref{eq:sigma-definition}, assuming moreover that \(\kappa_0>0\), \( \sigma^2 \) is 
		uniformly elliptic i.e.
	\(
	\exists\, \bar{\sigma}_0 > 0,\;\text{such that}\;\ \forall x \in \mathbb{R}_+,\ \sigma^2 (x) \geq \bar{\sigma}_0^2
	\), then \(\sigma\) is Lipschitz on \(I:=(0,+\infty)\) since
	
	\centerline{$ \forall x, y \in (0,+\infty), 
		|\sigma(x) - \sigma(y)| = \frac{|\kappa_1| |x-y|}{\sigma(x) + \sigma(y)} \leq \frac{|\kappa_1|}{2 \bar{\sigma}_0} |x-y|. $}
	\vspace{0.3cm}
	\noindent Consequently, if \(\kappa_0 + \kappa_1\;X_0 > 0\) and \( \mu, \lambda, \kappa_0, \kappa_1 \) are such that \( \kappa_0 + \kappa_1\;X_t > 0 \) for all \( t \) i.e.  any continuous solution \(X\) of equation~\eqref{eq:Volterrameanrevert} satisfies $\mathbb{P}( \sigma^2(X_t) > 0, \forall t \geq 0) = 1, \; \P-a.s $ , then by setting \(\kappa:=\left(\frac{\kappa_1}{2 \bar{\sigma}_0}\right)^2 > 0\), we have: 
	\begin{equation}\label{eq:lipschitz} 
		\forall t \geq 0 \quad	\E \left[ |\sigma(X_{t}) - \sigma(X^\prime_{t})|^2 \right]
		\leq \kappa \E \left[ | X_{t} - X^\prime_{t} |^2 \right].
	\end{equation}
	So that if $c$ is taken such that \(c \in (0, \frac1\kappa)\), the Proposition \ref{prop:Momentctrl_contraction_Hol} above applies as if \(\gamma=1\), in which case, there exists a continuous non-negative function \( \varphi_{\infty}^{\lambda,c, K}=: \varphi_{\infty}: \mathbb{R}^+ \to [0, 1] \), satisfying \( \varphi_\infty(0) = 1 \), \(  \lim_{t \to +\infty} \varphi_\infty(t) =  \frac{\ell_\infty^2}{1-\rho(1-\ell_\infty^2)} \), only depending on \( \lambda, c \) and \( K \), such that :
	\[
	\forall t \geq 0, \; \E\, \Big(\Big| X_t - X^\prime_t\Big|^2\Big) \leq \varphi_\infty(t) \E\, \Big(\Big| X_0 - X^\prime_0\Big|^{2}\Big) \;\text{or}\;  W_2([X_t'], [X_t]) \leq \varphi_\infty(t)^{1/2} W_2([X_0'], [X_0]).
	\]
	Consequently, whenever \(\ell_\infty=0\), we have \(\lim_{t \to +\infty} \varphi_\infty(t) = 0\) so that the confluence in Corollary \ref{corol:Momentctrl_contraction_Hol}~$(b)$ holds, that is  \( \mathbb{E} X_t' \to m_0 \), \( \text{Var}(X_t') \to v_0 \) as \(t \to +\infty\) (see, e.g. Figure~\ref{fig:confluent} further on): Thus, the process \(X^\prime\) mixes: as time increases, the
	random variable $X^\prime_t$ progressively forgets its initial mean and variance, and approaches a limit of a fake stationarity regime. 
	 Moreover, if X has a fake stationary regime of type II, its marginal distribution is unique.
	
	\subsection{Existence of limiting distributions}\label{subsect:LimitD}
	To establish the existence of limiting distributions for the inhomogeneous affine Volterra process~\eqref{eq:Volterrameanrevert}-\eqref{eq:Volterrameanrevert2} with diffusion coefficient~\ref{eq:sigma-definition}, it is sufficient to prove the convergence of its Fourier-Laplace transform, i.e. to show that the limit as \( t \to \infty \) in~\eqref{eq: MeasureFLplaceZero} and~\eqref{eq: MeasureFLplaceZero II} of Theorem~\ref{T:VolSqrt} (b) exists.  
    This requires that the function \(\psi\) obtained from~\eqref{eq:measureFLplce} is globally integrable in time; for example, it suffices that \(\psi \in L^1(\mathbb{R}_+; \mathbb{C}) \cap L^2(\mathbb{R}_+; \mathbb{C})\), a condition that has been established in Corollary~\ref{corol:riccati-extension}.  
    As a first step, we have the convergence of the Fourier-Laplace transform given below.
	
	\begin{Proposition}\label{prop: convergence FT}
	Let $X$ be the time-inhomogeneous affine Volterra Equation~\eqref{eq:Volterrameanrevert}-\eqref{eq:Volterrameanrevert2} with diffusion coefficient \(\sigma\) given by~\eqref{eq:sigma-definition} and
	  initial random state $X_0$. Suppose that Assumptions~\ref{ass:resolvent},~\ref{ass:int_holregul} and~\ref{ass:Mean_revertF} are in force
		and denote $\Xi_0 =
		\ell_\infty X_0  + (1-a) \frac{\mu_\infty}{\lambda}$. 
		Let $m \in \mathbb{M}_{\mathrm{lf}}^{-}$ be such that $|m|(\mathbb{R}_+) < \infty$. Then the Riccati-Volterra equation~\eqref{eq:measureFLplce2}-\eqref{eq:FVarsigma_infty} defined on \(\R_+\) has a unique global solution \(\psi =\psi(\cdot,m) \in L^1(\R_+; \mathbb{C}_-) \cap L^2(\R_+; \mathbb{C}_-)\) and 
		\begin{align}
			&\ \lim_{t \to \infty}
			\mathbb{E}_{}\Bigg[ \exp\Big( \int_0^t X_{t-s} \,m(ds) \Big) \Big| X_0 \Bigg]
			= \exp\left[  \Xi_0 \;m(\R_+) +  \frac{\varsigma^2_\infty}{2}\sigma^2\left( \Xi_0\right)\int_0^{\infty} \psi(s,m)^2 ds \right] \label{eq: FT convergence}
			\\ &\hspace{.1cm}= \exp\left[
			 m(\R_+) \ell_\infty X_0  + \frac{\varsigma^2_\infty}{2} \Big( \int_0^\infty \psi^2(s,m) \, ds \Big)\,\big(\kappa_0 + \kappa_1 \ell_\infty X_0 \big) +\Big( \int_0^{\infty}  \psi(s,m) \, ds\Big) \mu_\infty  \right] \label{eq: FT convergence II}
		\end{align}
		$\P\text{-a.s.}$, where $\varsigma^2_\infty$ is defined by Equation~\eqref{eq:FVarsigma_infty} and $\psi(\cdot,m) $ denoting the unique solution of~\eqref{eq:measureFLplce2}. 
	\end{Proposition}
	\begin{Example}[Limiting distributions in the stabilized Regime.]\label{Ex:longRunStabil}
		Let $X$ be the (stabilized) inhomogeneous affine Volterra equation~\eqref{eq:ConstMeanVar} (conditions~\eqref{eq:CondMean} and~\eqref{eq:MasterEq} are
		in force) with constant $\lambda, c > 0$, $m_0, \mu_\infty \in \mathbb{R}^*_+$ and diffusion coefficient \(\sigma\) given by~\eqref{eq:sigma-definition}. 
		 Noticing that in this regime, \( \varsigma^2_\infty:= \lim_{t \to +\infty} \varsigma^2(t) =\frac{c\lambda^2 (1-\ell_\infty^2)}{\|f_{\lambda}\|^2_{L^2(\text{Leb}_1 )}}\) owing to \cite[Lemma 3.9 ]{EGnabeyeu2025}, the identities ~\eqref{eq: FT convergence} and ~\eqref{eq: FT convergence II} become:
		{\small 
			\begin{align}
				&\ \lim_{t \to \infty}
				\mathbb{E}_{}\Big[ \exp\Big( \int_0^t X_{t-s} \,m(ds) \Big) \Big| X_0 \Big]= \exp\Big[ \Xi_0\, m(\mathbb R_+) +  \frac{1}{2} \frac{c\lambda^2 (1-\ell_\infty^2)}{\|f_{\lambda}\|^2_{L^2(\text{Leb}_1 )}}\sigma^2\left(\Xi_0\right)\int_0^{\infty} \psi^2(s,m) ds \Big] \label{eq: FT convergence_}
				\\ &\hspace{2cm}= \exp\Big[  m(\mathbb R_+) \, \ell_\infty X_0 +\Big( \int_0^{\infty}  \psi(s,m) \, ds\Big) \mu_\infty + \frac{1}{2}  \frac{c\lambda^2 (1-\ell_\infty^2)}{\|f_{\lambda}\|^2_{L^2(\text{Leb}_1 )}} \big(\kappa_0 + \kappa_1 \ell_\infty X_0 \big) \int_0^\infty \psi^2(s,m) \, ds \Big]. \notag
			\end{align}
		}
	\end{Example}
	For clarity and conciseness, the proof of proposition~\ref{prop: convergence FT} is deferred to Appendix \ref{app:lemmata}, where the main technical results are presented.
	From the convergence of the Fourier-Laplace transform we can now deduce convergence towards limiting distributions.
	The following is our main result on limiting distributions for the time-inhomogeneous affine Volterra process~\eqref{eq:Volterrameanrevert}, 
	 which, in contrast to the classical case, also depend on the initial state of the process.
	\begin{Theorem}[Limiting Distribution]\label{Theorem: limiting distribution}
		Let $X$ be the nonnegative time-inhomogeneous affine Volterra Equation~\eqref{eq:Volterrameanrevert}-\eqref{eq:Volterrameanrevert2} with diffusion coefficient \(\sigma\) given by~\eqref{eq:sigma-definition} and
		initial random state $X_0$. Suppose that Assumptions~\ref{ass:resolvent}, \ref{ass:int_holregul} and~\ref{ass:Mean_revertF} are in force
		and denote $\Xi_0 =
		\ell_\infty X_0  + (1-a) \frac{\mu_\infty}{\lambda}$. 
		Then, as $t\to\infty$, the conditional law of the random variable $X_t$ given
		$X_0$ converges weakly to a limiting random probability measure $\pi_{X_0}$, whose
		Fourier-Laplace transform is, for $u\in\mathbb{C}_-$, 
		\begin{align}
			&\ \int_{\R_+} \exp\left( u\,x \right)\pi_{X_0}(dx)= \exp\left[ u \,\Xi_0 +  \frac{\varsigma^2_\infty}{2}\sigma^2\left(\Xi_0\right)\int_0^{\infty} \psi^2(s,u\,\delta_0) ds \right]\label{eq:CondLawLaplace1}
			\\ &\ = \exp\left[u\, \ell_\infty X_0 +\left( \int_0^{\infty}  \psi(s,u\,\delta_0) \, ds\right) \mu_\infty + \frac{\varsigma^2_\infty}{2}  \big(\kappa_0 + \kappa_1 \ell_\infty X_0 \big) \int_0^\infty \psi^2(s,u\,\delta_0) \, ds \right]. \label{eq:CondLawLaplace2}
		\end{align}
		$\P\text{-a.s.}$, where $\varsigma^2_\infty$ is defined by Equation~\eqref{eq:FVarsigma_infty}. 
		Moreover, $\pi_{X_0}$ has a finite first moment $\P\text{-a.s.}$, with $
		\int_{\mathbb R_+}x\,\pi_{X_0}(dx)
		=\Xi_0
		,\; \text{a.s.}$
	\end{Theorem}
	\noindent {\bf Proof of Theorem~\ref{Theorem: limiting distribution}.}
	Consider $u \in \mathbb{C}_-$ and take \(m(ds) = u\,\delta_0(ds)\). Proposition \ref{prop: convergence FT} gives the pointwise
	convergence of the corresponding characteristic functions. That is $\P-a.s.$,
	\begin{align} \notag
		&\ \lim_{t \to \infty}
		\mathbb{E}_{} \left[
		\exp\left(u X_t \right) \big|X_0	\right]=:\Phi_{X_0}(u)
		\\ &\hspace{1.6cm}= \exp\left[ u \left(
		\ell_\infty X_0  + (1-a) \frac{\mu_\infty}{\lambda}\right) +  \frac{\varsigma^2_\infty}{2}\sigma^2\left(
		\ell_\infty X_0  + (1-a) \frac{\mu_\infty}{\lambda}\right)\int_0^{\infty} \psi(s,u\,\delta_0)^2 ds \right] \label{eq:CongceP2}
	\end{align}
		For \(t\geq0\), let $P_t(X_0,dx):=\mathbb P(X_t\in dx\mid X_0)$
		denote the conditional law of \(X_t\) given \(X_0\), i.e. $P_t(\omega,dx):=\mathcal{L}(X_t\mid X_0)(\omega,dx)$. Since \(X_t\geq0\),
		\(P_t(\omega,\cdot)\in\mathcal{P}(\mathbb{R}_+)\). Moreover,
		\vspace{-.1cm}
		\[
		\int_{\mathbb{R}_+}x\,P_t(\omega,dx)
		=
		\mathbb{E}[X_t\mid X_0](\omega)
		\longrightarrow
		\Xi_0(\omega),
		\qquad t\to\infty,
		\]
		almost surely by Lemma~\ref{lm:firstmoment}. By a standard countable dense-set argument, let $\Omega_0\subset\Omega$ with  $\P(\Omega_0)=1$, where $\Omega_0$ is the set of event on which the above
		convergence and~\eqref{eq:CongceP2}
		hold.  
		Fix \(\omega \in \Omega_0\). Then, for some \(t_0=t_0(\omega)\), $\sup_{t\geq t_0}
		\int_{\mathbb{R}_+}x\,P_t(\omega,dx)<\infty.$
		Hence, by Markov's inequality, $\sup_{t\geq t_0}P_t(\omega,[R,\infty))
		\leq
		\frac{1}{R}
		\sup_{t\geq t_0}
		\int_{\mathbb{R}_+}x\,P_t(\omega,dx)
		\longrightarrow0,
		\; R\to\infty.$
		Thus, \(\{P_t(\omega,\cdot):t\geq t_0\}\) is tight.
		Let \((t_n)_{n\geq1}\) be any sequence such that \(t_n\to\infty\).
		By Prokhorov's theorem, there exist a subsequence
		\((t_{n_k})_{k\geq1}\) and a probability measure
		\(\pi_\omega\in\mathcal{P}(\mathbb{R}_+)\) such that $P_{t_{n_k}}(\omega,\cdot)\Longrightarrow\pi_\omega.$
		Now, for \(u\in\mathbb{C}_-\), the mapping \(x\mapsto e^{ux}\) is
		bounded and continuous on \(\mathbb{R}_+\), and therefore, by weak convergence,
		\vspace{-.3cm}
		\[
		\int_{\mathbb{R}_+}e^{ux}\,\pi_\omega(dx)
		=
		\lim_{k\to\infty}
		\int_{\mathbb{R}_+}e^{ux}\,P_{t_{n_k}}(\omega,dx)
		=
		\lim_{k\to\infty}
		\mathbb{E}\!\left[e^{uX_{t_{n_k}}}\mid X_0\right](\omega)= \Phi_{X_0(\omega)}(u).
		\]
		owing to~\eqref{eq:CongceP2}. Consequently, any weak accumulation point has the same Fourier-Laplace transform.
		In particular, taking \(u=i\xi\), \(\xi\in\mathbb{R}\), shows that
		the characteristic function of any weak accumulation point is uniquely
		determined by \(X_0(\omega)\). Hence all weak accumulation points
		coincide. Since the family is tight, it follows that
		\begin{equation}\label{eq:Lim1}
		P_t(\omega,\cdot)\Longrightarrow\pi_{X_0(\omega)},
		\qquad t\to\infty,
	    \end{equation}
		for a unique probability measure \(\pi_{X_0(\omega)}\) on
		\(\mathbb{R}_+\). Consequently, $\mathcal{L}(X_t\mid X_0)\Longrightarrow\pi_{X_0}
		$ almost surely in $\mathcal{P}(\mathbb{R}_+)$.
		Finally, let \(u\in\mathbb C_-\). Since \(X_t\geq0\), the mapping
		\(x\mapsto e^{ux}\) is bounded and continuous on \(\mathbb R_+\).
		Hence, by the weak convergence established above,
		$\mathbb E\!\left[e^{uX_t}\mid X_0\right](\omega)
		=
		\int_{\mathbb R_+}e^{ux}\,P_t(\omega,dx)
		\longrightarrow
		\int_{\mathbb R_+}e^{ux}\,\pi_{X_0(\omega)}(dx).$
		On the other hand, \eqref{eq:CongceP2} still holds and yields $\mathbb E\!\left[e^{uX_t}\mid X_0\right](\omega)\to \Phi_{X_0(\omega)}(u)$.
		Therefore, by uniqueness of the limit,
		\vspace{-.3cm}
		\begin{equation}\label{eq:Lim2}
		\int_{\mathbb R_+}e^{ux}\,\pi_{X_0(\omega)}(dx)
		=\Phi_{X_0(\omega)}(u),
		\quad u\in\mathbb C_-.
		\vspace{-.2cm}
		\end{equation}
		 Since \(\omega\in\Omega_0\) was arbitrary, denoting the resulting
		random probability measure by \(\pi_{X_0}\), this yields~\eqref{eq:CondLawLaplace1} almost surely.
		Using~\eqref{eq: FT convergence II}, this may equivalently be written as.
		\[
		\int_{\mathbb{R}_+}e^{ux}\,\pi_{X_0}(dx)
		=\exp\Big[u \,\ell_\infty X_0   + \frac{\varsigma^2_\infty}{2}  \Big(\int_0^\infty \psi^2(s,u\,\delta_0) \, ds\Big) \big(\kappa_0 + \kappa_1 \ell_\infty X_0 \big) +\Big( \int_0^{\infty}  \psi(s,u\,\delta_0) \, ds\Big) \mu_\infty \Big].
		\]
		which yields~\eqref{eq:CondLawLaplace2}. It remains to identify the first moment of the limiting distribution. Let $m=u\delta_0$ with $, u\in \R_-$. We note that under condition~\eqref{eq:contf_lambda}, estimate~\eqref{eq:LpBounds} hold with $|m|(\mathbb{R}_+) = |u|$ and \(p=2\), showing that the right-hand side of~\eqref{eq:CondLawLaplace1} is continuous at $u = 0$. For $u<0$, we obtain
		$\Big|
		\int_0^\infty\psi^2(s,u\delta_0)\,ds
		\Big|
		\leq
		\|\psi(\cdot,u\delta_0)\|_{L^2(\mathbb R_+)}^2
		=O(|u|^2),
		\; u\uparrow0.$
		Therefore, for every \(\omega\in\Omega_0\),
		\vspace{-.2cm}
		\[
		\int_{\mathbb R_+}e^{ux}\,\pi_{X_0(\omega)}(dx)
		=
		\exp\left[
		-|u|\Xi_0(\omega)+O(u^2)
		\right]
		=
		1-|u|\Xi_0(\omega)+o(|u|).
		\vspace{-.2cm}
		\]
		On the other hand, by monotone convergence,
		\vspace{-.2cm}
		\[
		\int_{\mathbb R_+}x\,\pi_{X_0(\omega)}(dx)
		=
		\lim_{u\uparrow0}
		\int_{\mathbb R_+}\frac{1-e^{ux}}{|u|}\,\pi_{X_0(\omega)}(dx)
		=
		\lim_{u\uparrow0}
		\frac{1-\int_{\mathbb R_+}e^{ux}\,\pi_{X_0(\omega)}(dx)}{|u|}.
		\]
		Combining the last two displays gives $\int_{\mathbb R_+}x\,\pi_{X_0(\omega)}(dx)=\Xi_0(\omega)<\infty$, which yields the result.
	This completes the proof. \qed

	\smallskip
	\noindent {\bf Remark:}
	1. Note that as \(\Phi_{X_0(\omega)}(0)=1\), Equations~\eqref{eq:Lim1} and~\eqref{eq:Lim2} could also be obtained by Lévy's continuity theorem (see, e.g.,~\cite{friesen2022volterra}) for \(u:=i \xi, \, \xi \in \R\), provided that $\lim_{\xi\to0}\Phi_{X_0}(i\xi)=1$, which is true by the continuity of
	\(u\mapsto\psi(\cdot,u\delta_0)\) at \(u=0\) in
	\(L^2(\mathbb R_+)\) (see, e.g.,~\eqref{eq:LpBounds}).
	
	2. In the classical (Markovian) setting, one would typically use Feller
	semigroup theory to identify a limiting distribution as an invariant
	probability measure (stationary distribution). More precisely, if \(P_t\) denotes the transition
	semigroup and \(\pi\) is a limiting distribution, one seeks to verify
	that
	\vspace{-.2cm}
	\[
	\int_{\mathbb{R}_+} P_t f(x)\,\pi(dx)
	=
	\int_{\mathbb{R}_+} f(x)\,\pi(dx),
	\qquad t\geq0,
	\]
	for every bounded continuous function \(f\). Equivalently, taking
	\(f(x)=e^{iux}\), \(u\in\mathbb{R}\),
	\[
	\mathbb{E}_{\pi}\left[\exp\left(iuX_t\right)\right]
	=:\int_{\mathbb{R}_+}
	\mathbb{E}_{x}\left[\exp\left(iuX_t\right)\right]
	\,\pi(dx)
	=
	\int_{\mathbb{R}_+}\exp(iux)\,\pi(dx).
	\]
		This shows that \(X_0\sim\pi\) implies \(X_t\sim\pi\) for every
	\(t\geq0\), so that \(\pi\) is an invariant distribution, stationary for the process \((X_t)_{t\geq0}\) in the time-homogeneous Markov setting.
	However, in our setting, the process \(X\) does not possess the Markov property. Therefore,
	the usual invariant-measure argument based on a Markov semigroup does
	not directly apply. Instead, we follow the framework developed in
	\cite{friesen2022volterra}, which is specifically designed to handle
	non-Markovian dynamics, including Volterra-type processes.
	\subsection{Asymptotics: Long run functional weak behaviour and Stationary process}\label{subsect:longRun}
	In the next step, we construct the associated stationary process.
	To this end, we recall that for a real-valued stochastic process $X$ on $[0, T]$, the Kolmogorov continuity theorem states that if for some constants $p, a, C > 0$,
	\begin{equation}
		\mathbb{E}\left[ \left| X_t - X_s \right|^p \right]
		\leq C \cdot |t - s|^{1 + a}, 
	\end{equation}
	uniformly in $0 \leq s, t \leq T$, then the process has a $\theta$-H\"older continuous modification for all $0 < \theta < a/p$.
	
	To prove the existence of a stationary process associated with a limiting distribution, we cannot rely, as is typically done in the classical literature, on \textit{Feller semigroup} theory or standard \textit{Markovian arguments}. Instead, we adopt the alternative approach developed in~\cite{friesen2022volterra}, which is based on the exponential-affine transform formula presented in Theorem~\ref{T:VolSqrt}.
	Thanks to that result, we obtain explicit control over the process's finite-dimensional distributions.
	Consequently, we prove in this section that, as $u \to \infty$, the shifted process $(X^u_t)_{t \geq 0}$ defined by $X^u_t := X_{t+u}$ converges in law to a continuous process $X^{\infty}$. 
	From Theorem~\ref{Theorem: limiting distribution}, it follows that for any $t \geq 0$, the conditional law of $X_t$ given $X_0$, 
	 $\mathcal L(X_t\mid X_0) \longrightarrow \pi_{X_0}$ weakly as $t \to \infty$. Thus, $ \mathcal L(X_t^\infty\mid X_0) $ has distribution $\pi_{X_0}$ for each $t \geq 0$, which is therefore the desired stationary solution associated with the limiting conditional law.
	\begin{Theorem}[Long run theorem: Functional weak asymptotics and Stationary Process]\label{Theorem: stationary process}	
		Suppose that Assumptions~\ref{ass:resolvent}, \ref{ass:int_holregul}, \ref{ass:MasterEq} and~\ref{ass:Mean_revertF} are in force.
		Let $X$ be the (stabilized) nonnegative  time-inhomogeneous affine Volterra  equation~\eqref{eq:ConstMeanVar} (conditions~\eqref{eq:CondMean} and~\eqref{eq:MasterEq} are
		in force) with constant $\lambda, c > 0$, $m_0, \mu_\infty \in \mathbb{R}^*_+$ and diffusion coefficient \(\sigma\) given by~\eqref{eq:sigma-definition}. Let $X_0$ be the
		initial random state
		and denote $\Xi_0 =
		\ell_\infty X_0  + (1-a) \frac{\mu_\infty}{\lambda}$. 
		 Then the following 
		  hold:
		
		\medskip
		\noindent {(a)} Then, the family of shifted processes \( (X_{t+u})_{u \geq 0} \) is C-tight, uniformly integrable, and square uniformly integrable for \( p > 2 \) as \( t \to +\infty \). For any limiting distribution \( P \) on \( \Omega_0 := \mathcal{C}(\mathbb{R}_+, \mathbb{R}) \), the canonical process \( Y_t(\omega) = \omega(t) \) has a \( \left( \vartheta \wedge \frac{\beta-1}{2\beta} - \frac{1}{p} - \eta \right) \)-H\"older pathwise continuous \( P \)-modification for sufficiently small \( \eta > 0 \). 
		There exists a stationary process \( X^{\infty} \) with continuous sample paths such that
		\begin{equation}\label{eq:joint_convergence}
			(X_0,X_{u+\cdot})
			\Rightarrow (X_0,X^\infty)
			\quad\text{weakly in }
			\mathbb R_+\times\mathcal C(\mathbb R_+,\mathbb R_+),
			\qquad u\to\infty.
		\end{equation}
		
		\medskip
		\noindent {(b)} For every $t\geq 0$, the stationary process $X^{\infty}$ satisfies $\pi_{X_0}=\mathcal L(X_t^\infty\mid X_0),\, t\geq0,\; \mathbb P\text{-a.s.}$.  
		Moreover, its first moment is given by
		\begin{equation}\label{eq:FirstMmt}
	\forall\, t\geq0,\;	\E[X_t^{\infty}] =\mathbb{E}\big[\Xi_0\big]= \ell_\infty \mathbb{E}[X_0]  + (1-a) \frac{\mu_\infty}{\lambda}=\ell_\infty \mathbb{E}[X_0]  + m_0(1-\ell_\infty).
		\end{equation}
		while its autocovariance function is \(\text{Cov}(X_{t_1}^\infty, X_{t_2}^\infty) :=C_{f_{\lambda}}(t_1,t_2)\), for 
		\(0\leq t_1 \leq t_2 \), given by 
		\begin{equation}\label{eq:CovfunclongRun}
			\mathrm{cov}(X_{t_2}^{\infty}, X_{t_1}^{\infty}) =\ell_\infty^2 {\rm Var}(X_0) + \frac{c (1-\ell_\infty^2)}{ \|f_{\lambda}\|^2_{L^2(\text{Leb}_1 )}} \sigma^2\big( \mathbb{E}\big[\Xi_0\big] \big)\int_0^{+\infty}  f_{\lambda}(t_2-t_1+u)f_{\lambda}(u)du.
		\end{equation}
		so that  the process \((X^{\infty}_t)_{t \geq 0}\) is at least weak \( L^2 \)-stationary, square-integrable process with constant 
		 mean \( \ell_\infty \mathbb{E}[X_0]  + (1-a) \frac{\mu_\infty}{\lambda} \) and covariance function \( C_{f_\lambda}(s,t) \), for \( s, t \geq 0 \).
		
		\medskip
		\noindent {(c)} 
		 The finite dimensional distributions of $X^{\infty}$, conditional on \(X_0\) are determined for $ n \in \mathbb{N}, \; u_{1},\dots,u_{n} \in \mathbb{C}_{-}$ and $0\leq t_{1}<\dots<t_{n}$ by:
		 {\small 
		 	\begin{equation}\label{eq:FD_distr}
		 	\ \mathbb{E}_{}\Bigg[ \exp\Bigg( \sum_{i=1}^{n} u_{i}  X_{t_{i}}^{\infty} \Bigg) \Big| X_0\Bigg] =\exp\left[ \Xi_0 \sum_{i=1}^{n}  u_{i} +  \frac{c\lambda^2 (1-\ell_\infty^2)}{2 \|f_{\lambda}\|^2_{L^2(\text{Leb}_1 )}}\sigma^2\left(\Xi_0\right)\int_0^{\infty} \psi(s)^2 ds \right].
		 	\end{equation}
		 }
		 \noindent $\P\text{-a.s.}$, where
		 $\psi(\cdot) = \psi(\cdot, m_{t_1,\dots, t_n})$ denotes the unique solution of \eqref{eq:measureFLplce2} in \(\mathbb{R}_+\)
		 with  $m_{t_{1},\dots,t_{n}}(ds) \;=\; \sum_{j=1}^{n} u_{j} \, \delta_{\,t_{n}-t_{j}}(ds).$
	\end{Theorem}
	\noindent {\bf Remark:}(Functional asymptotics in the fake stationarity Regime.)\label{Thm:longRun} 
			\noindent Under the same conditions as in Theorem \ref{Theorem: stationary process},
			if the solution \((X_t)_{t \geq 0}\) of the Volterra equation~\eqref{eq:ConstMeanVar} has a fake stationary regime of type I, starting from a random variable \(X_0\) with mean $m_0$
			and variance \(v_0:=c\sigma^2\big(m_0\big)\), then, the family of shifted processes \( X_{t+\cdot}, t \geq 0 \), is \( C \)-tight as \( t \to +\infty \) and its (functional) limiting distributions are all at least weak \( L^2 \)-stationary
			processes with mean \(m_0\) and covariance function \( C_\infty \) given by:
			{\small
				\begin{equation}\label{eq:weak_cov}
					\forall t_1, t_2 \geq 0 \quad \text{with} \quad t_1 \leq t_2 \quad  C_\infty(t_1, t_2) =c \sigma^2\left(m_0\right)\Big(\ell^2_\infty + \frac{ (1-\ell_\infty^2) }{\|f_{\lambda}\|^2_{L^2(\text{Leb}_1 )}}\int_0^{\infty} f_{\lambda}(t_2 - t_1 + u) f_{\lambda}(u) \, du\Big).
				\end{equation}
			}
			\noindent Moreover, the finite dimensional limiting distributions conditional on \(X_0\) are still determined by~\eqref{eq:FD_distr} and in particular, the final distribution in the fake stationary regime does not depend on the initial distribution whenever \(\ell_\infty=0\) i.e.  \(m_0:= (1-a)\frac{\mu_\infty}{\lambda}\).
			
	\begin{Example}[Long-run behaviour of the stabilized \(\rho-\)exponential \(\alpha\)-fractional Cox-Ingersoll-Ross process]\label{Ex: LongR_VCIR 1dim}
		In this example, we briefly state our results when applied to the \textit{Volterra Cox-Ingersoll-Ross  process with a Gamma or an exponentially damped fractional integration kernel} obtained from
		\vspace{-.25cm}
		\begin{align}\label{eq: rough CIR}
			X_t &= x_0(t) + \int_0^t K(t-s) (\mu(s) -\lambda X_s) ds + \nu \int_0^t K(t-s)\varsigma_{\lambda,c}(s)\sqrt{X_s}dW_s,\quad X_0\perp\!\!\!\perp W.
		\end{align}
		where $K(t) = K_{1, \alpha, \rho}(t) = e^{-\rho t} \frac{t^{\alpha - 1}}{\Gamma(\alpha)} \mathbf{1}_{\mathbb{R}_+}(t)$,
		$\alpha \in (\frac12,1)$, $x_0 = X_0 \varphi$ with $\varphi$ in~\eqref{eq:CondMean}, $\lambda, \nu, X_0 \geq 0$, $(W_t)_{t \geq 0}$ a one-dimensional Brownian motion and \(\varsigma_{\lambda,c}\) given as the solution of \(\textit{($E_{\lambda, c}$)}\) in equation ~\eqref{eq:MasterEq} for a given $c>0$ and \(\mu:\mathbb{R}_+\to\mathbb{R}_+\).
		Letting $(X_t)_{t\geq0}$ be a weak solution of the stabilized Volterra Equation given by \eqref{eq: rough CIR} with constant $\lambda, c > 0$, $m_0, \mu_\infty \in \mathbb{R}^*_+$. We have the following claims:
		
		\medskip
		\noindent {(1).}	The shifted process $(X_{t+u})_{t\ge0}$  is \( C \)-tight as \( t \to +\infty \) and converges in law as $u \to \infty$ to a continuous stationary process $(X_t^{\infty})_{t \geq 0}$ whose finite dimensional distributions, conditional on \(X_0\) have the characteristic function
		\vspace{-.2cm}
		\begin{align*}
			\E_{}\Bigg[ \exp\Big( \sum_{i=1}^{n} u_{i}  X_{t_{i}}^{\infty} \Big) \Big|X_0 \Bigg] =\exp\left[\frac{m_0\big(\rho^{\alpha} + \lambda\big) X_0 + \mu_\infty \big(m_0- X_0\big)}{m_0\big(\rho^{\alpha} + \lambda\big)}\left(  \sum_{i=1}^{n}  u_{i}  +  \frac{\varsigma^2_\infty}{2}\nu^2\int_0^{\infty} \psi(s)^2 ds\right) \right], \; \text{where}
		\end{align*}
		$\varsigma^2_\infty := \frac{c \lambda^2 }{\|f_{\alpha, \rho, \lambda}\|^2_{L^2(\mathrm{Leb}_1)}} \frac{
			2m_0\mu_\infty\left(\rho^\alpha+\lambda\right)-\mu_\infty^2
		}{m_0^2\left(\rho^{\alpha}+\lambda\right)^2}$, $0 \leq t_1 < \dots < t_n$, $u_1,\dots, u_n \in \mathbb{C}_-$, and $\psi$ unique solution of
		\vspace{-.2cm}
		\begin{align*}
			\psi(t) &= \sum_{j=1}^n \mathbf{1}_{\{ t > t_n - t_j \}} \, K_{\alpha, \rho}(t - (t_n - t_j)) u_j  + \int_0^{t} K_{\alpha, \rho}(t - s) \left( -\lambda \psi(s) + \varsigma^2_\infty\frac{\nu^2}{2} \psi(s)^2 \right) ds.
			\vspace{-.2cm}
		\end{align*}
		\medskip
		\noindent {(2).} Moreover, the first moment and the autocovariance function of the stationary process satisfy  $\E[X_t^{\infty}] =\mathbb{E}[X_0]
			+
			\frac{\mu_\infty\bigl(m_0-\mathbb{E}[X_0]\bigr)}
			{m_0(\rho^\alpha+\lambda)}:=\bar m$ and for $t_1, t_2 \geq 0 ,\; t_1 \leq t_2 $
		\vspace{-.2cm}
		\begin{align}\label{asym_autocovariance}
			&\mathrm{Cov}(X_{t_2}^{\infty}, X_{t_1}^{\infty}) 
			= \frac{
				\left(m_0(\rho^\alpha+\lambda)-\mu_\infty\right)^2
			}{
				m_0^2(\rho^\alpha+\lambda)^2
			}\, \operatorname{Var}(X_0) + c \nu^2 \bar m
			\cdot \frac{
				\mu_\infty
				\left(
				2m_0(\rho^\alpha+\lambda)-\mu_\infty
				\right)}{m_0^2(\rho^\alpha+\lambda)^2 \|f_{\alpha, \rho, \lambda}\|^2_{L^2(\mathrm{Leb}_1)}}
			\cdot \notag \\
			& \hspace{5cm} \times e^{-\rho (t_2 - t_1)} \int_0^{+\infty} e^{-2\rho u} f_{\alpha, \lambda}(t_2 - t_1 + u) f_{\alpha, \lambda}(u) \, du.
		\end{align}
	\end{Example}
	\noindent 
	The particular form of the Fourier-Laplace transform for the limiting distribution in Theorem~\ref{Theorem: limiting distribution} and the stationary process $X^{\infty}$ in Theorem~\ref{Theorem: stationary process} give the following characterization for the independence on the initial condition $X_0 \in L^p(\P)$ for some suitable \(p>0\).
	\begin{Corollary}\label{corol:LimitDist}
		Let $X$ be the nonnegative time-inhomogeneous affine Volterra Equation~\eqref{eq:Volterrameanrevert}-\eqref{eq:Volterrameanrevert2} with diffusion coefficient \(\sigma\) given by~\eqref{eq:sigma-definition} and
		initial random state $X_0$ with $\operatorname{Var}(X_0) > 0$.
		Under the same assumptions as in Theorem~\ref{Theorem: stationary process}, the following are equivalent:
		\begin{enumerate}
			\item[(i)] The stationary process $X^{\infty}$ is independent of $X_0$;
			\item[(ii)] The random limiting probability measure $\pi_{X_0}$ is deterministic, i.e., there exists a probability measure $\pi$ such that
			\(    \pi_{X_0}=\pi,\quad \mathbb{P}\text{-a.s.};
			\)
			\item[(iii)] The random variable
			\(    \int_{\mathbb{R}_+}x\,\pi_{X_0}(dx)
			\) is constant $\, \mathbb{P}\text{-a.s.}$;
			\item[(iv)] $\ell_\infty =0$. 
		\end{enumerate}
		In particular, if $\ell_\infty=0$, the limiting distribution is deterministic and its Fourier--Laplace transform is 
		\begin{align*}
			\lim_{t \to \infty}
			\mathbb{E}\Big[ \exp\Big( \int_0^t X_{t-s} \,m(ds) \Big) \Big]
			&= \exp\Big[  (1-a) \frac{\mu_\infty}{\lambda}\;m(\R_+) +  \frac{\varsigma^2_\infty}{2}\sigma^2\Big( (1-a) \frac{\mu_\infty}{\lambda} \Big)\int_0^{\infty} \psi(s,m)^2 ds \Big]
			\\ &= \exp\Big[ \Big( \int_0^{\infty}  \psi(s,m) \, ds\Big) \mu_\infty + \frac{\kappa_0}{2} \varsigma^2_\infty \int_0^\infty \psi^2(s,m) \, ds \Big].
		\end{align*}
	\end{Corollary}
	\noindent {\bf Proof of Corollary~\ref{corol:LimitDist}.}
	\smallskip
	Since \(\pi_{X_0}\) is the 
	conditional marginal law of the stationary process \(X^\infty\), conditional on \(X_0\), 
	 it is clear that (i) implies (ii), while (ii) implies (iii).
	Suppose that (iii) holds. By the expression for the first moment of the limiting distribution given in Theorem~\ref{Theorem: limiting distribution},
	$$
	\int_{\mathbb{R}_+} y\,\pi_{X_0}(dy)
	=
	\Xi_0
	=
	\ell_\infty X_0
	+
	(1-a)\frac{\mu_\infty}{\lambda}.
	$$
	Since, by (iii), this quantity is constant $\mathbb{P}$-almost surely
	for every initial state $X_0$, we have $0=\operatorname{Var}(\Xi_0)
	=\ell_\infty^2\operatorname{Var}(X_0).$
	Since \(\operatorname{Var}(X_0)>0\), this immediately implies \(\ell_\infty=0\).
	This readily yields (iv). 
	 Finally, suppose that (iv) is satisfied. $
	\Xi_0=(1-a)\frac{\mu_\infty}{\lambda}
	$
	is deterministic and independent of the initial state $X_0$. Consequently, the Fourier--Laplace transform of $\pi_{X_0}$ given in Theorem~\ref{Theorem: limiting distribution} and the finite dimensional distributions of the stationary
	process $X^\infty$ given in Theorem~\ref{Theorem: stationary process}~$(c)$ are independent of $X_0$. Hence, by uniqueness of the Fourier--Laplace
	transform, the law of $X^\infty$ is independent of $X_0$, and therefore
	(i) holds.
	Consequently, (i)--(iv) are equivalent.
	\qed
	
	\medskip\noindent {\bf Remark:} As a consequence of that result in Example~\ref{Ex: LongR_VCIR 1dim}, we see that the stationary process $X^{\infty}$ is independent of the initial distribution of \(X_0\) if and only if $m_0 = \frac{\mu_\infty}{\rho^\alpha+\lambda}$. It boils down that \(\ell_\infty=0\) in accordance with Corollary ~\ref{corol:LimitDist}. This is the specific case of Stochastic Volterra Equations with \(\alpha-\)fractional integration kernel.
	Moreover,  for $\alpha=1$ and $\rho=0$ the autocovariance function satisfies
	\vspace{-.1cm}
	\begin{equation}
		\mathrm{cov}(X_{t_2}^{\infty}, X_{t_1}^{\infty}) = c \nu^2\,\frac{
			\mu_\infty(2m_0\lambda-\mu_\infty)
		}{
			m_0\lambda^2}\, e^{-\lambda (t_2-t_1)}.
		\vspace{-.1cm}
	\end{equation}
	which is that of the Markovian square root process with constant mean-reversion and stabilizing functions. In fact,
	If \( K = 1 \), \(\mu(t)=\mu_\infty\) for every \(t\geq0\) in the Volterra equation and \(c>0\) given, then \(R_\lambda(t) = e^{-\lambda t} \quad \text{and} \quad f_\lambda(t) = \lambda e^{-\lambda t},
	\quad \text{so that }\quad 
	\varsigma_{\lambda,c} = \sqrt{2\lambda c} \)
	and $X$ satisfy the following stochastic differential equation:
	\[dX_t = (\mu_\infty - \lambda X_t)\,dt + \varsigma_{\lambda,c} \nu \sqrt{X_t}\,dW_t, \quad t \geq 0.\]
	\noindent where $W$ is a $1$-dimensional Brownian motion on some filtered probability space $(\Omega, \mathcal{F}, \{\mathcal{F}_t\}_{t \geq 0}, \mathbb{P})$.

	\vspace{-.3cm}
	\section{Applications: The case of fake stationary Volterra CIR process.}\label{sect-Application}
	In this section, we provide a broad class of applications, starting from the fake stationary fractional-CIR process and its long run behaviour to the SPX volatility surface calibration in the model introduced in Section~\ref{subsect:FS-RH}.
	In this setting, the resulting Volterra SDE~\eqref{eq:Volterrameanrevert} admit (at least) a non-negative
	 weak solution (see, e.g., \cite{EGnabeyeuR2025}).
		Recall from Example~\ref{Ex:SolventGammaKernel} that these kernels can be
		written \( K(t) = K_{1, \alpha, 0}(t) = K_{\alpha}(t)= \frac{t^{\alpha - 1}}{\Gamma(\alpha)} \mathbf{1}_{\mathbb{R}_+}(t), \) $\alpha\in(\frac12,1)$ corresponding to the class of ``rough models".
		
	\subsection{A Numerical illustration: The fake stationary fractional-CIR process.}
		Let the diffusion coefficient be \(\sigma(x) = \nu \sqrt{x}\) for every $x\in\R_+$. For the \(\alpha-\) fractional kernels \( K(t)= K_{\alpha}(t) = \frac{t^{\alpha - 1}}{\Gamma(\alpha)} \mathbf{1}_{\mathbb{R}_+}(t) \)  with \(\alpha \in \left(\frac{1}{2}, 1\right)\), we have that \(a = 0\) in~\eqref{eq:limitResolGamma} and $f_{\alpha, \lambda}:=f_{1, \alpha, 0, \lambda}$ is a probability density called the Mittag-Leffler density function.
		  We further assume for simplification, that \(\mu(t) = \lambda m_0 = \mu_0\) for all  $t \geq 0$.
	In this case, the equation simplifies in the so-called fake stationarity regime as follows:
	\vspace{-.3cm}
	\begin{equation}\label{eq:volt}
		X_t= 	g_0(t)
		+  \frac{\nu}{\lambda}\int_0^t f_{\alpha, \lambda}(t-s)\varsigma(s)\sqrt{X_{s}}dW_s, \;
		\text{with}\;
		g_0(t)
		=
		m_0+(X_0-m_0)
		\Big(1-\int_0^t f_{\alpha,\lambda}(s)\,ds\Big).
		\vspace{-.3cm}
	\end{equation}
    \medskip\noindent
	We introduce an Euler-Maruyama scheme below on the time grid $t_k =t^n_k =\frac{kT}{n}, k=0, \dots, n$, for the semi-integrated form~\eqref{eq:volt}, which write recursively with $R_{\lambda}(t):=1-\int_0^t f_{\alpha,\lambda}(s)\,ds$:
	\vspace{-.3cm}
	{\small 
		\begin{align*} 
			\overline X_{t_{k}} 
			&= \frac{\mu_0}{\lambda} + \big(X_0 -\frac{\mu_0}{\lambda} \big)R_{\lambda}(t_k) + \sum_{\ell=1}^{k} \frac{\nu}{\lambda} \int_{t_{\ell-1}}^{t_{\ell}} f_{\lambda}(t_k-s) \, \varsigma(t_{\ell}) \sqrt{(\overline{X}_{t_{\ell-1}})_+}dW_s= g(t_k) + \frac{\nu}{\lambda} \sum_{\ell=1}^{k} \varsigma(t_{\ell})\sqrt{ (\overline{X}_{t_{\ell-1}})_+} I^{n,l}_k.
			\vspace{-.3cm}
		\end{align*}
	} 
	\noindent where $x_+:=\max\{x,0\}$ and the integrals $\left(I^{n,l}_k= \int_{t_{\ell-1}}^{t_{\ell}} f_{\lambda}(t_k-s) dW_s\right)_{k} $ can be simulated on the discrete  grid \((t^n_k)_{0\leq k\leq n}\) by generating an independent sequence of gaussian vectors \( G^{n,l}, l=1 \cdots n\) using the Cholesky decomposition of a well-defined covariance matrix \(C\).
	The reader is referred to \textit{Appendix A} of \cite{EGnabeyeu2025} for further details on the semi-integrated Euler scheme introduced in this context for the above equation.
	\vspace{-.3cm}
	\begin{figure}[H]
		\centering
		\includegraphics[width=0.93\linewidth]{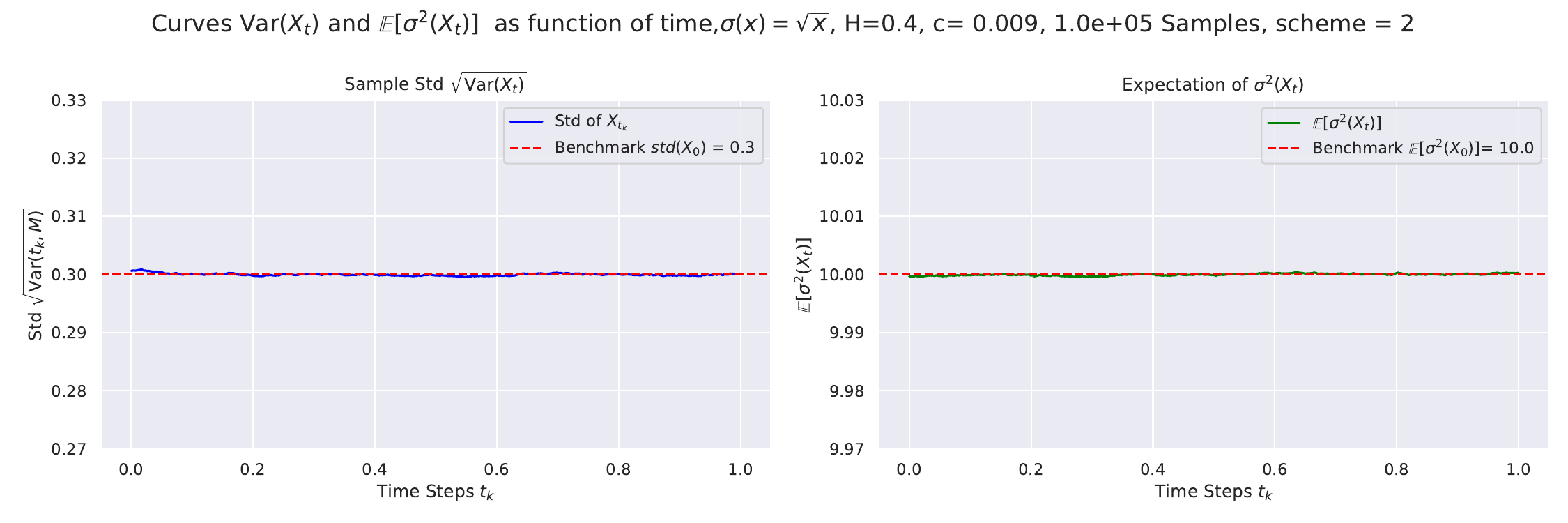}
		\caption{Graph of \( t_k \mapsto \text{StdDev}(t_k, M) \) and \( t_k \mapsto \mathbb{E}[\sigma^2(X_{t_k},M)] \) over the time interval \( [0, T] \), \( T = 1 \), \( H = 0.4 \), \( \mu_0 = 2 \), \( \lambda = 0.2 \), \( v_0 = 0.09 \), and \( \text{StdDev}(X_0) = 0.3 \), \( \nu = 1 \), \( n = 800 \), \( M = 100000 \).}
		\vspace{-.3cm}
	\end{figure}
	\vspace{-.3cm}
	\begin{figure}[H]
		\centering
		\includegraphics[width=0.93\linewidth]{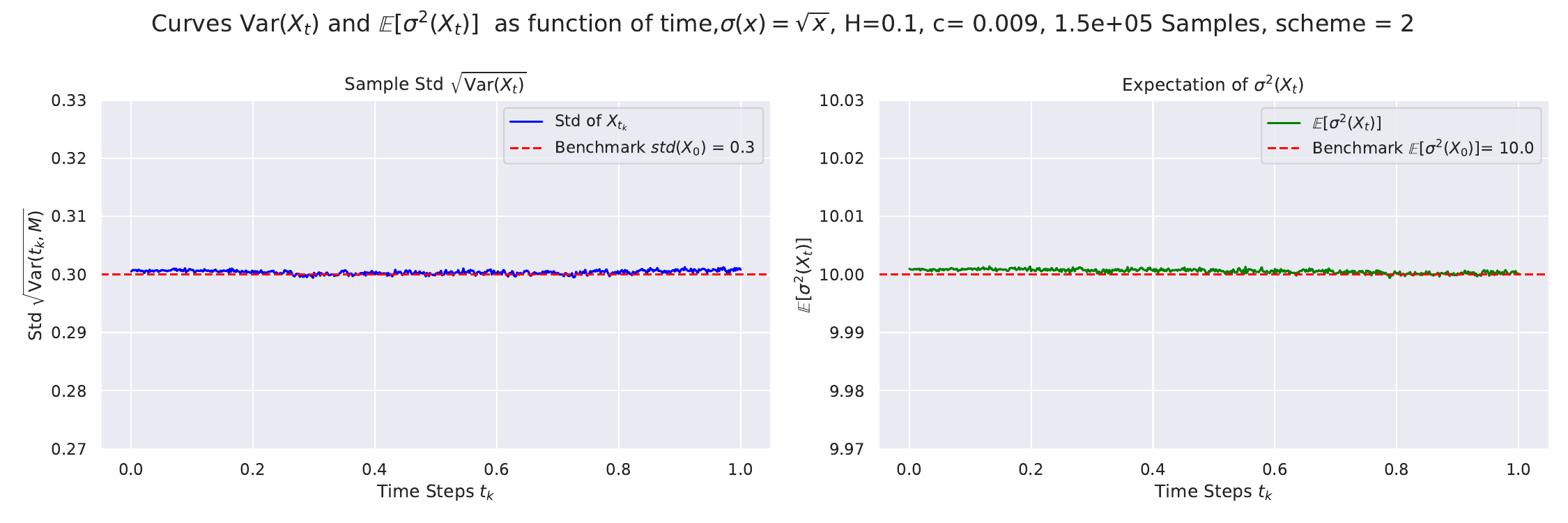
		}
		\caption{Graph of \( t_k \mapsto \text{StdDev}(t_k, M) \) and \( t_k \mapsto \mathbb{E}[\sigma^2(X_{t_k},M)] \) over the time interval \( [0, T] \), \( T = 1 \), \( H = 0.1 \), \( \mu_0 = 2 \), \( \lambda = 0.2 \), \( v_0 = 0.09 \), and \( \text{StdDev}(X_0) = 0.3 \), \( \nu = 1 \), \( n = 800 \), \( M = 100000 \).}
		\vspace{-.3cm}
	\end{figure}
	\medskip
	\noindent {\bf Long-Run \(\mathcal W_2\)-Confluence of the stabilized fractional CIR process~\eqref{eq:volt}:} Let $\lambda,c>0$. Here, we consider \(\ell_\infty=0\) i.e. $\lim_{t\to +\infty} h_{\lambda,c} (t) =0$ and \(\rho=0\) so that  \(a = 0\) in~\eqref{eq:limitResolGamma} (fractional kernel).  Using equations~\ref{eq:FirstMmt} and~\eqref{eq:conditionVar} we have that \( \mathbb{E} X_t \to m_0=\frac{\mu_0}{\lambda} \) and \( \text{Var}(X_t) \to v_0=c\nu^2 m_0 \) as \(t \to +\infty\).
	Moreover, since \(\ell_\infty=0\), Theorem~\ref{Theorem: limiting distribution} implies that the limiting distribution is independent of the initial state. Combined with the uniform \(L^p\)-moment estimates of Proposition~\ref{prop:Momentctrl_contraction_Hol}, for some \(p>2\), this yields convergence in \(\mathcal W_2\). Thus, the stabilized fractional CIR process exhibits long-run \(\mathcal W_2\)-confluence.
	\vspace{-.1cm}
	\begin{figure}[H]
		\centering
		\begin{minipage}{0.49\linewidth}
			\centering
			\includegraphics[width=1.1\linewidth]{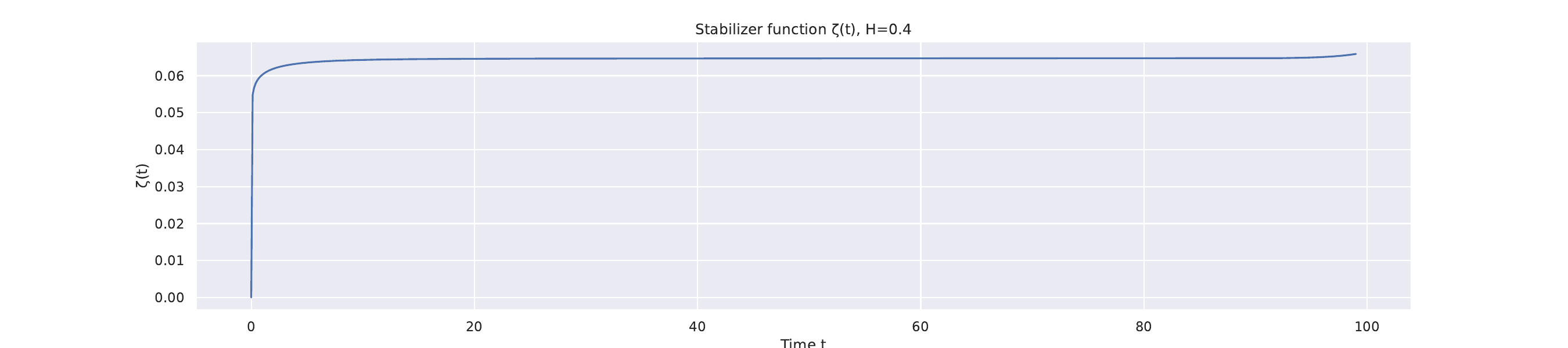}
			\caption{Graph of the stabilizer $ t \to \varsigma_{\alpha,\lambda,c}(t)$}
		\end{minipage}%
		\hfill
		\begin{minipage}{0.49\linewidth}
			\centering
			\includegraphics[width=.7\linewidth]{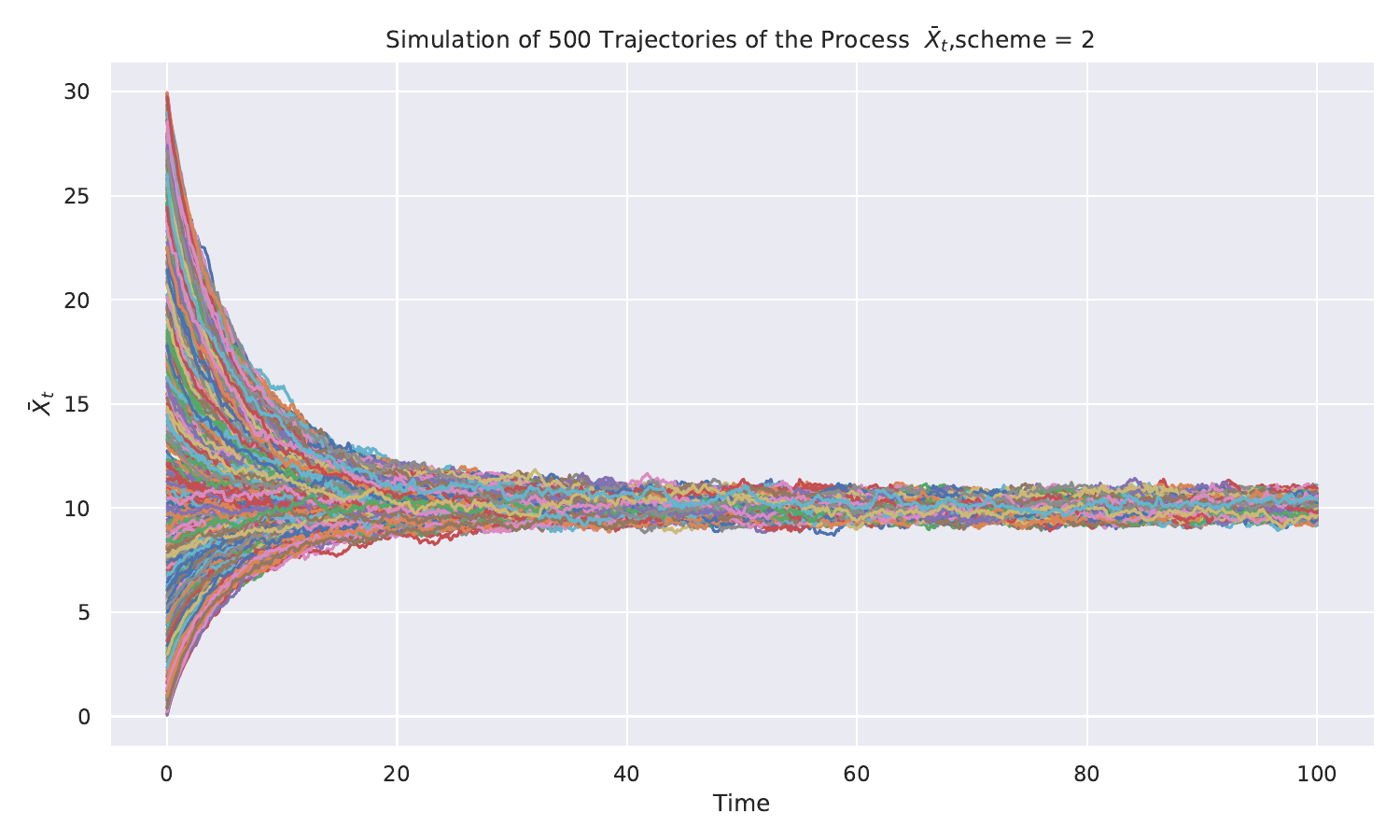}
			\caption{Confluence from $X_0\sim\mathcal{U}(]0,30])$}
		\end{minipage}
		\caption[Stabilization-induced confluence of trajectories.]
		{\textit{Stabilization-induced confluence of trajectories. Left: the stabilizing
				factor $t\mapsto\varsigma_{\alpha,\lambda,c}(t)$. Right: confluence of
				trajectories starting from $X_0$. 
				$H=0.4$, $\lambda=0.2$, $\mu_0=2$, $c=0.05$, and $T=100$.}}
		\label{fig:confluent}
	\end{figure}
	\vspace{-.5cm}
	\noindent Figure~\ref{fig:confluent} shows that the marginals of such a process when
	starting with various initial values are confluent in $\mathcal{W}_2$ as time goes to infinity.
	\subsection{An illustration and SPX fitting of the fake stationary rough Heston model}
	
	In practice, implied-volatility information is represented in several ways. Common approaches include non-parametric volatility grids with spline interpolation, direct modelling of the risk-neutral distribution of the underlying asset, and parametric descriptions of the volatility surface followed by data-driven calibration, as in \cite{GnabeyeuKarkarIdboufous2024}. Stochastic-volatility models provide another approach. The fake stationary rough Heston model belongs to this class; see \citep{EGnabeyeuR2025,Gnabeyeu2026b} and Figures~\ref{fig:Impl_Vol_Calib} and~\ref{fig:Impl_Vol_gen}.
	The model is specified by the parameter vector $\kappa = (\alpha,\lambda,\rho,c,\nu,m_0),
	\;
	m_0 := \frac{\mu_0}{\lambda}.$
	Under Proposition~\ref{prop:mainStab}, the variance of the initial variance variable is $v_0 = c\nu^2m_0.$
		For a maturity \(T\) and strike \(K\), the model-implied volatility
	\(\sigma_{\textsc{iv}}^{\mathrm{Model}}(\kappa,K,T)\) is the Black--Scholes
	volatility that reproduces the European option price generated by the model.
	In the figures below, option prices are computed using Lewis's Fourier
	pricing method (see, e.g.,
	\cite[Section~4.5]{Schmelzle2010}), which applies an inverse Fourier transform to the
	characteristic function of the log-price process obtained in~\eqref{eq:riccati_heston_final}-\eqref{eq:time_zero_frac_heston}.
	With non-zero interest rates or dividends, one includes the discount factor and uses the corresponding forward price
	 denoted \(F_T\) for expiry \(T\). The model's Riccati-Volterra equation, with \(p\) denoting the transform variable applied to the log-price, is
	 \vspace{-.2cm}
	\begin{equation}\label{eq:Char_p}
	\psi_p(t)=
	\int_0^t K_\alpha(t-s)
	\left[
	\frac12(p^2-p)
	+\bigl(\rho\nu\varsigma_{\lambda,c}(T-s)p-\lambda\bigr)\psi_p(s)
	+\frac{\nu^2}{2}\varsigma_{\lambda,c}^2(T-s)\psi_p(s)^2
	\right]ds .
	\vspace{-.2cm}
	\end{equation}
	For the calibration considered here, \(V_0\) follows a \textit{Gamma distribution} with shape and scale given as a function of the mean \(m_0\) and variance \(v_0=c\nu^2m_0\).
	The characteristic function of the forward log return \(\log(S_T/F_{T})\) follows from the affine representation, after preconditioning by $V_0\sim\Gamma\left(
	\frac{m_0}{c\nu^2},\,c\nu^2
	\right),$
	\vspace{-.1cm}
	\begin{align}
	\Phi_T(w)
	&:=\mathbb{E}\!\Big[ e^{iw \log(\frac{S_T}{F_{T}})} \Big]=
	\mathbb{E}\!\Big[
	\mathbb{E}\!\Big[ e^{iw \log(\frac{S_T}{F_{T}})}\mid V_0 \Big]
	\Big]=
	\exp\!\Big(\lambda m_0\int_0^T\psi_{iw}(s)\,ds\Big)
	\mathbb E\!\Big[
	\exp\!\Big(V_0(I^{1-\alpha}\psi_{iw})(T)\Big)
	\Big],\nonumber \\
	&\hspace{3.5cm}= 
	\exp\!\Big(\lambda m_0 \int_0^T\psi_{iw}(s)\,ds\Big)
	\left[
	1-c\nu^2(I^{1-\alpha}\psi_{iw})(T)
	\right]^{-\frac{m_0}{c\nu^2}}
	\end{align}
	on the principal branch of the logarithm, provided the term in square brackets remains in its domain of definition.
	 The resulting calibrated parameter values are reported in Table~\ref{RH:tab:params_nopena}.
	\vspace{-.2cm}
	\begin{table}[H]
		\centering
		\begin{tabular}{l||ccccccc}
			\toprule
			Model
			& $\alpha = H+\frac{1}{2}$
			& $\lambda$
			& $\rho$
			& $c$
			& $\nu$
			& $m_0$
			& $v_0$ \\
			\midrule \midrule
			Fake stationary rough Heston
			& $0.8999$
			& $2.1879$
			& $-0.9999$
			& $0.1074$
			& $0.8448$
			& $0.0450$
			& $0.0034503$ \\
			\bottomrule
		\end{tabular}
		\caption[Calibrated parameters for the fake stationary rough Heston model.]
		{\textit{Calibrated parameters for the fake stationary rough Heston model.}}
		\label{RH:tab:params_nopena}
	\end{table}
	\vspace{-.3cm}
		\begin{figure}[H]
		\centering
		\includegraphics[width=.9\linewidth]{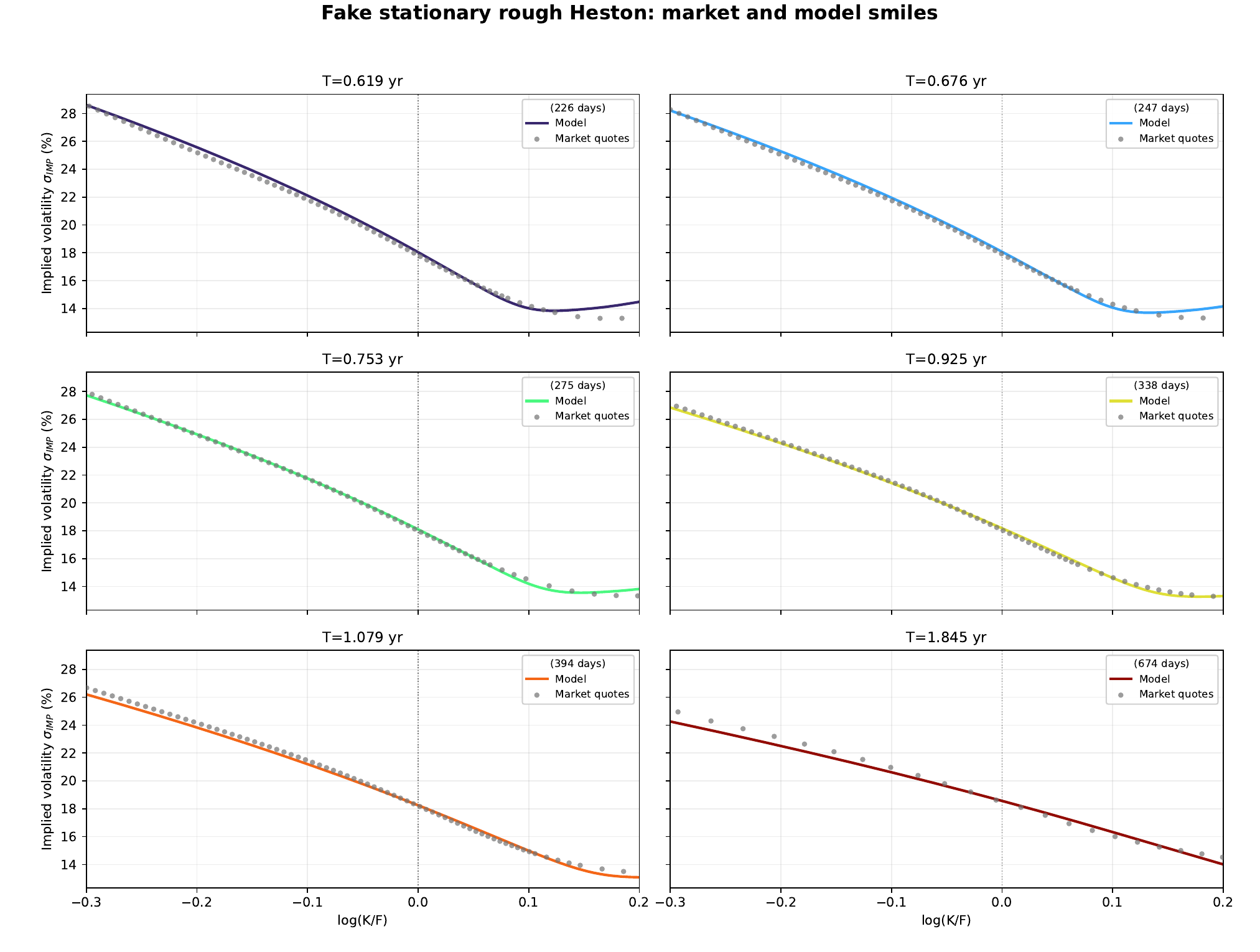}
		\caption{
			\textit{Market mid implied-volatility quotes (gray markers) observed on 15 February 2023 and corresponding calibrated smiles generated by the fake stationary rough Heston model (colored curves).
		}}
		\label{fig:Impl_Vol_Calib}
	\end{figure}
	Note that a \(\Gamma\)-law for $V_0$ is a modelling choice.
	The
	calibration expiries provide the direct comparison with observed option
	prices. It is clear in Figure~\ref{fig:Impl_Vol_Calib} that the fake stationary rough Heston model succeeds at producing or generating the desired smile for different 
	maturities. 
	The
	horizontal axis is log-forward moneyness, \(\log(K/F_T)\).
	Figure~\ref{fig:Impl_Vol_gen} displays the resulting implied-volatility
	smiles and surface for ten maturities, from $38$ days to two years with parameters in~\ref{RH:tab:params_nopena}. The log-forward moneyness, \(\log(K/F_T)\) are between $-0.6$ and $0.4$.
	\vspace{-.4cm}
	\begin{figure}[H]
		\centering
		\begin{minipage}[c]{0.48\linewidth}
			\centering
			\includegraphics[width=.8\linewidth]{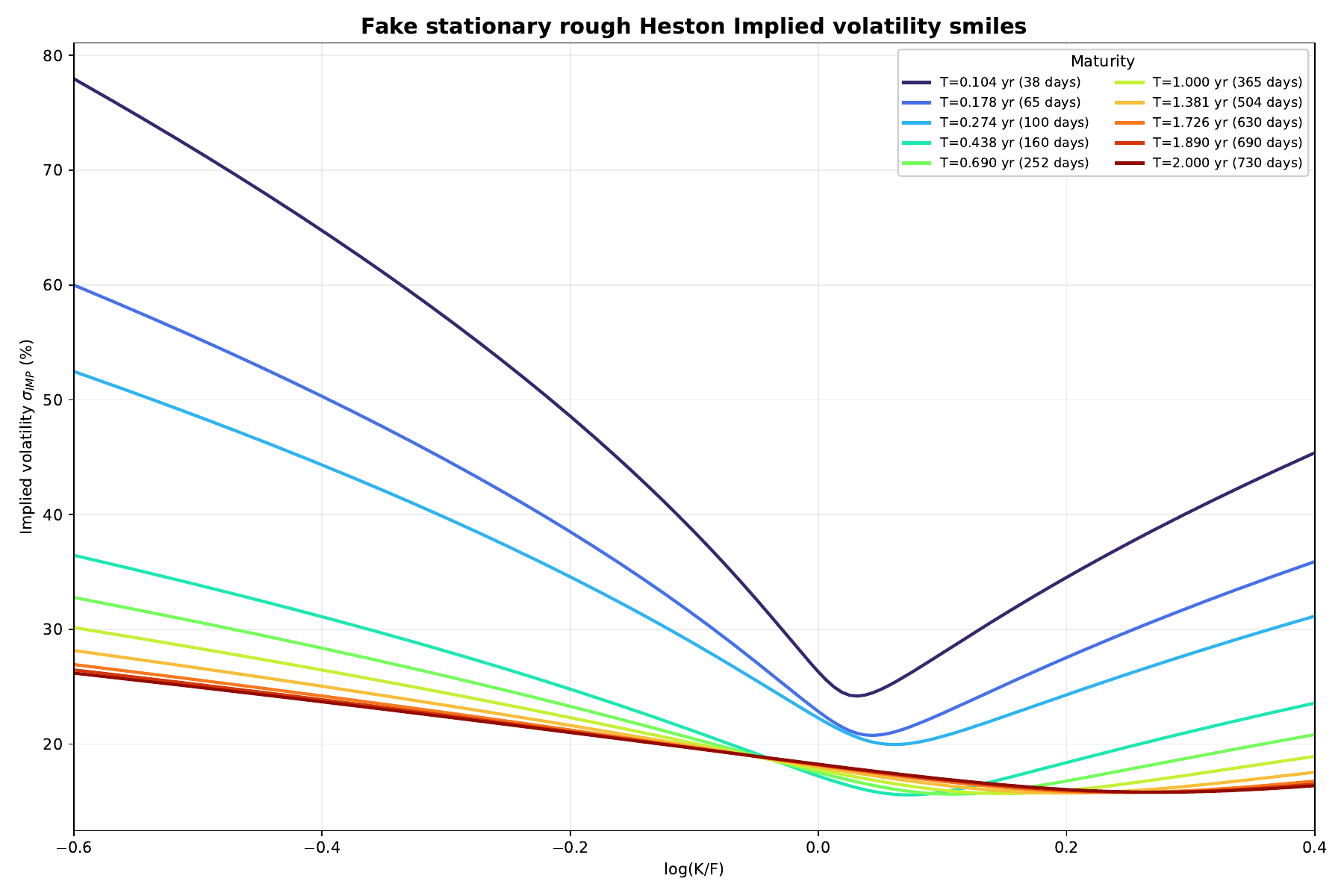}
		\end{minipage}
		\hfill
		\begin{minipage}[c]{0.48\linewidth}
			\centering
			\includegraphics[width=1.2\linewidth]{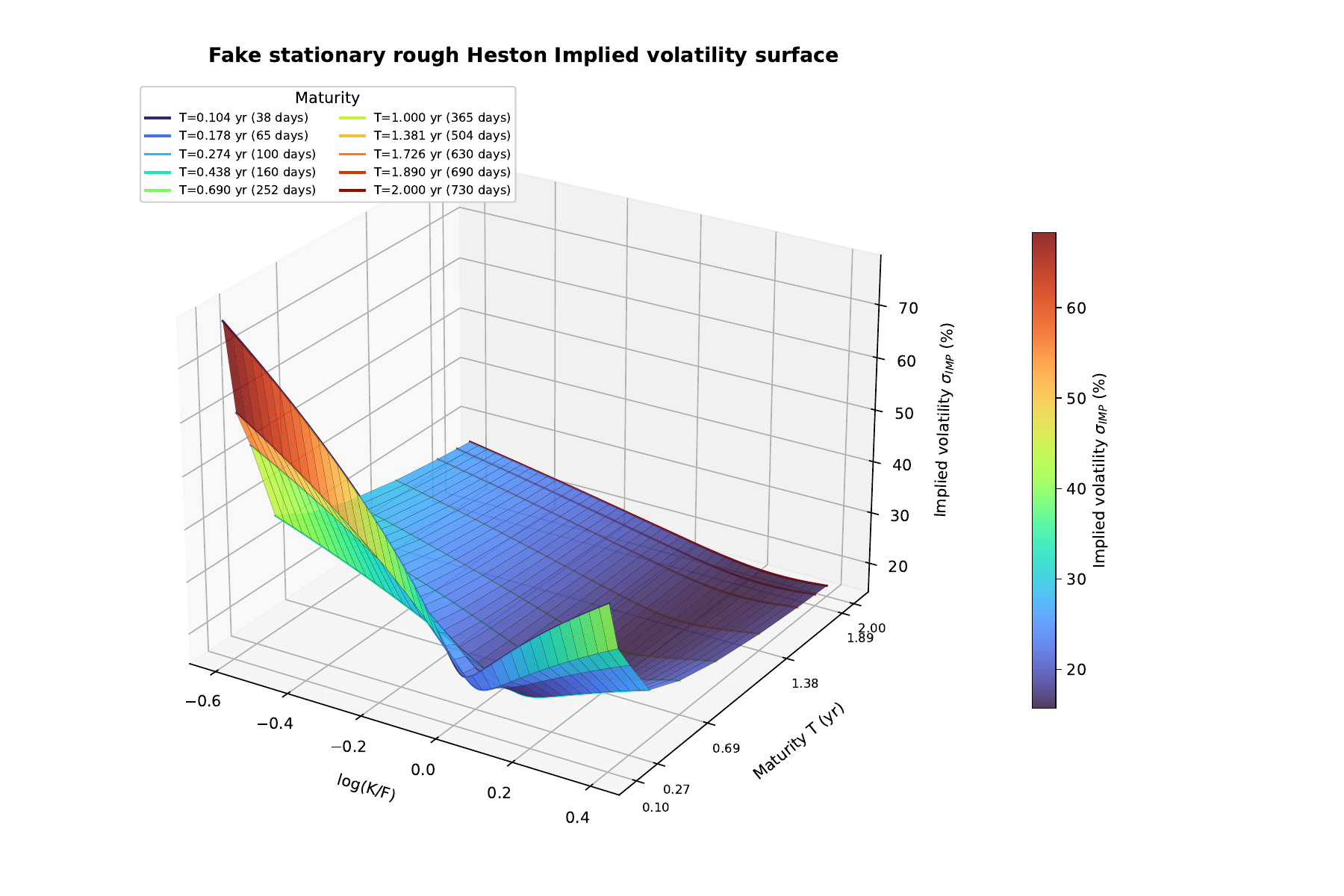}
		\end{minipage}
		\vspace{-.5cm}
		\caption{\textit{
				Left: model-implied volatility smiles as functions of
				log-forward moneyness, \(\log(K/F_T)\).
				Right: the corresponding implied-volatility surface, using the common parameter set across maturities.
				Both panels use the same calibrated parameter set and show maturities
				from 9 days to 730 days.
		}}
		\label{fig:Impl_Vol_gen}
		\vspace{-.4cm}
	\end{figure}
	Figure~\ref{fig:Impl_Vol_gen} illustrates how a single parameter set
	generates maturity-dependent, negatively skewed volatility smiles and
	their associated surface. The same set of model parameters reproduces important features and provide a tractable description of volatility across a large range of maturities (consistency).  The shortest-maturity smiles are particularly
	sensitive to the model parameters; moreover, maturities shorter than the
	first calibration expiry represent short-end extrapolations. The figures
	therefore illustrate the model's term-structure behaviour.

	\vspace{-.3cm}
	\section{Proofs of the main results}\label{sect:proofMresult}
	\vspace{-.2cm}
	\subsection{Proof of Theorems~\ref{theorem:riccati-extension} and~\ref{T:VolSqrt} }\label{subsect:proofMresult1} 
	\noindent {$\rhd$ \textit{ Preliminaries}. As a first preliminary,  we state the following lemma on the existence and uniqueness of a positive solution to a measure-forced linear Volterra equation, which will be instrumental in our subsequent analysis.}
	\begin{Lemma}\label{Lem:MeasLinVolterra}
		Let $K \in L^2_{\mathrm{loc}}(\mathbb{R}_+;\mathbb{R})$  be completely
		monotone, not identically zero, and satisfy Assumption~\ref{assump:kernelVolterra}~$(i)$. Let $v \in \mathbb{R}$, $\xi \in \mathbb{M}_{\mathrm{lf}}(\mathbb{R}_+;\mathbb{R})$ and $G \in
		L^2_{\mathrm{loc}}(\mathbb{R}_+;\mathbb{R})$. Then the linear Volterra equation
		\vspace{-.2cm}
		\begin{equation}\label{eq:MeasLinVolterra}
			\chi = v + \int_{[0,\cdot)} K(\cdot-s)\,\xi(\dd s) + \big(K * ( G\chi) \big) (\cdot)
			\vspace{-.2cm}
		\end{equation}
		has a unique solution $\chi \in L^2_{\mathrm{loc}}(\mathbb{R}_+;\mathbb{R})$. Moreover, if $v\ge0$ and $\xi$ is a nonnegative measure, then $\chi \geq 0$.
	\end{Lemma}
	\noindent\textbf{\textit{ Proof of Lemma~\ref{Lem:MeasLinVolterra}.}}
	Let $\rho(t)=e^{-t}\mathbf{1}_{\{t\geq 0\}}$ and for every $n\geq1$ consider a nonnegative mollifier $\rho_n(t)=n\rho(nt)= ne^{-nt}\mathbf{1}_{\{t\geq 0\}}$, shrinking in the sense of distributions to $\delta_0$ (adjusted to stay supported on $\mathbb{R}_+$). Define $f_n(t)=\int_{[0,t]}\rho_n(t-s)\,\xi(ds)$ for $t\geq0$ and let $\xi_n(\dd s)=f_n(s)\dd s$. We consider the linear Volterra equation
	\vspace{-.2cm}
	\begin{equation}\label{eq:MeasLinVolterra2}
		\chi^n = v + \int_{[0,\cdot)} K(\cdot-s)\,\xi_n(\dd s) + \big(K * ( G\chi^n) \big) (\cdot)
		\vspace{-.1cm}
	\end{equation}
	One checks that $\|f_n\|_{L^1([0,T])}\leq|\xi|([0,T])$ for every $0<T<\infty$, so that $f_n\in L^1_{\mathrm{loc}}(\mathbb{R}_+;\mathbb{R})$.
	We deduce from \cite[Corollary B.3]{abi2019affine}  (see also~\cite[Theorem 9.3.6]{gripenberg1990}) that the linear deterministic Volterra equation~\eqref{eq:MeasLinVolterra2} has a unique global solution $\chi^n \in L^2_{\mathrm{loc}}(\mathbb{R}_+;\mathbb{R})$.
	Fix $0<T<+\infty$ and note that since \( K \ast \xi\in L^2([0,T]; \mathbb{R})\), we have that \( \int_{[0,\cdot)} K(\cdot-s)\,\xi_n(\dd s) =  K \ast  \rho_n \ast\xi = \rho_n\ast(K \ast \xi)\to \int_{[0,\cdot)} K(\cdot-s)\,\xi(\dd s) \) in \( L^2([0,T]; \mathbb{R})\) by classical results from Fourier analysis (see, e.g.~\cite[ 
	Theorem 2.16]{LiebLoss2001}) .
	Define for short $g_n(t):=\int_{[0,t)}K(t-s)\,\xi_n(ds)$. By the elementary inequality $(a+b)^p \leq 2^{p-1}(a^p+b^p),\; a,b>0, p\geq1$, together with Cauchy--Schwarz, we obtain, for every $n,m\in\mathbb{N}$,
	\vspace{-.2cm}
	\begin{equation}\label{eq:CauchyDrift}
		|\chi^n-\chi^m|^2(t)
		\leq
		2|g_n-g_m|^2(t)
		+
		2\left(
		\int_0^t |K(t-s)|^2|G(s)|^2\,ds
		\right)
		\left(\int_0^t |\chi^n-\chi^m|^2(s)\,ds\right).
		\vspace{-.2cm}
	\end{equation}
	Integrating both side of~\eqref{eq:CauchyDrift} and applying Gr\"onwall's lemma yields
	\begin{align*}
		\int_{0}^{t}|\chi^n-\chi^m|^2(s)\,ds
		&\leq 2 \int_{0}^{t} |g_n-g_m|^2(s)\exp\Big(
		2 \int_s^t \int_0^u |K(u-v)|^2|G(v)|^2\,dv\,du
		\Big)\,ds\\
		& \leq 2 \exp\Big(
		2 \int_0^t |G(s)|^2  \int_0^{t-s} |K(r)|^2\,dr\,ds
		\Big)\int_{0}^{t} |g_n-g_m|^2(s)\,ds.
	\end{align*}
	where the last inequality follows from Fubini--Tonelli's theorem. Therefore, we deduce that 
	\[
	\|\chi^n-\chi^m\|_{L^2([0,T])}
	\leq
	\sqrt2\,e^{\|K\|_{L^2([0,T])}^2\|G\|_{L^2([0,T])}^2}
	\|g_n-g_m\|_{L^2([0,T])}.
	\]
	Since $(g_n)_{n\geq1}$ is Cauchy in $L^2([0,T])$, it follows immediately that $(\chi^n)_{n\geq1}$ is Cauchy in $L^2([0,T])$, then $\chi^n\to\chi$ in $L^2([0,T])$. 
	We now show that $\chi$ solves the equation~\eqref{eq:MeasLinVolterra} with measure forcing term on $[0,T]$.
	Since \(\chi^n\to\chi \) and \( \int_{[0,\cdot)} K(\cdot-s)\,\xi_n(\dd s) \to \int_{[0,\cdot)} K(\cdot-s)\,\xi(\dd s) \) in \( L^2([0,T]; \mathbb{R})\),
	it suffices to show that \(	K \ast ( G\chi^n) \to K \ast ( G\chi)\) in \( L^2([0,T]; \mathbb{R}).\)
	To this end, we apply Young's inequality
	to obtain
	\begin{align*}
		& \ \| K \ast ( G\chi^n) - K \ast ( G\chi) \|_{L^2([0,T])}
		\leq \| K \|_{L^2([0,T])} \int_0^T |G(s)| |\chi^n(s) - \chi(s)| \, ds \\
		& \hspace{5.5cm}\leq \| K \|_{L^2([0,T])} \|G \|_{L^2([0,T])}
		\| \chi^n - \chi \|_{L^2([0,T])}.
	\end{align*}
	where the last inequality comes form Cauchy-Schwarz inequality.
	Since the right-hand side tends to zero, it follows that passing to the limit in~\eqref{eq:MeasLinVolterra2} yields that $\chi$ solves the true, measure-forced equation~\eqref{eq:MeasLinVolterra}  on $[0,T]$.
	Finally, note that $K \ast \xi \in L^2_{\mathrm{loc}}(\mathbb{R}_+; \mathbb{R})$. Then by \cite[Theorem~B.1]{abi2019affine}, the equation \eqref{eq:MeasLinVolterra} admits a unique maximal solution. Since $\chi$ is a global solution, it must coincide with the unique maximal solution on $\mathbb{R}_+$.
	\noindent Now, if $v\ge0$ and $\xi$ is a nonnegative measure, then
	as $\rho_n\geq 0$, we have $f_n\geq 0$ so that $\xi_n(\cdot):=f_n(\cdot)\,\dd\cdot$ is a nonnegative-density measure. By \cite[Theorem C.2]{abi2019affine}, the unique solution to~\eqref{eq:MeasLinVolterra2} is nonnegative and
	$\chi^n \ge 0 \ \forall n$ passes to $\chi \ge 0$ a.e. under $L^2$ convergence (along a subsequence pointwise, hence a.e.). \qed
	
	\bigskip
	\textbf{\textit{Proof of Theorem~\ref{theorem:riccati-extension}.}} 
	For part~$(a)$, it suffices to prove the existence  of a noncontinuable solution\footnote{More generally, in \(  L^2([0,T]; \mathbb{C}) \), a non-continuable solution of~\eqref{RicVolSqrt2} is a pair \((\psi, T_{\max})\) with \(T_{\max} \in (0,\infty]\) and $\psi \in L^2_{\mathrm{loc}}([0, T_{\max})); \mathbb{C})$, such that \(\psi\) satisfies~\eqref{RicVolSqrt2} on 
		\([0, T_{\max})\) and \(\|\psi\|_{L^2(0, T_{\max})} = \infty\) whenever \(T_{\max} < T\). 
		If \(T_{\max} \geq T\), we call \(\psi\) a global solution of~\eqref{RicVolSqrt2} on $[0,T]$.
	}
	with a non-positive real part and the uniqueness of global solutions. 
	We proceed in three steps, starting with the existence of local solutions. 
	In contrast to~\cite{friesen2022volterra}, where the measure $m$ is approximated by a sequence of functions in $L^1_{\mathrm{loc}}(\mathbb{R}_+,\mathbb{C}_-)$, we work directly with the measure forcing term. The existence argument in~\cite{friesen2022volterra} relies on the relative compactness in $L^p(0,T)$ of bounded subsets of the Sobolev--Slobodeckij space $W^{\eta,p}(0,T)$ for any $T>0$; see, e.g., \cite[Theorem~2.1]{flandoli1995martingale}. Recall that relative compactness in $L^p$ spaces is classically characterized by the Kolmogorov--Riesz--Fréchet theorem; see, for instance, \cite[Theorem~4.26]{brezis2010functional}. Instead, we employ a Banach fixed-point argument  (see also, e.g., \cite{HorstXuZhang2024b}) to establish the existence of a non-continuable solution to \eqref{RicVolSqrt2} directly in the presence of the measure forcing term $m$. This is  obtained in the first two steps below and to this end,
	we shall first establish the existence and uniqueness of local solutions for which we introduce, for any $T>0$ and $M>0$ the set
	\begin{equation}\label{eqn.100}
		\mathcal{B}_{T,M}:=\big\{ \psi\in L^2\big([0,T],\mathbb{C} \big):\ 	\| \psi \|_{ L^2([0,T])}\leq M \big\}. 
	\end{equation}
	a closed ball in the Hilbert space $ L^2\big([0,T],\mathbb{C} \big)$, hence complete since a closed subset of a complete space is automatically complete in the induced metric.
	We also define for short
	\vspace{-.2cm}
	\begin{align*}\label{eq:Operator}
		g_m(t):=\frac1\lambda\int_{[0,t)} f_{\lambda}(t-s) \, \,m(\mathrm{d}s),\quad \Theta(t):=\frac{1}{\lambda}\frac{\kappa_1}{2}\varsigma^2(T-t) \quad\text{for every} \quad t\in [0,T].
		\vspace{-.2cm}
	\end{align*}
	{\bf Step 1: Existence and uniqueness of local solutions}.
	We will prove that for some $M>0$ and $\delta>0$ that will be specified later, there exists a local solution $\big(\delta,\psi(\cdot, m)\big)$ to \eqref{RicVolSqrt2} with  $\psi(\cdot, m) \in\mathcal B_{\delta,M}$. 
	For this, we introduce the mapping $\mathcal{T}_\delta$ that acts on locally integrable functions $\psi\in\mathcal{B}_{\delta,M}$ according to
	\vspace{-.2cm}
	\begin{equation}\label{eq:Operator}
		\mathcal{T}_\delta\big(\psi\big)(t):=  G_0 (t) + \int_0^t f_{\lambda}(t-s) \Theta(s) \psi^2(s) \, \mathrm{d}s,\quad \text{with}\quad G_0(t) := g_m(t),\quad t\in [0,\delta].
		\vspace{-.2cm}
	\end{equation}
	Since the set $\mathcal{B}_{\delta,M}$ is closed, bounded and convex subsets of the spaces $L^2([0,\delta];\mathbb{C})$, the existence of a local solution follows from Banach's fixed point theorem if we can prove that $\mathcal{T}_\delta$ is a contraction w.r.t.~the $L^2$-norm. 
	One checks that $\mathcal{T}_\delta$ is a self-map and a Contraction for a well chosen $M$. In fact, by Minkowski's inequality, $g_m\in L^2_{loc}(\R,\mathbb{C})$
	and for any $\psi\in\mathcal{B}_{\delta,M}$, Young's inequality gives
	\[
	\|\mathcal{T}_\delta(\psi)\|_{L^2(0,\delta)}
	\leq
	\|g_m\|_{L^2(0,\delta)}
	+
	\overline{\Theta}_\delta
	\|f_\lambda\|_{L^2(0,\delta)} M^2,
	\qquad
	\overline{\Theta}_\delta
	:=
	\frac{\kappa_1}{2\lambda}
	\|\varsigma^2\|_{\infty,T}.
	\]
	Moreover, for $\psi_1,\psi_2\in\mathcal{B}_{\delta,M}$, 
	Cauchy--Schwarz followed by Young's inequality gives
	\[
	\|\mathcal{T}_\delta(\psi_1)-\mathcal{T}_\delta(\psi_2)\|_{L^2(0,\delta)}
	\leq
	2\overline{\Theta}_\delta M
	\|f_\lambda\|_{L^2(0,\delta)}
	\|\psi_1-\psi_2\|_{L^2(0,\delta)}.
	\]
	Now, $\|f_\lambda\|_{L^2(0,\delta)}\to0$ and $\|g_m\|_{L^2(0,\delta)}\to0$ as $\delta\downarrow0$. We can therefore take $M_\delta=M:=2\|g_m\|_{L^2(0,\delta)}$
	and choose $\delta$ sufficiently small such that $4 \overline{\Theta}_\delta \|g_m\|_{L^2(0,\delta)} \|f_\lambda\|_{L^2(0,\delta)} <1$.
	This shows that the operator $\mathcal{T}_\delta$ maps the set $B_{\delta,M_\delta}$ into itself (i.e. $\mathcal{T}_\delta(B_{\delta,M_\delta})
	\subset B_{\delta,M_\delta}$) and that it is a contraction in the $L^2$-sense.
	In particular, for any two fixed points $\psi_1, \psi_2 \in B_{\delta,M_\delta}$, it holds that $\|\mathcal{T}_\delta(\psi_1)-\mathcal{T}_\delta(\psi_2)\|_{L^2(0,\delta)}=0$.
	As a result, $\mathcal{T}_\delta$ admits an a.e. uniquely defined fixed point $\psi(\cdot,m) \in B_{\delta,M_\delta}$. In particular, there exists a local
	solution $(\delta, \psi(\cdot,m))$ of the equation~\eqref{RicVolSqrt2} with the solution being unique up to a set of Lebesgue-measure zero. 
	
	\smallskip
	\noindent {\bf Step 2: Continuation dichotomy, existence and uniqueness of non-continuable solutions.} We are now going to show that any local solution can be extended to a local solution on a larger interval. 
	To this end, let $(t_0,\psi_0(\cdot,m))$ be a local solution of the equation~\eqref{RicVolSqrt2} with $\psi_0(\cdot,m) \in B_{t_0,M_{t_0}}$.  To extend this solution to a larger interval we set for $u\geq0$, $\widetilde{\psi}(u):=\psi(t_0+u)$ and we have that $\widetilde{\psi}$ satisfies the equation
	\vspace{-.2cm}
	\begin{equation}
		\widetilde{\psi}(u) = G_{t_0}(u) + \int_0^u f_{\lambda}(u-v) \Theta(t_0 + v) \widetilde{\psi}^2(v) \, \mathrm{d}v,
		\vspace{-.2cm}
	\end{equation}
	where the function $G_{t_0}: \mathbb R_+ \to \mathbb C$ is defined by
	\vspace{-.2cm}
	\begin{equation}
	G_{t_0}(u)
	=
	g_m(t_0+u)
	+
	\int_0^{t_0}
	f_\lambda(t_0+u-v)\Theta(v)\psi_0(v)^2\,dv.
	\vspace{-.2cm}
	\end{equation}
	One checks that $G_{t_0}\in L^2_{loc}(\R_+,\mathbb{C})$.
	Now, fix two constants $\delta_1\in(0,1)$ and $M_1>0$ to be specified later and introduce the operator $\mathcal{T}_{\delta_1}$
	that acts on the set $\mathcal{B}_{\delta_1,M_1}$ as in~\eqref{eq:Operator} with forcing term $G_{t_0}$.
	Following exactly the same fixed-point arguments as in Step 1 on $[0,\delta]$ and choosing first $\delta_1$ small enough such that $4 \overline{\Theta}_{\delta_1}\|G_{t_0}\|_{L^2(0,\delta_1)} \|f_\lambda\|_{L^2(0,\delta_1)} <1$ and then
	$M_{\delta_1}=M_1:=2\|G_{t_0}\|_{L^2(0,\delta_1)}$, we see that  $\mathcal{T}_{\delta_1}(B_{\delta_1,M_{\delta_1}})
	\subset B_{\delta_1,M_{\delta_1}}$ and that $\mathcal{T}_{\delta_1}$ is a contraction. We conclude that the operator $\mathcal{T}_{\delta_1}$ admits an a.e.~uniquely defined fixed point $\psi_1(\cdot,m)$ in the set $\mathcal{B}_{\delta_1,M_{\delta_1}}$, viewed as a closed, convex, bounded subset of $L^2([0,T];\mathbb{C})$.   Let us now set 
	\vspace{-.1cm}
	\begin{equation}
		\psi(t,m):=\psi_0(t,m)\cdot \mathbf{1}_{(0,t_0]}(t) + \psi_1(t -t_0,m)\cdot \mathbf{1}_{(t_0, t_0+\delta_1]}(t), \quad t\in(0,t_0+\delta_1] . 
	\end{equation}
	The function $\psi(\cdot,m)$ is a local solution to the equation~\eqref{RicVolSqrt2} on the interval $(0,t_0+\delta_1]$. Indeed, $\psi(t,m) = \psi_0(t,m)$ on $(0,t_0]$ while for $t\in (t_0,t_0+\delta_1]$,  
	\vspace{-.2cm}
	\begin{align*}
		\psi(t) & ~= \psi_1 (t-t_0) = \mathcal T_{\delta_1}\circ \psi_1(t-t_0) = G_{t_0}(t-t_0) + \int_0^{t-t_0} f_{\lambda}(t-t_0-s) \Theta(s) \psi_1^2(s) \, \mathrm{d}s,\\
		& ~=  g_m (t) + \int_0^t f_{\lambda}(t-s) \Theta(s) \psi^2(s) \, \mathrm{d}s.
		\vspace{-.2cm}
	\end{align*}
	In view of the above, we deduce that the largest interval on which local solutions can be defined is an open interval of the form  $(0,T_{\max})$ and there exists a function $\psi(\cdot,m)$ on $(0,T_{\max})$ such that for all $t_0<T_{\max}$ the function $\psi(\cdot,m)$ is local solution on $(0,t_0]$. 
	This yields that Equation~\eqref{RicVolSqrt2} admits a noncontinuable solution \((\psi(\cdot,m), T_{\max})\)
		having the standard continuation alternative on the maximal interval $(0,T_{\max})$: either $T_{\max}\geq T$ or $T_{\max}<T$ and $\limsup_{H\uparrow T_{\max}}\|\psi\|_{L^2([0,H])}=\infty.$ 
	
	\smallskip
	\noindent {\bf Step 3: Non-positivity statement for the real part of the solution of the Riccati--Volterra equation.}
	We now deal with the non-positivity of the real part of solutions to the deterministic Volterra equation~\eqref{RicVolSqrt2}.
	Let $\psi_r(\cdot, m)$ and $\psi_i(\cdot, m)$ denote the real and imaginary parts of $\psi(\cdot, m)$. They satisfy the equations
	\vspace{-.4cm}
	\begin{align}
		\qquad\; \psi_r(t)&= \frac1\lambda\int_{[0,t)} f_{\lambda}(t-s) \, \,\operatorname{Re}(m)(\mathrm{d}s)  + \int_0^t f_{\lambda}(t-s) \Theta(s) (\psi_r^2(s)-\psi_i^2(s)) \, \mathrm{d}s
		\label{eq:RealPart} \\
		\psi_i(t) &=  \frac1\lambda\int_{[0,t)} f_{\lambda}(t-s) \, \,\operatorname{Im}(m)(\mathrm{d}s) + 2\int_0^t f_{\lambda}(t-s) \Theta(s) \psi_r(s) \psi_i(s) \, \mathrm{d}s.
		\label{eq:ImPart}
	\end{align} 
	We first observe that,  
	on the interval $[0,T_{\rm max})$, the function   $\chi :=- \psi_r(\cdot, m) $ satisfies the linear equation
	\begin{equation}\label{eq:mod_chi}
		\chi(t) =   \int_{[0,t)} f_{\lambda} (t-s) \big( -\frac1\lambda\operatorname{Re}(m)(\mathrm{d}s) +  \Theta(s) \psi_i^2(s) \mathrm{d}s \big)  + \int_0^t f_{\lambda}(t-s) \Theta(s)\psi_r (s)  \chi(s) \mathrm{d}s .
	\end{equation}
	which has by Lemma~\ref{Lem:MeasLinVolterra} a unique solution  
	\(\chi\in L^2_{\text{loc}}(\mathbb{R}_+; \mathbb{R})\) with $\chi\ge0$, owing to assumption \ref{ass:int_holregul} on \(K\) and the fact that $\mathbb{M}_{\mathrm{lf}}(\mathbb{R}_+;\mathbb{R})\ni -\operatorname{Re}(m)(\cdot) +  \Theta(\cdot) \psi_i^2(\cdot) \geq0$, i.e. there exists an $\R_+$-valued continuous solution  $\chi$ to~\eqref{eq:mod_chi}.
	\noindent Then, the solution $\psi$ to the Riccati--Volterra equation \eqref{RicVolSqrt2} is such that $\operatorname{Re}(\psi(t,m)) \leq 0$ on $[0,T_{\rm max})$. 
	
	\smallskip
	\noindent {\bf Step~4: existence, uniqueness of global solutions and $L^p$-bounds.} 
	Our aim is to prove that $T_{\rm max} \geq T$ by showing that
	\vspace{-.2cm}
	\begin{equation}\label{eq:temptmax}
		\sup_{H< T_{\rm max}} \|\psi(\cdot,m)\|_{L^2([0,H])} < \infty.  
		\vspace{-.2cm}
	\end{equation}
	Let us define $h_m(t):=\frac1\lambda\int_{[0,t)} f_{\lambda}(t-s)\left|\operatorname{Im}(m)\right|(\mathrm{d}s)$ and $\ell_m(t):=\frac1\lambda\int_{[0,t)} f_{\lambda}(t-s)\operatorname{Re}(m)(\mathrm{d}s)  - \int_0^t f_{\lambda}(t-s) \Theta(s)h_m^2(s) \, \mathrm{d}s$ on $\mathbb{R}_+$ and let $g\in L^2_{\mathrm{loc}}((0,T_{\max});\mathbb{R}_+)$ be the unique solution (thanks to Lemma~\ref{Lem:MeasLinVolterra}) of the linear equation 
	\vspace{-.3cm}
	\begin{equation}
		g(t)= h_m(t) + 2\int_0^t f_{\lambda}(t-s) \Theta(s) \psi_r(s) g(s) \, \mathrm{d}s.
		\label{eq:gLin}
		\vspace{-.2cm}
	\end{equation} 
	Owing to~\eqref{eq:ImPart} and~\eqref{eq:gLin} and applying Lemma~\ref{Lem:MeasLinVolterra} to the equations satisfied by $g \pm \psi_i$ yields that $|\psi_i(\cdot,m)| \leq g$ on $(0,T_{\max})$.
	Noticing that $h_m$ is a supersolution of Equation~\ref{eq:gLin}
	and using that $\operatorname{Re}(m)\leq0$ in~\eqref{eq:RealPart} together with Step~3, we get
	that the real and imaginary parts of $\psi$ satisfy the inequalities
	\vspace{-.2cm}
	\begin{equation}\label{eq:IneqComp}
		\left|\operatorname{Im}(\psi(t,m))\right| \leq h_m(t)
		\qquad\text{and}\qquad \ell_m(t) \leq \operatorname{Re}(\psi(t,m)) \leq 0,
	\end{equation}
	Hence, for any horizon \(0<H\leq T<\infty\), as $\ell_m$ and $h_m$ are real-valued functions, we obtain owing to Young's inequality and for every $p\in[1,2\beta]$,
	\vspace{-.2cm}
	\begin{align}\label{eq:L^pbound}
		\left\|\psi(\cdot,m)\right\|_{L^p([0,H])}
		&\leq
		\|\ell_m\|_{L^p([0,H])}
		+
		\|h_m\|_{L^p([0,H])} \nonumber\\
		&\leq \frac2\lambda \|f_\lambda\|_{L^p([0,H])} |m| ([0,H]) +  \frac1{\lambda^2} \overline{\Theta}_T \|f_\lambda\|_{L^p([0,H])}  \|f_\lambda\|^2_{L^2([0,H])}\big(|m| ([0,H]) \big)^2 
	\end{align}
	which is finite, since by assumption, $f_\lambda$ has finite $L^p$-norm on any bounded interval for every $p\in[1,2\beta]$.
	This yields the desired bounds in part~(b). 
	Now, in view of the estimate established in~\eqref{eq:L^pbound}, any noncontinuable solution is a global solution which is locally square integrable.
	Indeed, if $T_{\rm max} < T < \infty$, then by~\eqref{eq:L^pbound}, $\|\psi(\cdot,m)\|_{L^2_{(0,T_{\rm max})}}<\infty$ which implies \eqref{eq:temptmax}.
	In this case $\psi(\cdot,m)$ is well-defined and locally square integrable in $(0,T_{\rm max}]$ and so the function $\psi(\cdot,m)$ can be extended to a larger interval, which contradicts the definition of $T_{\rm max}$. Hence, $T_{\max}\geq T$ and the solution is global i.e. $\psi\in L^2([0,T]; \mathbb{C})$.
	The uniqueness of global solution follows from the uniqueness of local solutions since any global solution is a local solution when restricted to finite time intervals. 
	This completes the proof of part~(a). 
	
	\smallskip
	\noindent {\bf Step~5: Continuity.} 
	The continuity here follows from the fact that
	for any two functions $f\in L^p(\mathbb{R}_+;\mathbb{R})$ and
	$h\in L^q(\mathbb{R}_+;\mathbb{R})$, with
	$p,q\in[1,\infty]$ satisfying $\frac{1}{p}+\frac{1}{q}=1,$
	the convolution $f*h$ is continuous on $\mathbb{R}_+$. In particular, to prove part (c), in view of ~\eqref{eq:measureFLplce}, it suffices to show that \( K * F(T-\cdot, \psi) \) is continuous on \( \mathbb{R}_+ \).
	By assumption, $K\in L^3_{\mathrm{loc}}(\mathbb R_+)$.
	Moreover, since $\beta\geq\frac32$, part (b)
	implies that
	$\psi\in L^3_{\mathrm{loc}}(\mathbb R_+)$.
	Using the quadratic growth bound \(|F(\cdot, \psi)| \leq (\lambda\vee \frac{\kappa_1}{2} \|\varsigma^2\|_\infty) (1 + |\psi|^2)\),
	we obtain
	$F(T-\cdot,\psi(\cdot))
	\in L^{3/2}_{\mathrm{loc}}(\mathbb R_+)$.
	Since $K\in L^3_{\mathrm{loc}}(\mathbb R_+)$
	and $F(T-\cdot,\psi(\cdot))
	\in L^{3/2}_{\mathrm{loc}}(\mathbb R_+)$,
	the continuity property of convolutions yields the result.
	This proves the assertion  
	in part (c). \qed
	
	\medskip
	\noindent {\bf Proof of Corollary~\ref{corol:riccati-extension}:} Applying Theorem~\ref{theorem:riccati-extension} on each finite interval
	with the constant coefficient $\varsigma_\infty^2$ yield existence and uniqueness of a solution
	$\psi_\infty\in L^2_{\mathrm{loc}}(\R_+;\mathbb{C})$
	to the limiting autonomous equation~\eqref{eq:measureFLplce2}. The estimate~\ref{eq:LpBounds} gives global
	$L^1\cap L^2$-integrability. 
	For the stability assertion, note that $\mbox{\bf 1}_{\{0\le t \le T\}}\varsigma^2(T-t)\to \lim_{T \to +\infty} \varsigma^2(T):= \varsigma^2_\infty $ for every $t\!\in \R_+$ as $T\to +\infty$.
	The comparison estimates~\eqref{eq:IneqComp} established in
	Step~4 of the proof of
	Theorem~\ref{theorem:riccati-extension}
	provide an $L^1\cap L^2$ majorant for
	$\psi_T$, independent of $T$.
	\vspace{-.1cm}
	\begin{equation}\label{eq:UniformM}
		|\psi_T(t)| \leq G_m(t)
		:= 2H_m(t)
		+ \frac{\kappa_1\|\varsigma\|_\infty^2}{2\lambda}
		(f_\lambda * H_m^2)(t), \quad\text{with}\quad H_m(t) := \frac{1}{\lambda}(f_\lambda * |m|)(t).
		\vspace{-.1cm}
	\end{equation}
	Since $f_\lambda \in L^1 \cap L^2$ and
	$|m|(\mathbb{R}_+) < \infty$, Young's inequality gives $H_m \in L^1 \cap L^2,
	\;
	H_m^2 \in L^1,
	\;
	G_m \in L^1 \cap L^2.$
	To prove the stability assertion, fix $R>0$ and let $T\geq R$.
	Set $\delta_T:=\psi_T-\psi_\infty$, $q_T(s):=\frac{\kappa_1}{2\lambda}\varsigma^2(T-s)$, $q_\infty:=\frac{\kappa_1}{2\lambda}\varsigma_\infty^2$ and $q_*:=\frac{\kappa_1}{2\lambda}\|\varsigma^2\|_\infty$.
	Subtracting the two Riccati--Volterra equations
	and applying Young's and Cauchy--Schwarz inequalities
	with exponential weights yields for every $\eta>0$
	\[
	\begin{aligned}
		\|e^{-\eta\cdot}\delta_T\|_{L^2(0,R)}
		\leq \|e^{-\eta\cdot}f_\lambda\|_{L^2(\R_+)}
		\Big(
		2q_*\|G_m\|_{L^2(0,R)}
		\|e^{-\eta\cdot}\delta_T\|_{L^2(0,R)}
		+\|(q_T-q_\infty)\psi_\infty^2\|_{L^1(0,R)}
		\Big).
	\end{aligned}
	\]
	Since $f_\lambda\in L^2(\R_+)$,
	$\|e^{-\eta\cdot}f_\lambda\|_{L^2(\R_+)}\to0$ as $\eta\to\infty$. Now, we
	choose $\eta>0$ large enough such that $2q_*\|G_m\|_{L^2(0,R)}
	\|e^{-\eta\cdot}f_\lambda\|_{L^2(\R_+)}<1.$
	Moreover, $q_T(s)\to q_\infty$ for each fixed
	$s\geq0$. Dominated convergence therefore gives $\|(q_T-q_\infty)\psi_\infty^2\|_{L^1(0,R)}
	\longrightarrow0.$
	Consequently, $\psi_T\longrightarrow\psi_\infty
	\;\text{in }L^2(0,R).$
	Finally, since by the  uniform majorant~\eqref{eq:UniformM} $|\mathbf1_{[0,T]}\psi_T-\psi_\infty|
	\leq 2G_m$
	and $G_m\in L^1(\R_+)\cap L^2(\R_+)$,
	a truncation argument yields $\mathbf1_{[0,T]}\psi_T
	\longrightarrow\psi_\infty
	\;
	\text{in }L^1(\R_+)\cap L^2(\R_+).$
	This proves~\eqref{eq:Stabil}.
	\qed
		\begin{Lemma}[Martingality]\label{lem:martingality}
		Let $T>0$ and $K\in L^2_{\mathrm{loc}}(\mathbb{R}_+)$, and let $X$ be the (unique in law) $\R_+$-valued solution on $[0,T]$ of equation~\eqref{eq:Volterrameanrevert} with diffusion coefficient $\sigma$ given by~\eqref{eq:sigma-definition} and mean-reversion speed $\lambda>0$.  Let $m\in\mathcal M^-_{\mathrm{lf}}(\mathbb{R}_+;\mathbb{C})$, and let $\psi=\psi(\cdot,m)\in L^2([0,T];\mathbb{C}_-)$ be the unique global solution of the measure-extended Riccati--Volterra equation~\eqref{eq:measureFLplce} given by Theorem~\ref{theorem:riccati-extension}. Define, for $t\in[0,T]$,
		\vspace{-.2cm}
		\begin{equation*}
			M_t := M_0\exp\Big(\int_0^t\psi(T-s)\,\sigma(s,X_s)\,dW_s - \frac12\int_0^t\psi^2(T-s)\,\sigma^2(s,X_s)\,ds\Big),\quad M_T=\exp\Big( \int_0^T X_{T-s} \,m( \mathrm{d}s) \Big).
			\vspace{-.2cm}
		\end{equation*}
		where $M_0$ is nonzero, $\mathcal F_0$-measurable,
		and chosen so that the terminal identity holds.  Then $(M_t)_{t\in[0,T]}$ is a true martingale.
	\end{Lemma}
	\noindent\textbf{\textit{Proof of Theorem~\ref{T:VolSqrt}.}} The proof closely follows the approach of \cite*[Theorem 4.3]{abi2019affine} for the convolution case and \cite*[Theorem 2.1]{ackermann2022} for the non-convolution and more general time-inhomogeneous setting, except that we work directly with the measure-extended conditional Fourier transform.
	
	\smallskip
	\noindent  {\sc Step~1} \textit{(The conditional Fourier transform of \(X_T\) ).} Let \( T > 0 \) and consider a measure $m \in \mathbb{M}^{-}_{\mathrm{lf}}(\R_+, \mathbb{C})$, by Theorem~\ref{theorem:riccati-extension} there exists a unique solution $\psi = 	\psi(\cdot,m) \in L^2([0,T];\mathbb{C}_-) $ to the measure-extended Riccati--Volterra equation ~\eqref{eq:measureFLplce}.
	Define 
	\vspace{-.2cm}
	\begin{align*}
		U_t &= 	\mathbb{E}_{}\left[ \int_0^{T} X_{T-s} \,m(ds) + \frac12 \int_t^T \sigma^2(s,X_s) \psi^2(T - s)\, ds \Big| \mathcal{F}_t \right]\\
		&= \int_{T-t}^{T} X_{T-s} \,m(ds) + \int_0^{T-t} \xi_t(T-s) \, m(\mathrm{d}s)+ \frac{1}{2}
		\int_t^T  \sigma^2(s,\xi_t(s))\, \psi^2(T - s)\, \mathrm{d}s.
	\end{align*}
	where $\sigma(\cdot,\cdot)$ is given by~\eqref{eq:sigma-definition} and the last equality comes from the affine nature of \(X\).
	Moreover, set \( M = \exp(U) \).
	\noindent Let \( X \) be a solution of equation~\eqref{eq:Volterrameanrevert}, with diffusion coefficient as in equation~\eqref{eq:sigma-definition}. 
	Then the process \((M_t )_{t\in[0,T]}\) is a
	local martingale on \( [0, T] \), and satisfies \(\frac{\mathrm{d}M_t}{M_t} = \psi(T - t)\, \sigma(t, X_t)\, \mathrm{d}W_t.\)
	In fact, by computing its dynamics using It\'o's formula, we can write:
	\vspace{-.3cm}
	\begin{align}\label{eq:MIto1D}
		\frac{d M_t}{M_t} = d U_t + \frac{1}{2} d\langle U \rangle_t.
		\vspace{-.3cm}
	\end{align}
	The dynamics of \( U \) can be obtained by recalling \( \xi_t(s) \) from \eqref{eq:DiffForward} and noting that for fixed \( s \), the dynamics of \( t \mapsto \xi_t(s) \) are given by  \(\mathrm{d}\xi_t(s) = \frac1\lambda f_\lambda(s - t)\, \sigma(t,X_t)\, \mathrm{d}W_t, \quad \text{for } t \le s.\)
	Since \( \xi_t(t) = X_t \), it follows that:
	\vspace{-.1cm}
	{\small
		\begin{align*}
			&\, dU_t =  \int_0^{T-t} d\xi_t(T-s) \, m(\mathrm{d}s) - \frac{1}{2} \sigma^2(t, X_t) \, \psi^2(T-t) \, dt 
			+ \frac{1}{2} \int_t^T \partial_x\sigma^2(s, \xi_t(s)) \, \psi^2(T-s) \, d\xi_t(s) \, ds \\[6pt]
			&= \int_0^{T-t} f_\lambda (T-s - t) \, m(\mathrm{d}s) \, \frac{\sigma(t, X_t)}{\lambda} \, dW_t
			- \frac{\sigma^2(t, X_t)}{2} \, \psi^2 (T - t) \, dt + \int_t^T k_1\, \varsigma^2(s) \, \psi^2 (T-s) \, f_\lambda (s - t) \, ds 
			\frac{\sigma(t, X_t)}{2\,\lambda} \, dW_t \\[6pt]
			&= -\frac{1}{2} \sigma^2(t, X_t) \, \psi^2 (T - t) \, dt 
			+ \Bigg( \int_0^{T-t} f_\lambda(T-t-s) \, m(\mathrm{d}s)+ \frac{1}{2} \int_t^T k_1\, \varsigma^2(s) \, \psi^2 (T - s) \, f_\lambda (s - t) \, ds \Bigg) 
			\frac{\sigma(t, X_t)}{\lambda} \, dW_t\\
			&= -\frac{1}{2} \sigma^2(t, X_t) \, \psi^2 (T - t) \, dt 
			+  \psi(T - t)\, \sigma(t, X_t)\, \mathrm{d}W_t.
		\end{align*}
	}
	\noindent where the last equality follows after a change of variable from the measure-extended Riccati--Volterra equation~\eqref{RicVolSqrt2} in Lemma \ref{lem:RiccatiWithF}. Thus, \(d\langle U \rangle_t =  \psi^2(T - t)\, \sigma^2(t, X_t) \, dt.\)
	Injecting the dynamics of \( dU_t \) and \( d\langle U \rangle_t \) into \eqref{eq:MIto1D}, we get \(	\frac{dM_t}{M_t} = \psi(T - t)\, \sigma(t, X_t)\, \mathrm{d}W_t.\)
	This shows that \( M \) is an exponential local martingale of the form  $M_t = M_0 \exp\Big( \int_0^t \psi(T - s)\, \sigma(s, X_s) \, dW_s - \frac{1}{2} \int_0^t \psi^2(T - s)\, \sigma^2(s, X_s)\, ds \Big).$
	\noindent To obtain \eqref{eq:laplace}, it suffices to prove that \( M \) is a martingale, which is given by Lemma~\ref{lem:martingality}. Consequently, the martingale property yields using that \(U_T = \int_0^T X_{T-s} \,m( \mathrm{d}s)\)
	\vspace{-.1cm}
	{\small
		\begin{align*}
			&\ \mathbb{E}_{}\left[ \exp\left( \int_0^T X_{T-s} \,m( \mathrm{d}s) \right)  \Big| \mathcal{F}_t \right]
			= \mathbb{E}_{}\left[ M_T \Big| \mathcal{F}_t \right] = M_t \\
			&\hspace{4.75cm}= \exp\left( \int_0^T  \mathbb{E}_{}[X_{T-s} \mid \mathcal{F}_t] \, m( \mathrm{d}s) + \frac{1}{2}
			\int_t^T \varsigma^2(s)\sigma^2( \mathbb{E}_{}[X_s \mid \mathcal{F}_t])\, \psi(T - s)^2\, \mathrm{d}s \right).
		\end{align*}
	}
	\noindent That is, as \( M \) is a true martingale, the measure-extended Fourier-Laplace transform of \( X_T \) is given by
		\vspace{-.2cm}
	{\small
		\begin{align*}\label{eq:resFourierLaplace}
			\mathbb{E}_{}\Big[ \exp\Big( \int_0^T  X_{T-s} \,m( \mathrm{d}s) \Big) \Big| \mathcal{F}_t \Big]
			= \exp\Big( \int_{T-t}^{T} X_{T-s} \,m(ds) + \int_0^{T-t} \xi_t(T-s) \,m( \mathrm{d}s) + \frac{1}{2}
			\int_t^T  \varsigma^2(s)\sigma^2(\xi_t(s))\, \psi(T - s)^2\, \mathrm{d}s \Big).
		\end{align*}
	}
	\noindent which yields \eqref{eq:laplace}.
	
	\smallskip
	\noindent  {\sc Step~2} \textit{(The Fourier-Laplace of \(X_T\) ).} 
	For any function \(g \in L^1_{\text{loc}}(\mathbb{R}; \mathbb{R})\), by regular Fubini's theorem and the measure-extended Riccati--Volterra equation~\eqref{RicVolSqrt2} in Lemma \ref{lem:RiccatiWithF}, we have that for all \(T>0\):
	{\small
		\begin{equation}\label{eq:DoubleConvol}  \int_0^T \int_0^{T-s}  f_{\lambda}(T-s-r)\, g(r) \, \mathrm{d}r\,m( \mathrm{d}s) 
			+ \int_0^T  \frac{\kappa_1}{2}\varsigma^2(s)\, \int_0^s  f_{\lambda}(s-r)\, g(r) \, \mathrm{d}r\,\psi(T-s)^2 \, \mathrm{d}s = \int_0^T \lambda \psi(T-r)\,g(r)\, \mathrm{d}r.\end{equation}
	}
	Evaluating the last equation of Step~1 
	 at $t = 0$, and using equation~\eqref{eq:DiffForward} for \(\xi_0(s) \;\forall \,s \leq T\), yields
	{\small
		\begin{align*}
			U_0 &=  \int_0^T \xi_0(T-s) \,m( \mathrm{d}s) + \frac{1}{2}
			\int_0^T  \varsigma^2(s)\sigma^2(\xi_0(s))\, \psi^2(T - s)\, \mathrm{d}s\\
			&=
			\int_0^T  x_0(T-s) \,m( \mathrm{d}s) +
			\int_0^T F(s, \psi(T-s))\, x_0(s) \, ds +
			\int_0^T \mu(s) \psi(T-s) \, ds + \frac{\kappa_0}{2} \int_0^T \varsigma^2(s) \psi^2(T-s) \, ds.
		\end{align*}
	} 
	\noindent where in the second line, we used the relation ~\eqref{eq:DoubleConvol} for both \(x_0\) and \(\mu\) and the measure-extended Riccati--Volterra equation~\eqref{eq:measureFLplce}. Note that the first equality in the claim follows directly by the fact that \(\forall s\geq0\,, \quad \xi_0(s) = \mathbb{E}_{} \left[X_s \big| X_0 \right] \).
	This provides the result. \qed
	
		\vspace{-.3cm}
	\subsection{Proof of Proposition~\ref{prop:Momentctrl_contraction_Hol} and Theorem~\ref{Theorem: stationary process} }\label{subsect:proofMresult2} 
	\noindent\textbf{\textit{Proof of Proposition~\ref{prop:Momentctrl_contraction_Hol}.}}
	\medskip
	\noindent {\bf (a) Moment Control.}
	Let $m_0 \in \mathbb{R}^*$ and $r>0$ a sufficiently large constant when $\gamma\in[1/2,1)$ ( if $\gamma=1$, take $r=1$). By the Hölder continuity of the state-dependent diffusion coefficient $\sigma$, together with \eqref{eq:subHolder}, it follows that for every $x\in\mathbb{R}$.
	\vspace{-.3cm}
	\[
	\sigma^2(x) \leq \Big(\sigma(m_0) + [\sigma]_{\text{H\"ol}} |x - m_0|^{\gamma}\Big)^2 \leq (1+\epsilon) |\sigma(m_0)|^2 + (1+\frac{1}{\epsilon})(1-\gamma)r^{\frac{1}{1-\gamma}} + (1+\frac{1}{\epsilon})\gamma\big(\frac{[\sigma]^2_{\text{H\"ol}}}{r}\big)^\frac{1}{\gamma} |x - m_0|^2\]
	 owing to the elementary inequality
	$(a+b)^2 \le (1+\epsilon) a^2+(1+\frac1\epsilon )b^2$ together with Young's inequality. Therefore, defining \(k_0 = k_0 (\epsilon):=(1+\epsilon) |\sigma(m_0)|^2 + (1+\frac{1}{\epsilon})(1-\gamma)r^{\frac{1}{1-\gamma}}\) and \(k_2= k_2(\epsilon):=(1+\frac{1}{\epsilon})\gamma\big(\frac{[\sigma]^2_{\text{H\"ol}}}{r}\big)^\frac{1}{\gamma}\), so that for \(r\) large enough \(ck_2(\epsilon)<1\), we can always assume, without loss of generality, that \textit{ $\sigma$ is sublinear}, namely,
		\vspace{-.3cm}
	\begin{equation}\label{eq:on_sigma}
		(SL_{\sigma}) : \exists k_0 = k_0 (\epsilon) \in \R_+, k_2= k_2(\epsilon) \in \R_+ \quad \textit{such that} \quad \forall x \in \mathbb{R}, |\sigma(x)|^2 \leq \kappa_0 + \kappa_2 |x -m_0|^2.  
	\end{equation}
	Using the reduced form~\eqref{eq:ConstMeanVar}, we have: $  \forall\, t\ge 0$,
	{\small
		\begin{align*}
			X_t - m_0 &= \big(X_0-m_0\big)h_{\lambda,c}(t) + \frac{1}{\lambda} \left(f_{\lambda} \stackrel{W}{*} \varsigma(\cdot)\sigma(X_{\cdot})\right)_t= \big(X_0-m_0\big)h_{\lambda,c}(t) + \frac{1}{\lambda} \int_0^t f_{\lambda}(t-s)\varsigma(s) \sigma(X_{ s}) \, dW_s.
		\end{align*}
	}
	Let $p\ge 2$. Define \( \rho_{p,r} := c\, (C_p^{BDG})^2\, \gamma\big(\frac{[\sigma]^2_{\text{H\"ol}}}{r}\big)^\frac{1}{\gamma} \) for sufficiently large \(r>0\) such that \(\rho_{p,r} \in(0,1)\). 
	Using the triangle inequality, the Burkholder--Davis--Gundy inequality applied to the local martingale
	$M_u = \int_0^u f_{\lambda}(t-s)\varsigma(s)\sigma(X_{s})dW_s$, $0\le s\le t$, (see~\cite[Proposition~4.3]{RevuzYor}) together with the generalized Minkowski inequality, yields
	\vspace{-.3cm}
	\begin{align*}
		\big\| |X_t-m_0|\big\|_p 
		&\le \big\|  |X_0-m_0|\big\|_p |h_{\lambda,c}(t)|  + \frac{C_p^{BDG}}{\lambda}\big\|\Big(f_{\lambda}^2 \stackrel{dt}{*} \varsigma^2(\cdot)|\sigma(X_{\cdot })|^2\Big)_t \big\|_{\frac p2}^{\frac 12}\\
		&\le  \big\|  |X_0-m_0|\big\|_p |h_{\lambda,c}(t)| +  \frac{C_p^{BDG}}{\lambda}\Big( \int_0^t f^2_{\lambda}(t-s)\varsigma^2(s) \big\||\sigma(X_{s})|^2\big\|_{\frac p2}\Big)^{\frac 12}.
	\end{align*}
	where equality holds in the first inequality when \(p=2\), in which case \(C_2^{\mathrm{BDG}}=1\). Fix now \(\epsilon\in(0,\rho_{p,r}^{-1}-1)\). Applying the elementary inequality $(a+b)^2 \le (1+\frac1\epsilon) a^2+(1+\epsilon )b^2\; $, we obtain
	\vspace{-.2cm}
	{\small
		\begin{align*}
			\Big\| |X_t-m_0|\Big\|_p^2 &\le \Big\|  |X_0-m_0|\Big\|_p^2 h^2_{\lambda,c}(t)(1+\frac1\epsilon)
			+ \frac{(C^{BDG}_p)^2}{\lambda^2}(1+\epsilon) \int_0^t f^2_{\lambda}(t-s)\varsigma^2(s) \||\sigma(X_{s})|^2\|_{\frac p2}ds
		\end{align*}
	}
	Likewise, introduce \( \widetilde{\rho}_{p,r} := c\, (C_p^{BDG})^2\, [\sigma]^2_{\text{H\"ol}}(1 + \varepsilon)= \rho_{p,r}(1 + \varepsilon) \in (0,\, 1)\) and  choose the parameter appearing in the sublinear growth condition \eqref{eq:on_sigma} as \( \frac{\widetilde{\rho}_{p,r}}{\eta}\) where \(\eta \in (0,\, 1 - \widetilde{\rho}_{p,r})\) is  a free parameter such that \( \widetilde{\rho}_{p,r} + \eta \in (0,\, 1) \).
	With this choice, the constants \( \kappa_i,\, i = 0, 2 \) in \eqref{eq:on_sigma} depending on \( \eta \) become
	\(k_0 = k_0 (\eta):=(1+\frac{\widetilde{\rho}_{p,r}}{\eta}) |\sigma(m_0)|^2+(1+\frac{\eta}{\widetilde{\rho}_{p,r}})(1-\gamma)r^{\frac{1}{1-\gamma}}\) and \(k_2= k_2(\eta):=\gamma\big(\frac{[\sigma]^2_{\text{H\"ol}}}{r}\big)^\frac{1}{\gamma}\) and satisfy \(c\, (C_p^{BDG})^2\,(1 + \varepsilon)\, \kappa_2 = \widetilde{\rho}_{p,r} + \eta < 1.\)
	Since $\frac p2\ge 1$, Jensen's inequality and \eqref{eq:on_sigma} imply that, for every \(s\geq0\), \(\;\big\||\sigma(X_{s})|^2\big\|_{\frac p2} \le \kappa_0 +\kappa_2\big\| |X_{ s}-m_0| \big\|^2_p\)
	which entails, combined with the identity  $f^2_{\lambda}*\varsigma^2 = c \lambda^2 (1-h^2_{\lambda,c})$,   that, for every $t\ge 0$, 
	{\small
		\begin{align}\nonumber 
			&\Big\| |X_t-m_0|\Big\|_p^2 \le  \Big\|  |X_0-m_0|\Big\|_p^2 h_{\lambda,c}^2(t) (1+\frac1\epsilon) \\ \label{eq:variancebound2}
			&\hspace{1.6cm}+ (C^{BDG}_p)^2(1+\epsilon)\Big( \kappa_0c\big( 1-h_{\lambda,c}^2(t) \big) + \frac{\kappa_2}{\lambda^2} \int_0^t f^2_{\lambda}(t-s) \varsigma^2(s)\Big\| |X_{ s}-m_0|\Big\|^2_p
			ds\Big).
		\end{align}
	}
	We now establish a uniform bound by a standard continuation argument. Let
	 $A> \bar A_{\eta,\epsilon,r}:=\frac{\kappa_0c(C^{BDG}_p)^2(1+\epsilon)}{1-\kappa_2c(C^{BDG}_p)^2(1+\epsilon)}\vee  \Big[ (1+1/\epsilon)\Big\|  |X_0-m_0|\Big\|_p^2\Big]$, $\delta>0$ and \(t_{\delta} = \inf\Big\{ t \geq 0: \Big\| |X_t-m_0|\Big\|_p^2 \ge A + \delta\Big\}.\)
	 Assume, for contradiction, that \(t_\delta<\infty\). By the continuity of $t\mapsto\Big\| |X_t-m_0|\Big\|_p^2$, we have $A+\delta =\Big\| |X_{ t_\delta}-m_0|\Big\|_p^2$. Moreover, for every \(t\leq t_\delta\), the Volterra identity for \(\varsigma\) gives and $h_{\lambda,c}$ gives
	\(\int_0^t f^2_{\lambda}(t-s) \varsigma^2(s)\Big\| |X_{s}-m_0|\Big\|_p^2ds\le  (A+\delta)c \lambda^2(1-h_{\lambda,c}^2(t)). \) Evaluating \eqref{eq:variancebound2} at \(t=t_\delta\), and using
	 \( A > \left\| X_0 - m_0\right\|_p^2 (1 + \frac{1}{\epsilon}) \), we obtain the inequalities:
	\begin{align*}
		A+\delta &= \Big\| |X_{ t_\delta}-m_0|\Big\|_p^2< A h^2_{\lambda,c}(t_\delta)+(C^{BDG}_p)^2(1+\epsilon)\big(\kappa_0 c +\kappa_2c (A+\delta)\big)\big(1-h^2_{\lambda,c}(t_\delta)\big)\\
		&< A h^2_{\lambda,c}(t_\delta) +  A(1-h^2_{\lambda,c}(t_\delta)) +   (C^{BDG}_p)^2(1+\epsilon)c\kappa_2 \delta\big((1-h^2_{\lambda,c}(t_\delta)\big)
		 < A+ \delta.
	\end{align*}
	where we used the bound  \( (C^{\mathrm{BDG}}_p)^2 (1 + \epsilon)\, c\, (\kappa_0 + \kappa_2 A) < A \)
	which holds by the definition of \( A \), whereas the last inequality follows from the fact that $ (C^{BDG}_p)^2(1+\epsilon)c\kappa_2<1$. This  yields a contradiction. Consequently, $t_{\delta}= +\infty$ and therefore $\Big\|  |X_{ t}-m_0|\Big\|_p^2  \le A+\delta$ for every $t\ge 0$. Letting $\delta \to 0$ and 
	$A\to \bar A_{\eta,\epsilon,r}$ successively, yields 
	\begin{equation}
		\sup_{t\ge 0}\Big\|  |X_t-m_0|\Big\|_p\leq \bar A^{\frac1 2}_{\eta,\epsilon,r}= (1+\frac1\epsilon)\Big\|  |X_0-m_0|\Big\|_p \vee \Big(\frac1\gamma\big(\frac{[\sigma]^2_{\text{H\"ol}}}{r}\big)^\frac{-1}{\gamma}\Big((1-\gamma)r^{\frac{1}{1-\gamma}}+\frac{\widetilde{\rho}_{p,r}|\sigma(m_0)|^2}{\eta}\Big)\frac{\eta+\widetilde{\rho}_{p,r}}{(1-\widetilde{\rho}_{p,r}-\eta)} \Big)^{\frac1 2},
	\end{equation}
	\noindent which is finite. It remains to optimize the bound with respect to the free parameter \(\eta\). Differentiating the right-hand side with respect to \(\eta\) shows that the minimum over
	\((0,1-\widetilde\rho_{p,r})\) of the mapping 
	\(\eta \mapsto \bar{A}_{\eta,\epsilon,r}\) is attained at
	\[\eta^* =\big((1-\gamma)r^{\frac{1}{1-\gamma}}+\widetilde{\rho}_{p,r}|\sigma(m_0)|^2\big)^{-1}\Big(-\widetilde{\rho}_{p,r}^2|\sigma(m_0)|^2+ \widetilde{\rho}_{p,r}|\sigma(m_0)|(\widetilde{\rho}_{p,r}|\sigma(m_0)|^2+(1-\gamma)r^{\frac{1}{1-\gamma}}(1-\widetilde{\rho}_{p,r}))^{\frac12}\Big).\] Substituting this value yields the minimum 
	\begin{equation}A^{\frac1 2}_{\eta^*,\epsilon,r}= (1+\frac1\epsilon)\Big\|  |X_0-m_0|\Big\|_p \vee \frac{\Big(\widetilde{\rho}_{p,r}|\sigma(m_0)|^2+|\sigma(m_0)|\big(\widetilde{\rho}_{p,r}|\sigma(m_0)|^2+(1-\gamma)r^{\frac{1}{1-\gamma}}(1-\widetilde{\rho}_{p,r})\big)^{\frac12}\Big)}{\gamma\big(\frac{[\sigma]^2_{\text{H\"ol}}}{r}\big)^\frac{1}{\gamma}(1-\widetilde{\rho}_{p,r})\big((1-\gamma)r^{\frac{1}{1-\gamma}}+\widetilde{\rho}_{p,r}|\sigma(m_0)|^2\big)^{\frac12}}
	\end{equation}	
	This completes the proof. 
	The stated result follows upon setting \(C^{\mathrm{BDG}}_p=2\sqrt{p}\), as given in \cite[Remark 2]{carlen1991lp}. 

	\medskip
	\noindent\textbf{(b) \(L^2\)-contraction.} In contrast to the proof of~\cite[Proposition~4.4]{EGnabeyeu2025}, which is based on a Banach fixed-point argument in the Lipschitz case (\(\gamma=1\)) for every \(L^p\) space with \(p>0\), we instead rely on It\^o's isometry and the first Dini lemma in this \(L^2\) setting.
	Let \((X_t)_{t\ge0}\) and \((X'_t)_{t\ge0}\) be two solutions of the Volterra stochastic equation~\eqref{eq:Volterrameanrevert} with initial conditions \(X_0,X'_0\in L^2(\mathbb P)\).
	Owing to Equation~\eqref{eq:subHolder} and the concavity of \(x \to x^\gamma\) for all \(\gamma \in [\frac12,1]\), from Equation~\eqref{eq:subHolder}, we may write
	\vspace{-.2cm}
	\begin{equation}\label{eq:subHolder2} 
		\E \left[|\sigma(X_{s}) - \sigma(X^\prime_{s})|^2\right]  
		\leq [\sigma]^2_{\text{H\"ol}} \E \left[| X_{s} - X^\prime_{s} |^2\right]^{\gamma}.
	\end{equation}
	Set $\Delta_t = X_t - X^\prime_t \in L^2 (\P)$ for every $t\ge 0$. One writes owing to the reduced form~\eqref{eq:ConstMeanVar}
	\begin{equation}
		X_t-X^\prime_t = \big(X_0 - X^\prime_0\big) h_{\lambda,c}(t) + \frac{1}{\lambda}\int_0^tf_{\lambda}(t-s)\varsigma(s) \big(\sigma(X_{s}) -\sigma(X^\prime_{s})\big)dW_s 
	\end{equation}
	Let \( \bar \delta_t = \Big\| |\Delta_t|\Big\|_2 \) for convenience. One checks that, \( t \mapsto \bar \delta_t \) is continuous. 
	Set $\rho = c[\sigma]^2_{\text{H\"ol}} \in (0, 1)$. Elementary computations and It\^o's Isometry show that for every $t\ge 0$
	\begin{align*}
		\E\, \big(\big| X_t - X^\prime_t\big|^2\big) = h^2_{\lambda,c}(t)\E\, \big(\big|X_0-X^\prime_0\big|^2\Big)+ \frac{1}{\lambda^2} \int_0^t f^2_{\lambda}(t-s) \varsigma^2(s) \E\, \big(|\sigma(X_{s}) -\sigma(X^\prime_{s})|^2\big) ds.
	\end{align*}
	which entails owing to equation ~\eqref{eq:subHolder2} :
	
	\begin{equation}\label{eq:l2_confl}
		\big\| |\Delta_t| \big\|_2^2 \le \big\| |\Delta_0| \big\|_2^2 h^2_{\lambda,c}(t) 
		+ \frac{[\sigma]^2_{\text{H\"ol}}}{\lambda^2} \int_0^t f_{\lambda}^2(t - s) \varsigma^2(s) 
		\big\| |\Delta_{s}| \big\|_2^{2\gamma} \, ds.
	\end{equation}
	From now on, we define $\Psi : x \to \left\{
	\begin{array}{ll}
		& x \text{ if } \gamma = 1 \\
		& 1 \vee x \text{ if } \gamma \in [\frac12,1[ 
	\end{array} .\right.$
	One checks that, \(\forall\; \alpha, x\in \R_+,\; \Psi^\alpha(x)= \Psi(x^\alpha) \), \(\Psi(x)\geq x\) and \(\Psi^{2\gamma}(x) \leq \Psi^{2}(x) = \Psi(x^2).\) (Note that \(\left(1\vee\bar\delta_0\right) \geq 1\) so that \(\left(1\vee\bar\delta_0\right)^{2\gamma} \leq \left(1\vee\bar\delta_0\right)^2\) as \(\gamma\le 1\))
	
	\noindent  {\sc Step~1} \textit{ Non-expansivity via a stopping-time argument:} Let \( \eta > 0 \) such that \( \rho(1 + \eta)^2 < 1 \). We define \( \tau_{\eta} := \inf\{ t \geq 0 : \bar \delta_t > (1 + \eta)  \Psi(\bar\delta_0) \} \).
	If \( \tau_{\eta} < +\infty \), we have \( \bar \delta_s \le (1 + \eta) \Psi(\bar\delta_0) \) for \( s \in (0, \tau_{\eta}) \) and by continuity, \( \bar \delta_{\tau_{\eta}}^2 = (1 + \eta)^2 \Psi(\bar\delta_0^2) \). Plugging this into equation ~\eqref{eq:l2_confl}, the inequality at time \( \tau_{\eta} \) combined with the identity  $f^2_{\lambda}*\varsigma^2 = c \lambda^2 (1-h^2_{\lambda,c})$ yields:
	\begin{align*}
		\bar \delta_{\tau_{\eta}}^2 &\le \Psi(\bar\delta_0^2) \left[ h^2_{\lambda,c}(\tau_{\eta}) + (1 - h^2_{\lambda,c}(\tau_{\eta})) \rho (1 + \eta)^{2\gamma} \right] \\
		&\leq\Psi(\bar\delta_0^2) 
		<  (1 + \eta)^2  \Psi(\bar\delta_0^2),
	\end{align*}
	which leads to a contradiction. Therefore  \( \tau_{\eta}  = +\infty \) i.e., \( \bar \delta_s \le (1 + \eta) \Psi(\bar\delta_0) \) for all \( s \ge 0 \). This holds for every \( \eta > 0 \), implying the non-expansive bound \( \bar \delta_t \le \Psi(\bar\delta_0) \) for all \( t \ge 0 \) when passing to the limit as \( \eta \) to 0. 
	
	\noindent  {\sc Step~2} \textit{Iteration and the Volterra map:}	
	Substituting this (i.e. \( \bar \delta_t \le \Psi(\bar\delta_0) \)) into~\eqref{eq:l2_confl} combined with the identity  $f^2_{\lambda}*\varsigma^2 = c \lambda^2 (1-h^2_{\lambda,c})$ gives, for all \( t > 0 \), 
	\[
	\bar\delta_t^2 \le \Psi(\bar\delta_0^2) \varphi_1(t), \quad \text{where} \quad \varphi_1(t) := h^2_{\lambda,c}(t) + (1 - h^2_{\lambda,c}(t)) \rho, \;\quad \lim_{t \to \infty}\varphi_1(t) =\ell_\infty^2 +(1-\ell_\infty^2)\rho.
	\]
	Note that \( \varphi_1(t) = \rho + h^2_{\lambda,c}(t)(1 - \rho) \) satisfies:
	\(\varphi_1(0) = 1, \; \varphi_1(t) \in [0, 1] \, \forall \, t >0, \; \varphi_1 \text{ is 
		continuous.}\)
	Substituting this upper bound (i.e. $\bar\delta_t^2 \le \Psi(\bar\delta_0^2) \varphi_1(t)$) into~\eqref{eq:l2_confl} yields 
	\[
	\bar\delta_t^2 \le \Psi(\bar\delta_0^2) \varphi_2(t), \quad \text{where} \quad \varphi_2(t) := h^2_{\lambda,c}(t) + \rho \int_0^t f_{\lambda}^2(t - s) \varsigma^2(s) \varphi_1^{\gamma}(s) \frac{ds}{\lambda^2 c}.
	\]
	Using the identity $f^2_{\lambda} * \varsigma^2 = c \lambda^2 (1 - h^2_{\lambda,c})$ and the definition of $\varphi_1$, we get
	\[
	\varphi_2(t) := \varphi_1(t) - \rho \int_0^t f^2_{\lambda}(t-s)\varsigma^2(s) \left(1 - \varphi_1^{\gamma}(s)\right) \frac{ds}{\lambda^2 c}.
	\]
	Thus, $0 \leq \varphi_2 < \varphi_1 < 1$ on $(0, +\infty)$. By induction, we show that \(\bar\delta^2_t \leq \Psi(\bar\delta_0^2) \varphi_k(t)\)
	with
	\begin{equation}\label{eq:IterVolteraMap}
	\varphi_k(t) = h^2_{\lambda,c}(t) + \rho \int_0^t f^2_{\lambda}(t-s) \varsigma^2(s) \varphi_{k-1}^{\gamma}(s) \frac{ds}{\lambda^2 c} = \varphi_1(t) - \rho \int_0^t f^2_{\lambda}(t-s) \varsigma^2(s) \left(1 - \varphi_{k-1}^{\gamma}(s)\right) \frac{ds}{\lambda^2 c}.
	\end{equation}
	where we used again the identity \eqref{eq:MasterEq} satisfied by $\varsigma^2$ and the definition of $\varphi_1$.
	
	\noindent  {\sc Step~3} \textit{ Monotonicity and Limit equation:} Consequently, starting from \( 0 \leq \varphi_2 < \varphi_1 < 1 \) on \( (0, +\infty) \), induction shows that $0 \leq \varphi_k < \varphi_{k-1} < 1$ on \( (0, +\infty) \) for every \( k \geq 2 \). Furthermore, we verify by induction that $\varphi_k$ is continuous, as by change of variable,
	\[
	\varphi_k(t) = \varphi_1(t) - \rho \int_0^t f^2_{\lambda}(s) \varsigma^2(t-s) \left(1 - \varphi_{k-1}^{\gamma}(t-s)\right) \frac{ds}{\lambda^2 c},
	\]
	where $\varsigma^2$ is continuous and bounded in \(\R_+^*\) by assumption.
	
	\noindent By the first Dini Lemma, we have $\varphi_k \downarrow \varphi_\infty \in \mathcal{C}(\mathbb{R}_+, \mathbb{R})$ uniformly on compact intervals with $\varphi_\infty(0) = 1$. In particular, $\varphi_\infty$ satisfies the functional equation
	\begin{equation}\label{eq:func}
		\varphi_\infty(t) = h^2_{\lambda,c}(t) + \rho \int_0^t f^2_{\lambda}(t-s) \varsigma^2(s) \varphi_\infty^{\gamma}(s) \frac{ds}{\lambda^2 c}.
	\end{equation}
	
	\noindent  {\sc Step~4} \textit{ Limit and asymptotic bound for $\ell_\infty \in (0,1]$:} Let $\underline{\ell}_\varphi, \overline{\ell}_\varphi := \liminf_{t \to +\infty} \varphi_\infty(t), \limsup_{t \to +\infty} \varphi_\infty(t) \in [0,1]$.  Now, $\underline{\ell}_\varphi, \overline{\ell}_\varphi \in [0,1]$ implies that for any $\eta > 0$, there exists $t_\eta \in \R^+$ such that for $t \geq t_\eta$, $ \underline{\ell}_\varphi - \eta \leq \varphi_\infty(t)\leq \overline{\ell}_\varphi + \eta $.
	On the first hand,
	\begin{align*}
		\int_0^t f^2_{\lambda}(t-s) \varsigma^2(s) \varphi_\infty^{\gamma}(s) \frac{ds}{\lambda^2 c} &\leq \frac{1}{c\lambda^2} \int_{t_\eta}^t f^2_{\lambda}(t-s) \varsigma^2(s) (\overline{\ell}_\varphi + \eta)^{\gamma} \, ds + \frac{1}{c\lambda^2} \int_{0 }^{t_\eta} f^2_{\lambda}(t-s) \varsigma^2(s)\varphi_\infty^{\gamma}(s)\, ds\\
		&\leq \frac{1}{c\lambda^2} \int_{t_\eta}^t f^2_{\lambda}(t-s) \varsigma^2(s) (\overline{\ell}_\varphi + \eta)^{\gamma} \, ds + \frac{1}{c\lambda^2} \int_{t - t_\eta}^t f^2_{\lambda}(u)\varsigma^2(t-u) \, du.
	\end{align*}
	where the second term on the right-hand side of the last inequality follows from the fact that \( \varphi_\infty(t - u) \leq 1\) for all \(u \leq t \leq t_\eta\) and vanishes as t goes to infinity. Since $f_\lambda \in L^2(\text{Leb}_1)$ and $\lim_{t \to +\infty} h^2_{\lambda,c}(t) = \ell^2_\infty$ both owing to Assumption \ref{ass:resolvent}, we conclude from  equation~\eqref{eq:func} and the identity $f^2_{\lambda}*\varsigma^2 = c \lambda^2 (1-h^2_{\lambda,c})$ :
	\[
	\overline{\ell}_\varphi =:\limsup_{t \to +\infty} \varphi_\infty(t)  \leq \ell^2_\infty + \rho(\overline{\ell}_\varphi + \eta)^{\gamma}(1-\ell^2_\infty),
	\]
	which implies $\overline{\ell}_\varphi \leq \ell^2_\infty +\rho(1-\ell^2_\infty)\overline{\ell}_\varphi^\gamma \quad \textit{by letting} \quad \eta \to 0.$
	On the other hand, 
	we also have:
	\begin{align*}
		\int_0^t f^2_{\lambda}(t-s) \varsigma^2(s) \varphi_\infty^{\gamma}(s) \frac{ds}{\lambda^2 c} &\geq \frac{1}{c\lambda^2} \int_{t_\eta}^t f^2_{\lambda}(t-s) \varsigma^2(s) (\underline{\ell}_\varphi - \eta)_+^{\gamma} \, ds + \frac{1}{c\lambda^2} \int_{t - t_\eta}^t f^2_{\lambda}(u) \, \varsigma^2(t-u)\varphi_\infty^{\gamma}(t-u)\, du\\
		&\geq \frac{1}{c\lambda^2} \int_{t_\eta}^t f^2_{\lambda}(t-s) \varsigma^2(s) (\underline{\ell}_\varphi - \eta)_+^{\gamma} \, ds.
	\end{align*}
	where $x_+:=\max\{x,0\}$. Therefore, still with the fact that $f_\lambda \in L^2(\text{Leb}_1)$ and $\lim_{t \to +\infty} h^2_{\lambda,c}(t) = \ell^2_\infty$, we obtain from equation~\eqref{eq:func} and the identity $f^2_{\lambda}*\varsigma^2 = c \lambda^2 (1-h^2_{\lambda,c})$ :
	\[
	\underline{\ell}_\varphi =:\liminf_{t \to +\infty} \varphi_\infty(t)  \geq \ell^2_\infty + \rho(\underline{\ell}_\varphi - \eta)^{\gamma}(1-\ell^2_\infty) \quad \implies \underline{\ell}_\varphi \geq \ell^2_\infty+\rho(1-\ell^2_\infty)\underline{\ell}_\varphi^\gamma  \quad \textit{as} \quad \eta \to 0.
	\] The function $G(x):=\ell_\infty^2+\rho(1-\ell_\infty^2)x^\gamma-x$
	admits a unique zero \(\ell_\varphi\in[0,1]\), with \(G(x)\geq0\) for \(x\leq\ell_\varphi\) and \(G(x)\leq0\) for \(x\geq\ell_\varphi\). We conclude that $\underline{\ell}_\varphi
	=
	\overline{\ell}_\varphi
	=
	\ell_\varphi.$ Now,
	if \(\gamma=1,\) then \(\ell_\varphi = \frac{\ell^2_\infty}{1-\rho(1-\ell^2_\infty)}:=\overline{\ell}_{\varphi}^\star\). Furthermore, if 
	\(\gamma=\frac12,\) then \(\ell_\varphi=\underline{\ell}_{\varphi}^\star\) defined in the proposition
	and if \(\gamma \in ]\frac12,1[\), as \(\ell_\varphi \in (0,1)\), we have \(\ell_\varphi \leq \ell_\varphi^\gamma\leq \sqrt{\ell_\varphi}\) yielding the inequality $\overline{\ell}_{\varphi}^\star\leq \ell_\varphi \leq\underline{\ell}_{\varphi}^\star$.
	
	\smallskip
	\noindent  {\sc Step~5} \textit{Case \(\ell_\infty=0\).} When \(\ell_\infty=0\), \(\ell_\varphi\) is a fixed point of the function \(x\to\rho x^\gamma\) i.e. \(\ell_\varphi\) is either \(0\) or \(\rho^{\frac{1}{1-\gamma}} \). Note that the only fixed point is \(0\) when \(\gamma =1\). If \(\gamma \in [\frac12,1[\), \(\ell_\varphi \in \big\{0,\rho^{\frac{1}{1-\gamma}}\big\}.\) Moreover, since the sequence \((\varphi_k)_{k\geq1}\) is non-increasing and \(\forall\,k\geq2,\;\lim_{t \to \infty}\varphi_k(t)\leq \lim_{t \to \infty}\varphi_1(t) =\rho \). 
	We now show that the positive fixed point is attained.
	Passing
	\eqref{eq:IterVolteraMap} to the limit
	as $t\to\infty$, we obtain the scalar recursion \(L_k = \rho L_{k-1}^{\gamma}\) where \(L_k:=\lim_{t\to\infty}\varphi_k(t).\)
	Since \(L_1=\rho>0\), it follows by induction that \(L_k= \rho^{1+\gamma+\cdots+\gamma^{k-1}}
	\downarrow \rho^{\sum_{j=0}^{\infty}\gamma^j}
	=\rho^{\frac{1}{1-\gamma}}>0.
	\) 
	To transfer this lower bound to the limiting function
	$\varphi_\infty$, using~\eqref{eq:IterVolteraMap}
	and the master equation, one checks by induction that
	$\varphi_k(t)\geq \rho^{\frac{1}{1-\gamma}}$ for every $k\geq1$
	and $t\geq0$. Hence, letting $k\to\infty$
	pointwise yields $\varphi_\infty(t)\geq \rho^{\frac{1}{1-\gamma}}$ for every $t\geq0$. Consequently, \(\ell_\varphi\geq\rho^{1/(1-\gamma)}>0\).
Combining this lower bound with the fixed-point equation yields \(\ell_\varphi=\rho^{\frac{1}{1-\gamma}},
	\; \gamma\in[\frac12,1).\) 
	\noindent This completes the proof and we are done. \qed
	
	\noindent {\bf Remark:} The function \( \varphi_\infty\) quantifies the decay over time of the expected squared difference between two solutions of the SVIE~\eqref{eq:ConstMeanVar} with different initial values. 
	As \(\varsigma\) is bounded, if both \([\sigma]^2_{\text{H\"ol}} < \frac{\lambda^2}{\| \varsigma^2 \|_\infty \int_0^{+\infty} f_\lambda^2(u) \, du}\) and \(h_{\lambda,c}  \in L^2(\text{Leb}_1)\), 
	then one derives from equation~\eqref{eq:func} and using Fubini-Tonelli's theorem that:
	
	\centerline{$
		\int_0^{+\infty} \varphi_\infty(s) \, ds \leq \frac{\lambda^2 \int_0^{+\infty} h_{\lambda,c}^2(t) \, dt}{\lambda^2 - [\sigma]^2_{\text{H\"ol}} \| \varsigma^2 \|_\infty \int_0^{+\infty} f_\lambda^2(u) \, du} < +\infty.
		$}

	\bigskip
		\noindent\textbf{\textit{Proof of Theorem~\ref{Theorem: stationary process}.}}
	\noindent Let $p\!\in [2,+\infty)$ be such that $p(\vartheta\wedge \frac{\beta-1}{2\beta})>1$. Hence, Proposition \ref{prop:Momentctrl_contraction_Hol} implies that $ \sup_{t\ge 0}\big\|  |X_t-m_0|\big\|_p < +\infty$.
	As a consequence of the linear growth of \(\sigma\), ( there is a constant \(C_{\sigma} > 0\) such that $|\sigma(x)| \leq C_{\sigma}(1 + |x|)$,
	for all  $t \in [0, T]$ and $ x \in \mathbb{R}$), we derive that
	 \( \sup_{t \geq 0} \big\| |\sigma(X_t)|^2 \big\|_p < +\infty \).\\
	\noindent  {\sc Step~1} \textit{ (C-tightness of time-shifted processes and weak convergence ).} 
	Now, we can argue the C-tightness by the Kolmogorov criterion. One writes  for $s,\,t\ge 0$ with $s\le t$ and owing to equation~\ref{eq:ConstMean}:
	\begin{equation*}\label{eq:Tight} X_t-X_s = \big(X_0 - m_0\big) \big( h_{\lambda,c}(t) - h_{\lambda,c}(s) \big) + \frac{1}{\lambda}\big(I(t)-I(s)\big).
	\end{equation*}
	where we set $I(t):=\int_0^tf_{\lambda}(t-u)\varsigma(u)\sigma(X_{u})dW_u$.
	First, using the master equation~\eqref{eq:MasterEq}
	and the fact that $h_{\lambda,c}\geq0$ together with the elementary inequality $|\sqrt{a}-\sqrt{b}|\leq \sqrt{|a-b|}$ for every \(a,b\geq0\),
	we obtain on the first hand, for \(0\leq s<t\)
	\begin{equation*}
		\begin{aligned}	&\, \Big\| \big(X_0 - m_0\big)  \big( h_{\lambda,c}(t)-h_{\lambda,c}(s) \big) \Big\|_p \le \big\| X_0 - m_0 \big\|_p \frac{\|\varsigma\|_{\infty}}{\lambda\sqrt{c}}  \Big(\int_s^{t}f_{\lambda}^2(t-u)\, du  +  \int_0^s \big| f_{\lambda}^2(t-s + u) - f^2_{\lambda}(u) \big|\, du \Big)^\frac12  \\
			&\hspace{.8cm}\le \big\| X_0 - m_0 \big\|_p \frac{\|\varsigma\|_{\infty}}{\lambda\sqrt{c}}    \times \Big( \big(\int_0^{+\infty} f_{\lambda}^{2\beta}(u)du \big)^{\frac{1}{\beta}}  (t-s)^{1-\frac{1}{\beta}} + \int_0^{t-s} f_\lambda(u)^2\,du
			-\int_s^{t}f_\lambda(u)^2\,du \Big)^\frac12 \\
			&\hspace{.8cm}\le \big\| X_0 - m_0 \big\|_p 	\frac{\sqrt{2}\|\varsigma\|_\infty}
			{\lambda\sqrt c}
			\|f_\lambda\|_{L^{2\beta}(\R_+)} \cdot|t-s|^{\frac12(1-\frac{1}{\beta})}.
		\end{aligned}
		\vspace{-.1cm}
	\end{equation*}
Here, the penultimate inequality follows from the complete monotonicity of \(f_\lambda\), which implies that \(f_\lambda^2\) is nonnegative and nonincreasing, while the last inequality follows from Hölder's inequality.
	On the other hand, combining the $L^p$-BDG and the generalized Minkowski inequality, one derives that,
	\vspace{-.1cm}
	\begin{equation*}
		\begin{aligned}
			&\ \big\| J(t)-J(s) \big\|_p \le \Big\| \int_s^{t}f_{\lambda}(t-u)\varsigma(u)\sigma(X_{u})\mathrm{d}W_u \Big\|_p + \Big\| \int_0^s \big( f_{\lambda}(t-u) - f_{\lambda}(s-u) \big)\varsigma(u)\sigma(X_{u})dW_u  \Big\|_p\\
			&\leq C^{BDG}_p \|\varsigma\|_{\infty} \sup_{u\ge 0}\big\| |\sigma^2(X_{u})| \big\|_p^{\frac12} \Big[\Big(\int_0^{+\infty}f^{2\beta}_{\lambda}(u) du\Big)^{\frac{1}{2\beta}} |t-s|^{\frac{\beta-1}{2\beta}} + \Big(\int_0^{+\infty} \big(f_{\lambda}(t-s+u)-f_{\lambda}(u)\big)^2du\Big)^{\frac12}\Big] \\
		\end{aligned}
		\vspace{-.1cm}
	\end{equation*}
	Finally, using Assumption~\ref{ass:int_holregul} and putting all these estimates together, since $ \frac{\beta-1}{2\beta}\le \frac{\beta-1}{\beta} $ we have the existence of a real constant $0<C_{p, X_0, \varsigma, \beta,\lambda, f_{\lambda}}<\infty$ such that:
	\vspace{-.3cm}
	\begin{equation}\label{eq:cont}\E\, \left(|X_t-X_s| \right)^p \le C_{p, X_0, \varsigma, \beta,\lambda, f_{\lambda}}|t-s|^{p(\vartheta\wedge \frac{\beta-1}{2\beta})}\quad \text{for every} \quad t,s\geq0\quad \text{with } 0 \leq t - s \leq 1.
		\vspace{-.2cm}
	\end{equation}
	And thus, \( t \mapsto X_t \) admits a H\"older continuous modification (still denoted \( X \) in lieu of \( \tilde{X} \) up to a \( \mathbb{P} \)-indistinguishability), so that the process has the announced H\"older pathwise regularity, i.e. \( t \mapsto X_t \) has a \( \big(\vartheta\wedge \frac{\beta-1}{2\beta}-\frac1p - \eta \big) \)-H\"older pathwise continuous \( \mathbb{P} \)-modification for sufficiently small \( \eta > 0 \).
	
	Define for \( u \geq 0 \) the process \( X^u \) by \( X^u_t = X_{t+u} \), where \( t \geq 0 \). Then \( X^u \) has continuous sample paths and satisfies for a constant $C(p)>0$ (this follows from equation~\eqref{eq:cont}):
	\[
	\sup_{u \geq 0} \mathbb{E}[|X^u_t - X^u_s|^p] \leq C(p)|t-s|^{p(\vartheta\wedge \frac{\beta-1}{2\beta})}
	\quad \text{for all} \quad t,s\geq0\quad \text{with } 0 \leq t - s \leq 1.\]
	Now, as $p\geq2$ is chosen sufficiently large so that $p(\vartheta\wedge \frac{\beta-1}{2\beta})>1$, it follows from  Kolmogorov's $C$-tightness criterion (see ~\cite[Theorem 2.1, p. 26, 3rd edition]{RevuzYor} 
	or \cite[Lemma 44.4, Section IV.44, p.100]{RogersWilliamsII}), that the family of  shifted processes $X_{t+\cdot}$, $t\ge 0$, is $C$-tight i.e. \( (X^u)_{u \geq 0} \) is tight on \( \mathcal{C}(\mathbb{R}_+; \mathbb{R}) \) (hence the existence of a weak
	continuous accumulation point thanks to Prokhorov's theorem) with limiting distributions P under which the canonical process has the announced H\"older pathwise regularity. 
	
	\smallskip
	\noindent  {\sc Step~2} \textit{(Moment and autocovariance function ).} For the first moment, equation~\eqref{eq:Forward2}, together with
	Assumptions~\ref{ass:resolvent}--\ref{ass:Mean_revertF} and the dominated
	convergence theorem, yields (see also~\cite[Lemma 3.1]{EGnabeyeu2025})
		\[
	\mathbb{E}[X_t] =  \mathbb{E}[X_0](\varphi(t)- \int_0^t f_{\lambda}(t-r) \varphi(r) \, dr ) + \frac{1}{\lambda} \int_0^t f_{\lambda}(t-r) \mu(r) \, \mathrm{d}r \longrightarrow \ell_\infty  \mathbb{E}[X_0]  + (1-a) \frac{\mu_\infty}{\lambda} \quad \text{as } t \to \infty.
	\]
	On the other hand, for every \(t\geq0\), $X_{t+u}\Longrightarrow X_t^\infty,
	\; u\to\infty,$
	and the uniform \(L^2\)-bound $\sup_{s\geq0}\mathbb E[|X_s|^2]<\infty$
	implies that \((X_{t+u})_{u\geq0}\) is uniformly integrable. Hence
	\[
	\mathbb E[X_t^\infty]
	=
	\lim_{u\to\infty}\mathbb E[X_{t+u}]
	=
	\ell_\infty\mathbb E[X_0]
	+
	(1-a)\frac{\mu_\infty}{\lambda},
	\qquad t\geq0.
	\]
	Moreover, by Theorem~\ref{Theorem: limiting distribution}, $\int_{\mathbb R_+}x\,\pi_{X_0}(dx)=\Xi_0,
	\; \mathbb P\text{-a.s.},$
	so that $	\mathbb E[X_t^\infty]
	=
	\mathbb E[\Xi_0]
	=
	\mathbb E\!\left[
	\int_{\mathbb R_+}x\,\pi_{X_0}(dx)
	\right].$
	This proves the desired first moment formula for the asymptotic process.
	
	\noindent (a) Now let us consider the asymptotic  covariance between $X_{t+t_1}$ and $X_{t+t_2}$, $0<t_1<t_2$ when $X_t$ starts for $X_0$, $t\ge 0$ with \(v_0:={\rm Var}(X_0)\).
	Noting from equation ~\eqref{eq:Volterrameanrevert2} that :
	\begin{align*}
		& \; X_{t_2} - \E[X_{t_2}]
		= (X_{0} - \E[X_{0}])h_{\lambda,c}(t_2)
		+ \frac{1}{\lambda} \Big(
		\int_0^{t_1} f_{\lambda}(t_2-s)\varsigma(s)\sigma(X_{s})\,dW_s
		+ \int_{t_1}^{t_2} f_{\lambda}(t_2-s)\varsigma(s)\sigma(X_{s})\,dW_s
		\Big), \\
		& \; \text{and} \; \qquad X_{t_1} - \E[X_{t_1}]
		= (X_{0} - \E[X_{0}])h_{\lambda,c}(t_1)
		+ \frac{1}{\lambda}\int_0^{t_1} f_{\lambda}(t_1-s)\varsigma(s)\sigma(X_{s})\,dW_s.
	\end{align*}
	Using $\text{Cov}(aU + b, cV + d) = ac \, \text{Cov}(U, V)$, we find that the autocovariance function for $X$ is given by:
	\begin{align*}
		&{\rm Cov}(X_{t_1}, X_{t_2})
		= {\rm Var}(X_0) h_{\lambda,c}(t_1) \left(h_{\lambda,c}(t_2)\right)+ \frac{1}{\lambda^2}\E\left[\int_0^{t_1} f_{\lambda}(t_2-s)f_{\lambda}(t_1-s)\varsigma^2(s)\sigma^2(X_{s})ds\right]\\
		&\hspace{.02cm}=  {\rm Var}(X_0) h_{\lambda,c}(t_1) \left(h_{\lambda,c}(t_2)\right)+ \frac{1}{\lambda^2}\int_0^{t_1} f_{\lambda}(t_2-t_1+u)f_{\lambda}(u)\varsigma^2(t_1-u) \sigma^2(\E\big[ X_{t_1-u} \big])du 
	\end{align*}
	where in the last equality, we have used the particular affine form of $\sigma^2(x)$.
	Consequently, the autocovariance function of the asymptotic (stationary) process is given by:
	\begin{align*}
		\mathrm{cov}(X_{t_2}^{\infty}, X_{t_1}^{\infty}) &= \lim_{t \to \infty}{\rm Cov}(X_{t+t_1}, X_{t+t_2}) ={\rm Var}(X_0) \lim_{t \to \infty} h_{\lambda,c}(t+t_1) \left(h_{\lambda,c}(t+t_2)\right) \\
		&\hspace{2cm}+ \frac{1}{\lambda^2}\lim_{t \to \infty}\int_0^{t+t_1} f_{\lambda}(t_2-t_1+u)f_{\lambda}(u)\varsigma^2(t+t_1-u)\sigma^2(\E\big[X_{t+t_1-u}\big])du
		\\ &= \ell^2_\infty {\rm Var}(X_0) + \frac{1}{\lambda^2}\varsigma^2_\infty\;\sigma^2\big(
		\ell_\infty \mathbb{E}[X_0]  + (1-a) \frac{\mu_\infty}{\lambda}\big)\int_0^{+\infty}  f_{\lambda}(t_2-t_1+u)f_{\lambda}(u)du.
	\end{align*}
	where the last equality follows from the fact that $ f_{\lambda}(t_2-t_1+\cdot)f_{\lambda}\!\in {\cal L}^1({\rm Leb}_1)$ since $f_{\lambda}\!\in {\cal L}^2({\rm Leb}_1)$ (Assumption~\ref{ass:int_holregul}~$(i)$), $\mbox{\bf 1}_{\{0\le u \le t+t_1\}}\varsigma^2(t+t_1-u)\to \lim_{t \to +\infty} \varsigma^2(t) =:\varsigma^2_\infty$ for every $u\!\in \R_+$ as $t\to +\infty$ (owing to \cite[Lemma 3.9]{EGnabeyeu2025}) and $\lim_{t\to+\infty}h_{\lambda,c}(t)= \ell_\infty$. 
	
	\smallskip
	\noindent  {\sc Step~3.} 
		To strengthen the functional convergence to convergence jointly with the
		initial state, take $n\in\N$ and fix $0\leq t_{1}<\dots<t_{n}$. Applying Proposition
	\ref{prop: convergence FT} for the particular choice $m_{t_{1},\dots,t_{n}}(ds)=\sum_{i=1}^{n}u_{i}\delta_{t_{n}-t_{i}}(ds)$,
	where $u_{1},\dots,u_{n}\in\mathbb{C}_{-}$,	we find owing to \( m_{t_{1},\dots,t_{n}}(\R_+)=\sum_{i=1}^{n}  u_{i}\), that as \(u\to\infty\),
	\[
	\mathbb E\!\left[
	e^{\sum_{j=1}^n u_jX_{u+t_j}}
	\mid X_0
	\right]
	\rightarrow \exp\Big[ \Xi_0\sum_{i=1}^{n}  u_{i} +  \frac{\varsigma^2_\infty}{2}\sigma^2\big(\Xi_0\big)\int_0^{\infty} \psi^2(s,m_{t_{1},\dots,t_{n}}) ds \Big]:=
	\Phi_{X_0}^{t_1,\ldots,t_n}(u_1,\ldots,u_n),
	\; \P\text{-a.s.}
	\]
	Therefore, for every bounded continuous function $G$, dominated convergence gives
	\begin{equation}\label{eq:CgceMix}
		\mathbb E\!\left[
		G(X_0)
		e^{\sum_{j=1}^n u_jX_{u+t_j}}
		\right]
		=
		\mathbb E\!\left[
		G(X_0)
		\mathbb E\!\left[
		e^{\sum_{j=1}^n u_j X_{u+t_j}}
		\mid X_0
		\right]
		\right]
		\longrightarrow
		\mathbb E\!\left[
		G(X_0)
		\Phi_{X_0}^{t_1,\ldots,t_n}(u_1,\ldots,u_n)
		\right].
	\end{equation}
	Since the law of $X_0$ is fixed and the family
	$(X_{u+\cdot})_{u\geq 0}$ is tight in
	$\mathcal C(\mathbb R_+,\mathbb R_+)$, the family
	$(X_0,X_{u+\cdot})_{u\geq 0}$ is jointly tight in
	$\mathbb R_+\times\mathcal C(\mathbb R_+,\mathbb R_+)$.
	Moreover, \eqref{eq:CgceMix} uniquely identifies
	the finite-dimensional distributions of every joint
	accumulation point. Consequently, on a suitable
	extension of the probability space, there exists a
	continuous process $X^\infty$ such that~\eqref{eq:joint_convergence} holds.
	In particular, this yields stable convergence with
	respect to $\sigma(X_0)$.
	To establish conditional stationarity, we combine the joint convergence~\eqref{eq:joint_convergence}, the identification of the limiting conditional Fourier-Laplace transforms in~\eqref{eq:CgceMix} by Proposition~\ref{prop: convergence FT}. Indeed, for every fixed \(h\ge0\), by the joint weak convergence~\eqref{eq:joint_convergence} and an arbitrary subsequence $\{h_k\}_k$,
		\begin{align}
		&\	\E\Big[ G(X_0)\exp\Big( \sum_{i=1}^{n} u_{i}  X_{t_{i}+h}^{\infty} \Big) \Big]
		= \lim_{k \to \infty} \E\Big[ G(X_0) \exp\Big(  \sum_{i=1}^{n} u_{i} X_{t_{i}+h}^{h_{k}}  \Big) \Big] = \lim_{k \to \infty} \E\Big[ G(X_0)\exp\Big( \sum_{i=1}^{n} u_{i} X_{h_{k}+h+t_{i}} \Big) \Big] \nonumber \\
		&\hspace{1cm}= \lim_{k \to \infty} \E\Big[ G(X_0) \exp\Big(  \int_{[0,\, h_{k}+h+t_{n}]} X_{h_{k}+h+t_{n}-s} \, m_{t_{1},\dots,t_{n}}(ds) \Big) \Big] =\E\left[
		G(X_0)\,
		\Phi_{X_0}^{t_1,\ldots,t_n}
		(u_1,\ldots,u_n)
		\right].\nonumber 
	\end{align}
where we applied Proposition~\ref{prop: convergence FT} as in~\eqref{eq:CgceMix}.
	As the preceding identity holds for every bounded
	continuous $G$, uniqueness of conditional
	expectations yields \(\mathbb P\)-almost surely,
	\begin{equation}\label{eq:weak_conv_finite_II}
	\E\left[
	\exp\left(
	\sum_{i=1}^n u_iX^\infty_{t_i+h}
	\right)\Bigm|X_0
	\right]
	=
	\Phi_{X_0}^{t_1,\ldots,t_n}
	(u_1,\ldots,u_n).
	\end{equation}
	\noindent Since \eqref{eq:weak_conv_finite_II} is independent of $h$ ($m_{t_1+h,\ldots,t_n+h}
	=
	m_{t_1,\ldots,t_n}$), the uniqueness of Fourier-Laplace transforms 
	therefore implies $\mathcal L\left(
	(X^\infty_{t_1+h},\ldots,X^\infty_{t_n+h})
	\mid X_0
	\right)
	=
	\mathcal L\left(
	(X^\infty_{t_1},\ldots,X^\infty_{t_n})
	\mid X_0
	\right).$
	Consequently, $X^\infty$ is strictly stationary,
	even conditionally on $X_0$.
	
	Finally, taking $n=1$ and invoking Theorem
	\ref{Theorem: limiting distribution}, we obtain
	$\mathcal L(X_t^\infty\mid X_0)=\pi_{X_0},
	\; t\ge0,$
	and consequently $	\mathbb E[X_t^\infty\mid X_0]
	=
	\int_{\mathbb R_+}x\,\pi_{X_0}(dx)
	=
	\Xi_0,
	\; \mathbb P\text{-a.s.}$
	The result follows by noticing that in the stabilized regime, \( \varsigma^2_\infty:= \lim_{t \to +\infty} \varsigma^2(t)= \lim_{t \to +\infty} \varsigma_{\lambda,c}^2(t) =\frac{c\lambda^2 (1-\ell_\infty^2)}{\|f_{\lambda}\|^2_{L^2(\text{Leb}_1 )}}\) owing to \cite[Lemma 3.9 ]{EGnabeyeu2025}.
	This completes the proof and we are done. \qed
	
	\bigskip
	\noindent {\bf Acknowledgement:}  The authors thank J-F. Chassagneux for insightful discussions, help and comments.
	
	\vspace{-.5cm}
	\bibliographystyle{alpha}
	\bibliography{Bibliography}
		
		\vspace{-.5cm}
	\appendix
	\section{Supplementary material and Proofs.}\label{app:lemmata}
	\noindent {\bf Proof of Proposition \ref{lm:fwdprocess}:}
	The first claim is a straightforward consequence of \cite[Equivalence Wiener-Hopf transform, Proposition 2.8]{EGnabeyeuR2025} (Equation~\eqref{eq:L^p-supBound}), under assumption \eqref{eq:hypoRlambda} \(({\cal K})\).
	
	\noindent For every \( s \in [0, T] \), define the process \( M^s = (M^s_t)_{t \in [0, T]} \) by
	\(
	M^s_t = \int_0^{t \wedge s} K(s- r) \sigma(r, X_r) \, \mathrm{d}W_r, \quad t \in [0, T]
	\). By the linear growth in Assumption~\ref{assump:kernelVolterra}(ii) and \cite[equation 4.58, Theorem 4.10]{EGnabeyeuR2025}, we have
	\[
	\mathbb{E}[\langle M^s \rangle_T] = 
	\int_0^{T \wedge s} |K(s-r)|^2 \mathbb{E} \left[| \sigma(r, X_r)|^2\right] \, \mathrm{d}r
	\leq C(1 + \| \varphi\|_{\infty,T}^2\mathbb{E}[|X_0|^2])
	\left( \int_0^s K(s-r)^{2} \, \mathrm{d}r \right) < \infty.
	\]
	\noindent Hence, \( \forall \, s \in [0, T] \), the process \( M^s \) is a martingale. Therefore, in the case \( s \geq t \), it holds by ~\eqref{eq:Volterrameanrevert} that
	\[
	\begin{aligned}
		\mathbb{E}[X_s \mid \mathcal{F}_t]
		&= x_0(s) + \mathbb{E} \left[ \int_0^s K(s-r) \left(\mu(r)-\lambda X_r\right) \, \mathrm{d}r \, \Big| \, \mathcal{F}_t \right]
		+ \mathbb{E} \left[ \int_0^s K(s-r) \sigma(r, X_r) \, \mathrm{d}W_r \, \Big| \, \mathcal{F}_t \right] \\
		&= X_0\varphi(s) +
		\int_0^s K(s-r) \left(\mu(r)-\lambda\mathbb{E}[X_r \mid \mathcal{F}_t]\right) \, \mathrm{d}r +
		\int_0^{t} K(s-r) \sigma(r, X_r) \, \mathrm{d}W_r.
	\end{aligned}
	\]
	The third equation in~\eqref{eq:Forward2} follows as a direct consequence of the Wiener-Hopf equation ~\eqref{eq:Volterrameanrevert2}. Alternatively, using the representation~\eqref{eq:Volterrameanrevert2}, one may derive this equation by considering the process \(M_t := \int_0^t f_\lambda(T - s)\, \sigma(s, X_s)\, \mathrm{d}W_s, \quad t \in [0, T],\)
	which is a local martingale. Its quadratic variation satisfies the estimate
	
	\centerline{$\mathbb{E}\left[\langle M \rangle_T\right] 
		\leq \int_0^T |f_\lambda(T - s)|^2\, \mathbb{E}[|\sigma(s, X_s)|^2]\, \mathrm{d}s 
		\leq C\,\|f_\lambda\|_{L^2(0,T)}^2\,\left(1 + \|\varphi\|_{\infty,T}^2\, \mathbb{E}[|X_0|^2]\right).$}
	\noindent which is finite under Assumption~\ref{ass:int_holregul}~$(i)$ in~\ref{ass:resolvent} together with the moment bound in \cite[Theorem 4.2]{EGnabeyeuR2025}. Therefore, \(M\) is a true martingale, and taking \(\mathcal{F}_t\)-conditional expectations concludes the proof and we are done.  \qed
	
	\medskip
	\noindent {\bf Proof of the Martingality Lemma~\ref{lem:martingality}.}
	Unlike the real-valued case, $M$ is $\mathbb{C}$-valued, so Fatou's lemma cannot be applied to $M$ directly. We first isolate a genuine nonnegative real local martingale by splitting $\psi$ into its real and imaginary parts, and treat the two separately.
	Write, for $t\in[0,T]$, $a(t):=\operatorname{Re}\psi(T-t)\le 0$ and $b(t):=\operatorname{Im}\psi(T-t)$; since $\psi\in L^2([0,T];\mathbb{C}_-)$ we have $a,b\in L^2(0,T;\mathbb{R})$. Write $\sigma_s:=\sigma(s,X_s)$. Since \(M_0\) is nonzero and \(\mathcal F_0\)-measurable, we set
	\begin{equation}\label{eq:R-def}
		R_t:=\exp\Big(\int_0^t a(s)\sigma_s\,dW_s-\frac12\int_0^t a(s)^2\sigma_s^2\,ds\Big),
	\end{equation}
	a nonnegative real local martingale with $R_0=1$, and $N_t:= M_t/M_0 R_t$. Since $\psi(T-\cdot)-a=ib$ and $\psi(T-\cdot)^2-a^2=2iab-b^2$, substituting into the exponent of $ M$ gives
	\begin{equation}\label{eq:N-formula}
		N_t=\exp\Big(i\int_0^t b(s)\sigma_s\,dW_s - i\int_0^t a(s)b(s)\sigma_s^2\,ds+\frac12\int_0^t b(s)^2\sigma_s^2\,ds\Big),\qquad M=M_0 RN.
	\end{equation}
	
	\smallskip
	\noindent {\sc Step~1} \textit{($R$ is a true martingale on $[0,T]$)}. Being a nonnegative local martingale, $R$ is a supermartingale by Fatou's lemma; it thus suffices to show $\mathbb{E}^{\mathbb P}[R_T]\ge 1$. Let $\tau_n:=\inf\{t\ge 0: X_t\ge n\}\wedge T$. On $[0,\tau_n]$, $\sigma_s^2\le \varsigma^2(s)(\kappa_0+\kappa_1 n)\le \Sigma_n:=(\kappa_0+\kappa_1 n) \sup_{s\le T}\varsigma^2(s)<\infty$, so
	$\int_0^{\tau_n} a(s)^2\sigma_s^2\,ds \ \le\ \Sigma_n\,\|a\|_{L^2(0,T)}^2\ <\ \infty,$ deterministically
	and Novikov's condition holds for $a\mathbf 1_{[0,\tau_n]}\sigma$ on $[0,T]$. Hence $R^{\tau_n}$ is a true, uniformly integrable martingale for each \( n \), and we may define a measure change $d\mathbb Q^n:=R_{\tau_n}\,d\mathbb P$. By Girsanov's theorem, $W_t^n:=W_t-\int_0^{t\wedge\tau_n}a(s)\sigma_s\,ds$ is a $\mathbb Q^n$-Brownian motion on $[0,T]$, and
	\begin{equation}\label{eq:X-under-Qn}
		X_t = x_0(t)+\int_0^t K(t-s)\big(\mu(s)-\lambda X_s+\mathbf 1_{[0,\tau_n]}(s)\,a(s)\,\sigma^2(s,X_s)\big)ds+\int_0^t K(t-s)\sigma(s,X_s)\,dW_s^n.
	\end{equation}
	 Markov's inequality then gives 
	\begin{align}
		&\, \mathbb{Q}^n(\tau_n < T) = \mathbb{Q}^n\Big( \sup_{t \in [0,T]} X_t \ge n \Big) \leq  \frac{1}{n} \mathbb{E}_{}^{\mathbb{Q}^n} \Big( \sup_{t \in [0,T]} |X_t| \Big) \leq \frac{C}{n} 
		\to 0 \quad \text{as } n \to \infty,\label{eq:BoundLim}
	\end{align}
	provided that  $\sup_n\mathbb{E}^{\mathbb Q^n}\big[\sup_{t\le T}X_t\big]\le C$ for a constant $C=C(K,T,\kappa_0,\kappa_1,\varsigma,x_0,\mu,a)$ independent of $n$ (see Step~2 further on). Whence
	$$\mathbb{E}^{\mathbb P}[R_T]\ \ge\ \mathbb{E}^{\mathbb P}\big[R_T\mathbf 1_{\{\tau_n=T\}}\big]=\mathbb{E}^{\mathbb P}\big[R_{\tau_n}\mathbf 1_{\{\tau_n=T\}}\big]=\mathbb{Q}^n(\tau_n=T)\ \longrightarrow\ 1,\quad \text{as } n \to \infty.$$
	Thus $\mathbb{E}^{\mathbb P}[R_T]\ge 1$; combined with the supermartingale property this forces $\mathbb{E}^{\mathbb P}[R_t]=1$ for all $t\in[0,T]$, hence $R$ is a true martingale on $[0,T]$.

	\smallskip
	\noindent {\sc Step~2} \textit{(Proof of the desired boundedness result)}.
	To obtain~\eqref{eq:BoundLim}, it thus suffices to show that $\sup_n\mathbb{E}^{\mathbb Q^n}\big[\sup_{t\le T}X_t\big]\le \infty$. For this, we use the representation~\eqref{eq:X-under-Qn} under the law $\mathbb{Q}^n$, which can equivalently be written (see e.g., Proposition~\ref{lm:fwdprocess}) with \(R_\lambda:=\varphi- \int_0^\cdot f_{\lambda}(\cdot-s)\varphi(s)ds\).
	\begin{equation}
		X_t= X_0R_\lambda(t) +\frac{1}{\lambda}\int_0^t f_{\lambda}(t-s)\big(\mu(s)+\mathbf 1_{[0,\tau_n]}(s)\,a(s)\,\sigma^2(s,X_s)\big)ds + \frac{1}{\lambda}\int_0^t f_{\lambda}(t-s)\sigma(s, X_s) \,dW_s^n.
	\end{equation} 
	Since $X\ge0$, bounding $\mathbf 1_{[0,\tau_n]}(s)|a(s)|\le|a(s)|$ and for $\theta>0$, multiplying by $e^{-\theta t}$ and writing $X_s=e^{\theta s}\cdot e^{-\theta s}X_s$:
	$$e^{-\theta t}\!\int_0^t\!f_{\lambda}(t-s)\varsigma^2(s)|a(s)|X_s\,ds \ \le\ \Big(\sup_{r\le T}e^{-\theta r}X_r\Big)\int_0^tf_{\lambda}(t-s)e^{-\theta(t-s)}\varsigma^2(s)|a(s)|\,ds.$$
	By Cauchy--Schwarz, the second factor is bounded, uniformly in $t\in[0,T]$ and in $n$, by
	$$\epsilon(\theta):=\|\varsigma^2\|_{\infty,T}\,\|a\|_{L^2(0,T)}\Big(\int_0^Te^{-2\theta u}f_{\lambda}(u)^2\,du\Big)^{1/2},$$
	which tends to $0$ as $\theta\to\infty$ by dominated convergence owing to~\eqref{eq:contf_lambda}($f_{\lambda}\in L^2(0,T)$, $e^{-2\theta u}\to0$ a.e.). Hence, 
	setting
	$$D(\theta):=\sup_{t\le T}e^{-\theta t}\Big(R_\lambda(t)+\frac1\lambda\int_0^tf_\lambda(t-s)\big(|\mu(s)|+\kappa_0\varsigma^2(s)|a(s)|\big)ds\Big)<\infty,$$
	(finite and independent of $n$, by continuity of $R_\lambda$ and of the 
	 convolution on compact $[0,T]$), we get
	$$\sup_{t\le T}e^{-\theta t}X_t \ \le\ D(\theta)\big(1+X_0\big)+\frac{\kappa_1}{\lambda}\epsilon(\theta)\sup_{t\le T}e^{-\theta t}X_t+ \frac1\lambda\sup_{t\le T}\Big|e^{-\theta t}\!\int_0^tf_{\lambda}(t-s)\sigma(s,X_s)\,dW_s^n\Big|.$$
	Fix $\theta^*$ such that $\frac{\kappa_1}{\lambda}\epsilon(\theta^*)\le\frac12$, this choice depends only on $\varphi, f_{\lambda},\varsigma,a$, hence is ``uniform in $n$". Absorbing the term gives, since $t\mapsto e^{-\theta^*t}X_t$ is a.s.\ continuous on $[0,T]$ and hence finite:
	\begin{equation}\label{eq:absorbed}
		\sup_{t\le T}e^{-\theta^*t}X_t \ \le\ 2D(\theta^*) \big(1+X_0\big) + \frac2\lambda\sup_{t\le T}\Big|e^{-\theta^*t}\!\int_0^tf_\lambda(t-s)\sigma(s,X_s)\,dW_s^n\Big| \quad\mathbb Q^n\text{-a.s.}
	\end{equation}
	It remains to bound the martingale term, uniformly in $n$. Since $\sigma^2(s,X_s)
	=
	\varsigma^2(s)(\kappa_0+\kappa_1X_s)$, 
	 	let $p$ be given by Assumption~\ref{ass:int_holregul} such that $p(\frac{\beta-1}{2\beta} \wedge \vartheta)>1$ and
	 	set $Y_{n,p}:=\E^{\mathbb Q^n}[\sup_{t\le T}e^{-p\theta^*t}X_t^p]$,
	 	$Z_t^n
	 	:=
	 	e^{-\theta^*t}
	 	\int_0^t f_\lambda(t-s)\sigma(s,X_s)\,dW_s^n.$
	 	Using the Burkholder--Davis--Gundy inequality together with
	 	Assumption~\ref{ass:int_holregul}, we obtain, 
	 	\[
	 	\mathbb E^{\mathbb Q^n}
	 	\left[
	 	|Z_t^n-Z_s^n|^p
	 	\right]
	 	\le
	 	C_{T,p} |t-s|^{p(\frac{\beta-1}{2\beta} \wedge \vartheta)} \big(1+Y_{n,p}^{1/2}\big),
	 	\qquad 0\le s\le t\le T,
	 	\]
	 	where \(C_{T,p}\) is independent of \(n\), but depends on the deterministic regularity constants of \(f_\lambda\). Therefore, by the
	 	Garsia--Rodemich--Rumsey inequality, for every
	 	\(a<\frac{\beta-1}{2\beta} \wedge \vartheta-1/p\),
	 	\[
	 	\mathbb E^{\mathbb Q^n}
	 	\left[
	 	\|Z^n\|_{C^a([0,T])}^p
	 	\right]:=
	 	\mathbb E^{\mathbb Q^n}
	 	\left[\sup_{0 \leq s< t \leq T} \frac{|Z^n_t - Z^n_s|^p}{|t - s|^{pa}}\right]
	 	\le
	 	C_{T,p,a}\big(1+Y_{n,p}^{1/2}\big).
	 	\]
	 	Since \(Z_0^n=0\), raising \eqref{eq:absorbed} to the power of \(p\) and taking
	 	\(\mathbb Q^n\)-expectations therefore gives
	 	\[
	 	Y_{n,p}
	 	\le
	 	\big(2D(\theta^*)\big)^p \big(1+\E^{\mathbb Q^n}\big[|X_0|^p\big]\big)
	 	+
	 	\frac{2^p}{\lambda^p}
	 	\E^{\mathbb Q^n}\Big[\sup_{t\le T}|Z_t^n|^p\Big]
	 		\le
	 		A+B\,Y_{n,p}^{1/2},
	 	\] 
	 	with constants \(A,B\) independent of \(n\) since $\mathbb E^{\mathbb Q^n}[|X_0|^p]
	 	=
	 	\mathbb E^{\mathbb P}[|X_0|^p]
	 	<\infty$. Hence $\sqrt{Y_{n,p}}
	 	\le
	 	\frac{B+\sqrt{B^2+4A}}{2}<\infty,$
	 	and consequently $\sup_{n\ge1}Y_{n,p}<\infty$ and $\sup_{n\ge1}Y_{n,1}<\infty.$
	Finally, multiplying by $e^{\theta^*T}$ concludes
	\vspace{-.1cm}
	\begin{equation}
	\sup_{n\ge1}
	\mathbb E^{\mathbb Q^n}
	\Big[\sup_{t\le T}X_t\Big]
	\le
	e^{\theta^*T}\sup_{n\ge1}Y_{n,1}
	<\infty,
	\vspace{-.1cm}
	\end{equation}
	which proves~\eqref{eq:BoundLim}, a bound independent of $n$.
	
	\smallskip
	\noindent {\sc Step~3} \textit{($M_0 N$ is a true $\mathbb Q$-martingale)}. Since $R$ is now a global true martingale, $d\mathbb Q:=R_T\,d\mathbb P$ defines a probability measure on $\mathcal F_T$, and $W_t^{\mathbb Q}:=W_t-\int_0^t a(s)\sigma_s\,ds$ is a $\mathbb Q$-Brownian motion on $[0,T]$. Note that under $\mathbb{Q}$ the process $X$ solves \eqref{eq:X-under-Qn} without truncation, with mean-reversion rate $\lambda(s):=\lambda-\kappa_1\varsigma^2(s)a(s)\ge\lambda$ (as $a\le0$) for $s\in[0,T]$. Substituting $dW_s=dW_s^{\mathbb Q}+a(s)\sigma_s\,ds$ into \eqref{eq:N-formula}, the term $i\int_0^t ab\sigma^2\,ds$ cancels exactly, leaving
	\vspace{-.1cm}
	\begin{equation}\label{eq:N-under-Q}
		N_t=\exp\Big(i\int_0^t b(s)\sigma_s\,dW_s^{\mathbb Q}+\frac12\int_0^t b(s)^2\sigma_s^2\,ds\Big)\qquad\text{under }\mathbb Q,
		\vspace{-.1cm}
	\end{equation}
	manifestly a $\mathbb Q$-local martingale by Itô's formula. Since the stochastic integral in \eqref{eq:N-under-Q} is real, $\exp\big(i\int_0^t b\sigma\,dW^{\mathbb Q}\big)$ has modulus exactly $1$, so
	\vspace{-.1cm}
	\begin{equation}\label{eq:modulus-exact}
		|N_t| = \exp\Big(\frac12\int_0^t b(s)^2\sigma^2(s,X_s)\,ds\Big) \qquad\text{pathwise, for every }t\in[0,T],
		\vspace{-.1cm}
	\end{equation}
	an a.s.\ nondecreasing process; in particular $\sup_{t\le T}|N_t|=|N_T|$ a.s. 
	Since \(M_0\) is \(\mathcal F_0\)-measurable and \(N\) is a \(\mathbb Q\)-local martingale, \(M_0N\) is also a \(\mathbb Q\)-local martingale.
	Since \(X\geq0\) and \(\Re m\leq0\), we have $|M_T|
	=
	\exp\Big(
	\int_0^T X_{T-s}\,\Re m(ds)
	\Big)
	\leq1.$
	On the other hand, $|M_T|
	=
	R_T\,|M_0 N_T|.$
	Therefore, $\mathbb E^{\mathbb Q}[|M_0N_T|]
	=\mathbb E[R_T|M_0N_T|]
	=\mathbb E[|M_T|]
	\le1.$
	Consequently, $\mathbb E^{\mathbb Q}
	\Big[\sup_{t\le T}|M_0N_t|\Big]
	=\mathbb E^{\mathbb Q}[|M_0N_T|]<\infty$ by \eqref{eq:modulus-exact}, so $M_0 N$ is a $\mathbb Q$-local martingale dominated on $[0,T]$ by an integrable random variable; dominated convergence along any localising sequence shows $M_0N$ is a true, uniformly integrable $\mathbb Q$-martingale on $[0,T]$.
	
	\smallskip
	\noindent {\sc Step~4} \textit{(Conclusion)}. By the generalised Bayes formula, for $0\le s\le t\le T$,
	$$\mathbb{E}^{\mathbb P}[M_t\mid\mathcal F_s] = \mathbb{E}^{\mathbb P}[R_t M_0 N_t\mid\mathcal F_s] = R_s\,\mathbb{E}^{\mathbb Q}[ M_0 N_t\mid\mathcal F_s] = R_s M_0 N_s = M_s,$$
	using Step~3 for the second equality. Hence $M$ is a true martingale on $[0,T]$. \qed
	
	\medskip
	\noindent {\bf Proof of Proposition \ref{prop:existence_VolterraHestonS_V}.}
	The first claim~\ref{T:VolHeston:1} of the proposition is a direct consequence of \cite[Theorem 3.2]{EGnabeyeuR2025} combined with the fact that $S$ is fully determined by $V$. The uniqueness in law comes from the existence of a solution to the Riccati Volterra equation together with Theorem~\ref{T:VolSqrt}. 
	Let $u:=1-\frac{\mu}{\lambda m_0}$ such that  $u\in L^2_{\rm loc}([0,\infty))$
	by assumption. Then $\varphi$ in~\eqref{eq:CondMean} reads $\varphi=1+\lambda K*u$ and satisfies $\varphi(0)=1$ and, for every $T>0$
	and $0\le s<t\le T$, we have
	\vspace{-.1cm}
	\begin{align*}
		|\varphi(t)-\varphi(s)|
		&\;\le\;\lambda\,\|u\|_{L^2(0,T)}\,\Big( \int_0^s\!\big(K(t-r)-K(s-r)\big)^2dr+\int_s^t K(t-r)^2dr\Big)^\frac12\\
		&\;\le\; C_{u,\lambda,T, K}\,|t-s|^{\widehat\theta\wedge\theta}
	\end{align*}
	where we used Cauchy--Schwarz, together with Assumption~\ref{assump:kernelVolterra}~$(i)$. Consequently, the H\"older continuity with exponent
	$\widehat\theta\wedge\theta$, as stated in the first assertion
	of Proposition~\ref{prop:existence_VolterraHestonS_V}-\eqref{T:VolHeston:1}, follows by an argument
	analogous to that in~\cite[Proof of Theorem 3.2]{EGnabeyeuR2025}.
	The claims in part \ref{T:VolHeston:2} regarding existence, uniqueness of \((\psi_1,\psi_2)\) in~\eqref{eq:riccati_heston_gen}, and non-positivity of ${\rm Re\,}\psi_2$ follow from \cite[Proposition 4.2]{ackermann2022}.
	Indeed, fix \(T>0\).
	The first component of~\eqref{eq:riccati_heston_gen} is explicit. Once \(\psi_1\) is given, the second component is the
	scalar Riccati--Volterra equation
	\begin{equation}\label{eq:RcctiPsi2}
		\psi_2(t)
		=
		u_2K(t)
		+
		\int_0^t K(t-s)
		\Big(
		g_T(s)
		+q_T(s)\psi_2(s)
		+\frac{\nu^2}{2}\varsigma^2(T-s)\psi_2(s)^2
		\Big)ds,
	\end{equation}
	where $g_T(s)
	:=
	f_2(s)
	+\frac12\bigl(\psi_1(s)^2-\psi_1(s)\bigr),
	\;
	q_T(s)
	:=
	\rho\nu\varsigma(T-s)\psi_1(s)-\lambda.$
	Since \(0\leq\Re\psi_1\leq1\) and \(\Re f_2\leq0\), we can verify that $	\Re g_T(s)\leq0$ and also that the closed left half-plane \(\mathbb C_-:=\{z\in\mathbb C:\Re z\leq0\}\) is invariant for the scalar
	Riccati--Volterra equation~\eqref{eq:RcctiPsi2} for \(\psi_2\).
	
	Moreover, under Assumption~\ref{ass:int_holregul} and the conditions $0\leq\Re\psi_1\leq1,\;
	\Re u_2\leq0,\;
	\Re f_2\leq0,$
	the proof of Theorem~\ref{theorem:riccati-extension} carries over to the
	present equation and yields a unique solution
	$\psi_2\in L^2([0,T];\mathbb C),\;
	\Re\psi_2\leq0$ 	(alternatively, one may invoke
	\cite[Proposition~4.2]{ackermann2022}). 
	Finally, with $	m(ds)=u\,\delta_0(ds)+f(s)\,ds,$
	the pair \(\psi=(\psi_1,\psi_2)\) is precisely the solution of the
	measure-extended Riccati--Volterra equation associated with the affine
	Volterra system~\eqref{eq:Systemlog_S_V}. Therefore,Theorem~\ref{T:VolSqrt}~$(2)$ applies directly and
	yields the exponential-affine representation
	\eqref{eq: MeasureFLplaceZero}--\eqref{eq: MeasureFLplaceZero II}(see also, e.g., \cite[Proposition 4.3]{ackermann2022}).
	
	Finally, note that the	representation~\eqref{eq:formula_S_inhom_Heston} follows directly from
	\eqref{eq:log_S_t} so that $S$ given by~\eqref{eq:formula_S_inhom_Heston} clearly solves the first equation in~\eqref{E:HestonlogS_V}.
	$B:=\sqrt{1-\rho^2}\,W^{(1)}+\rho W^{(2)}$ is a Brownian motion by Lévy's characterization 
	and $S_t=S_0\mathcal E\left(\int_0^\cdot\sqrt{V_s}\,dB_s\right)_t,$
	so that \(S\) is a nonnegative local martingale. Applying
	Part~\ref{T:VolHeston:2} with \(u=(1,0)\) and \(f\equiv0\) gives
	\(\psi_1\equiv1\), while the second Riccati--Volterra equation reduces to~\eqref{eq:Char_p} with \(p=1\), \(K_\alpha\equiv K\) and clearly admits
	\(\psi_2\equiv0\) as its unique solution. Hence the affine transform
	formula yields $\mathbb E[S_T]=S_0,\; T\geq0.$
	Since a nonnegative local martingale is a supermartingale, the constancy
	of its expectation implies that \(S\) is a true martingale. 
	\qed

	\medskip
	\noindent {\bf Proof of Corollary \ref{corol:StabilizedRHeston}.}
	Indeed, the fake stationary rough Heston model arises as a specific instance of the inhomogeneous Volterra-Heston framework analyzed in \cite[Corollary 4.5.]{ackermann2022} when the volatility coefficient is constant, i.e., $\eta(t) \equiv 1$, the parameters $\mu$, and $\sigma(\cdot,\cdot)$ are time-dependent, and the kernel is chosen to be fractional.
	Consider \(\mathbb{M}_{\mathrm{lf}}(\R_+, (\mathbb C^2)^*) \ni m(\dd s) := u\,\delta_0(\dd s) +f(s)\,\mathrm{Leb}_{\mathbb{R}_+}(\mathrm{d}s)\), the existence follows from the second claim of Proposition~\ref{prop:existence_VolterraHestonS_V}.
	Consequently, formula in equation~\eqref{eq: MeasureFLplaceZero II} of Theorem~\ref{T:VolSqrt} (b) applies for $(\log S, V)$ and yields the expression in   \eqref{eq:time_zero_frac_heston_all} since on the first hand:
	\vspace{-0.3cm}
	\begin{align*}
		\int_0^T F(s, \psi(T-s))\, \begin{bmatrix}
			\log(S_0) \\  V_0 \varphi(s)
		\end{bmatrix}\, \dd s &= \int_0^T \begin{bmatrix}
			\log(S_0) & V_0 \varphi(T-s)
		\end{bmatrix} \left(\psi(s) \begin{bmatrix}
			0 & -\frac12 \\ 0 & -\lambda
		\end{bmatrix} + \frac12 \begin{bmatrix}
			0 \\ \psi(s) A A^{\top} \psi^{\top}(s)
		\end{bmatrix} \right)\, \dd s \\
		&= \chi(T) V_0
	\end{align*}
	where \( A = A(T-s) = \begin{pmatrix}
		\sqrt{1 - \rho^2} & \rho \\
		0 & \nu \varsigma(T-s)
	\end{pmatrix}\), so that on the second hand, the full Fourier--Laplace representation for the integrated log-price and variance:
		\vspace{-.3cm}
		\begin{align*}
			\begin{split}
				& \E_{(\log(S_0),V_0)}\left[ \exp\left( u_1 \log(S_T) + u_2 V_T + \int_0^T f_1(T-u) \log(S_u) \,du   + \int_0^T f_2(T-u) V_u \,du \right) \right] \\
				& = 
				\exp\left( \phi(T) + \left(u_1 + \int_0^T f_1(s)\,\dd s \right)\log(S_0) + \left(u_2\,\varphi(T) + \int_0^T f_2(T-s)\, \varphi(s) \,\dd s + \chi(T)\right) V_0 \right).
			\end{split}\\
			&\ \forall\, t \in [0,T],\; \chi(t) = \int_0^t \varphi(t-s)\;\Big(  \frac{1}{2} (\psi^2_1(s)-\psi_1(s)) + \left( \rho \nu \varsigma(t-s) \psi_1(s) - \lambda \right)  \psi_2(s) + \frac{\nu^2}{2} \varsigma^2(t-s) \psi_2^2(s) \Big) \dd s,\\
			&\quad\quad= \int_0^t \Big(  \frac{1}{2} (\psi^2_1(s)-\psi_1(s)) + \left( \rho \nu \varsigma(t-s) \psi_1(s) - \lambda \right)  \psi_2(s) + \frac{\nu^2}{2} \varsigma^2(t-s) \psi_2^2(s) \Big) \dd s -\lambda \Bigg[
			\int_0^t \Big( \frac{\mu(u)}{\lambda m_0} - 1 \Big) \\
			& \quad \times \Bigg(
			\int_0^{t-u} K_{\alpha}(t-u-s) 
			\Big( \frac{1}{2} (\psi_1^2(s)-\psi_1(s)) 
			+ (\rho \nu \varsigma(t-s) \psi_1(s) - \lambda) \psi_2(s) 
			+ \frac{\nu^2}{2} \varsigma^2(t-s) \psi_2^2(s) \Big) \dd s
			\Bigg)
			du
			\Bigg].
		\end{align*}
		where we used the expression of \(\varphi\) in~\eqref{eq:CondMean} 
		together with Fubini-Tonelli's theorem and a change of variables, so that we have
		\begin{align*}
			&\,u_2\,\varphi(t) + \int_0^t f_2(t-s)\, \varphi(s) \,\dd s + \chi(t) = (I^{1-\alpha} \psi_2)(t)
			\\
			&\qquad\qquad\qquad -\lambda \Bigg[
			\int_0^t \Big( \frac{\mu(u)}{\lambda m_0} - 1 \Big) \Bigg( u_2\,K_{\alpha}(t-u) + \int_0^{t-u} K_{\alpha}(t-u-s)f_2(s)\, \varphi(s) \,\dd s \\
			& \quad +
			\int_0^{t-u} K_{\alpha}(t-u-s) 
			\Big( \frac{1}{2} (\psi_1^2(s)-\psi_1(s)) 
			+ (\rho \nu \varsigma(t-s) \psi_1(s) - \lambda) \psi_2(s) 
			+ \frac{\nu^2}{2} \varsigma^2(t-s) \psi_2^2(s) \Big) \dd s
			\Bigg)
			du
			\Bigg] \\
			&\qquad\qquad\qquad= (I^{1-\alpha} \psi_2)(t) -\lambda \int_0^t \Big( \frac{\mu(u)}{\lambda m_0} - 1 \Big) \psi_2(t-u) \,du.
		\end{align*}
		still owing to Fubini-Tonelli's theorem.
		The expression in~\eqref{eq:time_zero_frac_heston} follows directly from standard results in fractional calculus, in conjunction with equation~\eqref{eq:riccati_heston_gen}.
		\qed
	\begin{Lemma}\label{lm:firstmoment} 
		Let \(g: \R_+\to \R\) be a continuous function such that \(\lim_{t \to +\infty}g(t) =g_\infty \in\R\). $\varsigma$ the borel bounded continuous function given in~\eqref{eq:sigma-definition} such that $\varsigma^2_\infty:= \lim_{t \to +\infty} \varsigma^2(t)$. Let Assumptions \ref{ass:resolvent}-\ref{ass:int_holregul} and~\ref{ass:Mean_revertF} be in force.
		Consider a measure \( m \in \mathbb{M}^{-}_{\mathrm{lf}}(\R_+, \mathbb{C}) \) with \( |m|(\mathbb{R}_+) < \infty \) and let \(  \psi_\infty(\cdot, m):=\psi_\infty \in L^1(\mathbb{R}_+; \mathbb{C}_{-}) \cap L^2(\mathbb{R}_+; \mathbb{C}_{-}) \), the solution to the  measure-extended Riccati--Volterra equation~\eqref{eq:measureFLplce2}-\eqref{eq:FVarsigma_infty}, $\psi_T := \psi_T(\cdot,m)$ that of the Riccati-Volterra equation~\eqref{eq:measureFLplce}-\eqref{RicVolSqrt2}
		on $[0,T]$ for every $T>0$.
		Then,
		\begin{enumerate}
		\item[(i)]	The expected process at future time defined in Equation~\eqref{eq:DiffForward} satisfied in the long run  
			\begin{align*}
				\lim_{t \to +\infty} \mathbb{E}_{}\big[X_t | X_0\big] =	\lim_{t \to +\infty}	\xi_0(t) = \ell_\infty X_0 + (1-a) \frac{\mu_\infty}{\lambda} := \Xi_0.
			\end{align*}
			\item[(ii)] $\lim_{t \to \infty}\int_0^t g(t-s) \,m(ds) = g_\infty\,m(\R_+)$ and $\lim_{t \to \infty} \int_0^t g(s)\, \psi_t(t - s)\, ds 
			= g_\infty\int_0^{\infty} \psi_\infty(s) ds$.
			\\
			\item[(iii)] $\lim_{t \to \infty} \int_0^t \varsigma^2(s)\, g(s)\, \psi_t(t - s)^2\, ds 
			= \varsigma^2_\infty g_\infty\int_0^{\infty} \psi_\infty(s)^2 ds$.
		\end{enumerate}
	\end{Lemma}
	\noindent {\bf Proof:}
	Note that \( f_\lambda \) is integrable, hence we can pass to the limit \( t \to \infty \) in~\eqref{eq:DiffForward}.
	{\small
		\begin{align*}
			\lim_{t \to +\infty} \mathbb{E}_{}\big[X_t | X_0\big] = \lim_{t \to +\infty}	\xi_0(t) &= X_0\lim_{t \to +\infty} \Big(\varphi(t) - \int_0^t f_{\lambda}(t-r) \varphi(r)\, dr\Big)  + \lim_{t \to +\infty}\frac{1}{\lambda} \int_0^t f_{\lambda}(t-r) \mu(r) \, \mathrm{d}r\\
			& = \ell_\infty X_0  + (1-a) \frac{\mu_\infty}{\lambda},
		\end{align*}
	}
	\noindent where we used Assumptions \ref{ass:resolvent} and~\ref{ass:Mean_revertF} together with the dominated convergence theorem.
	This proves the desired convergence. For the first claim in $(ii)$, we have
	\begin{align*} 
		\int_0^t g(t-s) \,m(ds) - g_\infty\,m(\R_+) = \int_0^t(g(t-s)-g_\infty) \,m(ds) - g_\infty \int_t^{+\infty} m(ds).
	\end{align*} 
	so that, by the triangle inequality,
	\(\left|	\int_0^t g(t-s) \,m(ds) - g_\infty m(\R_+)\right| \leq \left| \int_0^t(g(t-s)-g_\infty) \,m(ds)  \right| + g_\infty\, |m| ((t,\infty))\).
	First note that we can split the first integral as follows:
	\vspace{-0.3cm}
	\begin{align*}
		\left| \int_0^t(g(t-s)-g_\infty) \,m(ds)  \right|  &\leq \left|\int_0^{t-A_{\epsilon}} (g(t-s)-g_\infty) \,m(ds)\right|
		+ \left|\int_{t-A_{\epsilon}}^t(g(t-s)-g_\infty) \,m(ds)\right| \\
		&\leq  \epsilon|m|(\R_+)+ \Big(\sup_{s \geq 0}|g(s)|+|g_\infty|\Big) |m|((t-A_{\epsilon},t]) .
	\end{align*}
	where \( A_{\epsilon} \) is chosen such that \( \forall\,s \ge A_{\epsilon} \), we have \( \big|g(s) - g_\infty\big| \le \epsilon \). Moreover, \( \forall\,s \in (0, t - A_{\epsilon}) \), we have \( t - s \ge A_{\epsilon} \) for \( t \) large enough (say \( t \ge 2A_{\epsilon} \)), and hence this implies that $\left| g(t-s) - g_\infty \right| \le \epsilon, \; \forall s \in (0, t - A_{\epsilon}).$
	We thus have:
	\vspace{-0.4cm}
	\[
	\left|	\int_0^t g(t-s) \,m(ds) - g_\infty\,m(\R_+)\right| \le\epsilon|m|(\R_+)+ \Big(\sup_{s \geq 0}|g(s)|+|g_\infty|\Big) |m|((t-A_{\epsilon},t]) + 2 |g_\infty| \,|m| ((t-A_{\epsilon},\infty)).
	\]
	Since \( |m|(\mathbb{R}_+) < \infty \), we have \(|m|((t-A_{\epsilon},t]) \to 0\) and \( |m|((t-A_{\epsilon}, \infty)) \to 0 \) as \( t \to \infty \), which proves the first assertion in (ii).
	\noindent Likewise, the proofs of the second claim in $(ii)$ and of $(iii)$ rely on similar arguments; we therefore only prove $(iii)$. We have
	\vspace{-0.3cm}
	\begin{align*}
		 \Big|\int_0^t \varsigma^2(s)\, g(s)\, \psi_t(t - s)^2\, ds 
	&\,	- \varsigma^2_\infty g_\infty \int_0^{\infty} \psi_\infty(s)^2 ds\Big|
		\leq \int_0^t \left| \varsigma^2(t - s)\, g(t-s)\, 
		- \varsigma^2_\infty  g_\infty \right| |\psi_\infty(s)|^2 ds \\ & + \int_0^t  \varsigma^2(t - s)\, |g(t-s)|\, \left|\psi_t(s)^2- \psi_\infty(s)^2 \right| ds + \varsigma^2_\infty |g_\infty|\int_t^{\infty} |\psi_\infty(s)|^2 ds.
	\end{align*}
	The third term tends to zero since $\psi_\infty \in L^2(\R_+; \mathbb C)$, likewise for the second term since $g \varsigma $ is bounded by Cauchy--Schwarz,$\|\mathbf1_{[0,t]}\psi_t^2
	-\psi_\infty^2\|_{L^1}\leq
	\left(
	\|\mathbf1_{[0,t]}\psi_t\|_{L^2}
	+\|\psi_\infty\|_{L^2}
	\right)
	\|\mathbf1_{[0,t]}\psi_t-\psi_\infty\|_{L^2}
	\longrightarrow0,$
	owing to~\eqref{eq:Stabil}. For the first term, choose \( A_{\epsilon} \) such that for all \( s \ge A_{\epsilon} \), we have \( \left| \varsigma^2(s)\, g(s)\, 
	- \varsigma^2_\infty  g_\infty \right| \le \epsilon \). Moreover, for all \( s \in (0, t - A_{\epsilon}) \), we have \( t - s \ge A_{\epsilon} \) for \( t \) large enough (say \( t \ge 2A_{\epsilon} \)), and hence this implies that $\Big| \varsigma^2(t - s)\, g(t-s)\, 
	- \varsigma^2_\infty  g_\infty \Big| \le \epsilon, \; \forall s \in (0, t - A_{\epsilon}).$
	\vspace{-0.4cm}
	\begin{align*}
		&\ \int_0^t \left| \varsigma^2(t - s)\, g(t-s)\, 
		- \varsigma^2_\infty  g_\infty \right| |\psi_\infty(s)|^2 ds =\int_0^{t-A_{\epsilon}}\left| \varsigma^2(t - s)\, g(t-s)\, 
		- \varsigma^2_\infty  g_\infty \right| |\psi_\infty(s)|^2 ds\\
		&\hspace{5cm}+ \int_{t-A_{\epsilon}}^t\left| \varsigma^2(t - s)\, g(t-s)\, 
		- \varsigma^2_\infty  g_\infty \right|
		|\psi_\infty(s)|^2 ds \\
		&\hspace{1.6cm}\leq  \epsilon \int_0^{+\infty} |\psi_\infty(s)|^2 \, ds + 
		\int_0^{\infty} \mathbf{1}_{[t-A_{\epsilon},\, t]}(s)\, \left( \sup_{s \geq 0}|g(s)|\,\varsigma^2(t - s) + \varsigma^2_\infty |g_\infty| \right) \, |\psi_\infty(s)|^2 \, \mathrm{d}s\\
		&\hspace{1.6cm}\leq  \epsilon \int_0^{+\infty} |\psi_\infty(s)|^2 \, ds + \left( \|\varsigma^2\|_\infty \sup_{s \geq 0} |g(s)|+  \varsigma^2_\infty |g_\infty| \right)
		\int_0^{\infty} \mathbf{1}_{[t-A_{\epsilon},\, t]}(s)\, |\psi_\infty(s)|^2 \, \mathrm{d}s.
	\end{align*}
	By the dominated convergence theorem, the second term vanishes as \( t \to \infty \).
	Since \(\varepsilon>0\) was arbitrary, this establishes (iii) and completes the proof.
	\qed
	
	\medskip
	\noindent {\bf Proof of Proposition \ref{prop: convergence FT}.}
	For every $T>0$,  assuming that $\psi_T := \psi_T(\cdot,m)$ satisfies the Riccati-Volterra equation~\eqref{eq:measureFLplce}-\eqref{RicVolSqrt2}
	 on $[0,T]$,
	we set 
	\(F_\infty(\psi):=	\lim_{t\to +\infty} F(t, \psi) \) given by~\eqref{eq:FVarsigma_infty} and we let $\psi_\infty := \psi_\infty(\cdot,m)$ denote
	the unique solution of the limiting autonomous equation
	\eqref{eq:measureFLplce2}.
	Let for short \(R_\lambda(t):=\varphi(t)- \int_0^t f_{\lambda}(t-s)\varphi(s)ds\). One checks using Equation~\eqref{eq:DiffForward} together with~\eqref{eq:DoubleConvol} that
	\vspace{-.2cm}
	\begin{align*}
		&\, H_t := \int_0^T \xi_0(T-s) \,m( \mathrm{d}s) + \frac{1}{2}
		\int_0^T  \varsigma^2(s)\sigma^2(\xi_0(s))\, \psi_T(T - s)^2\, \mathrm{d}s\\
		&=  X_0\int_0^T R_\lambda(T-s) \big(m( \mathrm{d}s)  +  \frac{\kappa_1}{2} \varsigma^2(T-s)\psi_T^2(s) \dd s\big)
		+ \int_0^T\psi_T(T-s)\, \mu(s)\, \dd s
		+ \frac{\kappa_0}{2} \int_0^T \varsigma^2(s) \psi_T^2(T-s) \dd s
	\end{align*}
	As $f_{\lambda}\!\in {\cal L}^2({\rm Leb}_1)$ , $\mbox{\bf 1}_{\{0\le s \le t\}}\varsigma^2(t-s)\to \lim_{t \to +\infty} \varsigma^2(t):= \varsigma^2_\infty $ for every $s\!\in \R_+$ as $t\to +\infty$ (owing to \cite[Lemma 3.9 ]{EGnabeyeu2025}) and likewise $\lim_{t\to+\infty}(\varphi - f_{\lambda} * \varphi)(t)= \lim_{t\to+\infty} R_\lambda(t)= \ell_\infty$, we have \(F_\infty(\psi):= -\lambda \psi + \frac{\kappa_1}{2} \varsigma^2_\infty\psi^2\).
	Evaluating the equation~\eqref{eq:laplace} at $t = 0$, and using equation~\eqref{eq:DiffForward} for \(\xi_0(s)\; \forall s \leq T\), we find that~\eqref{eq: MeasureFLplaceZero}, namely
	\vspace{-.2cm}
	\begin{align}
		 &\, \mathbb{E}_{}\Big[ \exp\big( \int_0^T X_{T-s} \,m(\dd s) \big) \Big| X_0\Big]
		= \exp\left( \int_0^T \xi_0(T-s)m( \mathrm{d}s) + \frac{1}{2}
		\int_0^T \varsigma^2(T-s)\sigma^2(\xi_0(T-s))\, \psi_T(s)^2\, \mathrm{d}s \right)\label{eq:FirstIneq}\\
		&=\exp\Big( X_0\int_0^T R_\lambda(T-s) \big(m( \mathrm{d}s)  +  \frac{\kappa_1}{2} \varsigma^2(T-s)\psi_T^2(s) \dd s\big)
		+ \int_0^T\psi_T(T-s)\, \mu(s)\, \dd s\nonumber\\
		&\hspace{6cm}
		+ \frac{\kappa_0}{2} \int_0^T \varsigma^2(s) \psi_T^2(T-s) \dd s\Big)\label{eq:SecondIneq}.
		\vspace{-.2cm}
	\end{align}
	It is worth mentioning that, at this stage, substituting the expression for \(R_\lambda\) into~\eqref{eq:SecondIneq} and using once again Equation~\eqref{eq:DoubleConvol} yields the second identity in Theorem~\ref{T:VolSqrt}. 
	Furthermore, from the affinity of the squared diffusion coefficient \(\sigma^2\), $\mbox{\bf 1}_{\{0\le s \le t\}}\sigma^2(\xi_0(t-s))\to  \sigma^2(\lim_{t \to +\infty}	\xi_0(t)) $ for every $s\!\in \R_+$ as $t\to +\infty$ by the continuity of \(\sigma\)
	where the quantity \( \lim_{t \to +\infty}	\xi_0(t) = \ell_\infty X_0  + (1-a) \frac{\mu_\infty}{\lambda}\) as given by Lemma~\ref{lm:firstmoment}.
	\noindent To establish the first equality in~\eqref{eq: FT convergence}, we appeal to the first identity from Theorem~\ref{T:VolSqrt} (equivalently~\eqref{eq:FirstIneq}) and pass to the limit as \( t \to \infty \), leveraging the continuity of the map \( x \mapsto \exp(x) \). It therefore suffices to show that:
	\begin{enumerate}
		\item[(i)] $\lim_{T \to \infty}\int_0^T \mathbb{E}_{}\big[X_{T-s} | X_0\big] \,m(ds) = \left(\ell_\infty X_0 + (1-a) \frac{\mu_\infty}{\lambda}\right)m(\R_+).$
		\\
		\item[(ii)] $\lim_{T\to \infty} \int_0^T \varsigma^2(s)\, \sigma^2(\mathbb{E}_{}\big[X_s\big| X_0\big])\, \psi_T(T - s)^2\, ds 
		= \varsigma^2_\infty\sigma^2\left( \ell_\infty X_0  + (1-a) \frac{\mu_\infty}{\lambda}\right)\int_0^{\infty} \psi_\infty(s)^2 ds$.
	\end{enumerate}
		As $
	\Xi_0:=\ell_\infty X_0+(1-a)\frac{\mu_\infty}{\lambda}
	$,
	By Lemma~\ref{lm:firstmoment}~(i),
	$
	\lim_{t\to\infty}\mathbb{E}[X_t\mid X_0]=\Xi_0.
	$
	For the first term, Lemma~\ref{lm:firstmoment}~(ii), applied with
	$g(t)=\mathbb{E}[X_t\mid X_0]$,  which is bounded since $\sup_{s \geq 0}\mathbb{E}_{}\big[|X_{s}|\,| X_0\big]<\infty$.
	This
	gives
	$
	\lim_{T\to\infty}
	\int_0^t\mathbb{E}[X_{T-s}\mid X_0]\,m(ds)
	=
	\Xi_0\,m(\mathbb{R}_+).
	$
	
	\noindent Likewise for (ii) we still apply the same Lemma, but with \(g(\cdot):=\sigma^2(\mathbb{E}_{}\big[X_\cdot | X_0\big])\) which is bounded since $\sigma^2
	\big(\sup_{s \geq 0}\mathbb{E}_{}\big[|X_{s}|\,| X_0\big]\big)<\infty$ and \(g_\infty:=\sigma^2\left(  \Xi \right)\) and leveraging the fact that 
	 $\psi_\infty \in L^2(\R_+; \mathbb C)$.  Hence Lemma~\ref{lm:firstmoment}~(iii), yields
	 $
	 \lim_{T\to\infty}
	 \int_0^T
	 R_\lambda(T-s)\varsigma^2(T-s)\psi_T^2(s)\,ds
	 =
	 \ell_\infty\varsigma_\infty^2
	 \int_0^\infty\psi_\infty^2(s)\,ds
	 $ and thereby completes the proof of the first claim in~\eqref{eq: FT convergence}.\\
	
	\noindent For the second identity in~\eqref{eq: FT convergence II}, one uses Equation~\eqref{eq:SecondIneq} and pass to the limit as $t \to \infty$, namely, we show that\\
	 1.  $\lim_{T \to \infty}\int_0^T R_\lambda(T-s) \,m(ds) = \ell_\infty  m(\R_+)$;
		\; 
		2. $\lim_{T \to \infty} \int_0^t R_\lambda(T-s) \varsigma^2(T-s)\psi_T^2(s)\, ds 
		= \ell_\infty \varsigma^2_\infty \int_0^{\infty}\psi_\infty(s)^2\, ds  $;\\
		3. $\lim_{T\to \infty} \int_0^T \mu(s)\, \psi_T(T-s)\, ds  
		= \mu_\infty \int_0^{\infty}  \psi_\infty(s) \, ds $;\; 4. $\lim_{T \to \infty} \frac{\kappa_0}{2} \int_0^T \varsigma^2(s)\, \psi_T^2(T-s)\, ds  
		= \frac{\kappa_0}{2} \varsigma^2_\infty \int_0^\infty \psi_\infty^2(s)\, ds$.
	
	\smallskip\noindent
	Assertions 1., 2. and 4. follow directly from the arguments presented above still owing to Lemma~\ref{lm:firstmoment}~$(ii)$ and~$(iii)$, as Corollary~\ref{corol:riccati-extension} ensures that the desired integrability still holds for \( T = \infty \), namely \(\int_0^{\infty} \left( |\psi_\infty(t)| + |\psi\infty(t)|^2 \right) dt < \infty.\)
	Finally, for the claim 3.,
	set $w_\lambda:=f_\lambda*\mu$. Under Assumption~\ref{ass:Mean_revertF}, $w_\lambda$ is continuous, bounded, and $w_\lambda(t)\rightarrow
	w_{\lambda,\infty}
	=(1-a)\mu_\infty.$
	Convolving the Riccati--Volterra equation with $\mu$
	and using associativity yields
	\[
	\begin{aligned}
		\int_0^T\mu(s)\psi_T(T-s)\,ds
		=\frac1\lambda(w_\lambda*m)(T)+\frac{\kappa_1}{2\lambda}
		\int_0^T
		w_\lambda(T-s)\varsigma^2(T-s)
		\psi_T(s)^2\,ds.
	\end{aligned}
	\]
	Applying Lemma~\ref{lm:firstmoment}~(ii)--(iii)
	with $g=w_\lambda$, we obtain
	\[
	\begin{aligned}
		\lim_{T\to\infty}
		\int_0^T\mu(s)\psi_T(T-s)\,ds=(1-a)\mu_\infty
		\left[
		\frac{m(\R_+)}{\lambda}
		+\frac{\kappa_1\varsigma_\infty^2}{2\lambda}
		\int_0^\infty\psi_\infty(s)^2\,ds
		\right]=\mu_\infty\int_0^\infty\psi_\infty(s)\,ds,
	\end{aligned}
	\]
	where the last equality follows by integrating
	the limiting Riccati--Volterra equation and using
	$\|f_\lambda\|_{L^1(\R_+)}=1-a$.
	Combining the above limits with Equation~\eqref{eq:SecondIneq} establish the
	second convergence identity in~\eqref{eq: FT convergence II}.
	\noindent This completes the proof and we are done. \qed	
\end{document}